\documentclass[aomart]{article}

\usepackage{amsmath,amsthm,amssymb,bm}
\usepackage{mathabx}
\usepackage{authblk}
\usepackage{color}
\usepackage{tikz}
\usepackage{esint}
\usetikzlibrary{calc,decorations.markings,positioning}
\allowdisplaybreaks
\usepackage{mathtools}

\usetikzlibrary{arrows,positioning}

\theoremstyle{plain}
\newtheorem{theo}{Theorem}

\newtheorem{coro}[theo]{Corollary}
\newtheorem{proposition}[theo]{Proposition}

\newtheorem{lemma}[theo]{Lemma}
\theoremstyle{definition}
\newtheorem{define}[theo]{Definition}

\usepackage[colorlinks=true,urlcolor=blue,pdftex,pdfpagelabels=true,breaklinks=false]{hyperref}
\usepackage{xcolor}
\hypersetup{
    colorlinks,
    linkcolor={red!20!black},
    citecolor={blue!50!black},
    urlcolor={blue!80!black}
}
\usepackage{url}

\usepackage[utf8]{inputenc}
\usepackage[style=alphabetic, backend=biber,maxbibnames=99]{biblatex} 
\newbibmacro{string+doi}[1]{%
  \iffieldundef{doi}{#1}{\href{http://dx.doi.org/\thefield{doi}}{#1}}}
\DeclareFieldFormat{title}{\usebibmacro{string+doi}{\mkbibemph{#1}}}
\DeclareFieldFormat[article]{title}{\usebibmacro{string+doi}{\mkbibquote{#1}}}
\DeclareFieldFormat[incollection]{title}{\usebibmacro{string+doi}{\mkbibquote{#1}}}
\DeclareFieldFormat[book]{title}{\usebibmacro{string+doi}{\mkbibquote{#1}}}
\DeclareSourcemap{
  \maps[datatype=bibtex]{
    \map{
      \step[fieldset=issn, null]
    }
  }
}
\renewbibmacro{in:}{%
  \ifentrytype{article}{}{\printtext{\bibstring{in}\intitlepunct}}}

\DeclareMathOperator{\Hom}{Hom}
\DeclareMathOperator{\Ker}{Ker}
\DeclareMathOperator{\End}{End}
\DeclareMathAlphabet{\mathpzc}{OT1}{pzc}{m}{it}
\newcommand{\E}{\mathcal{E}}
\newcommand{\e}{\mathbf{e}}

\newcommand{\Fe}{\mathpzc{e}}
\newcommand{\FE}{\mathpzc{E}}
\newcommand{\Feu}{\underline{\mathpzc{e}}}
\newcommand{\calS}{\mathcal{S}}
\newcommand{\TN}{\mathrm{TN}}
\newcommand{\tslash}{{\slash\!\!\!t}}
\newcommand{\ii}{\mathrm{i}}
\newcommand{\p}{\partial}

\newcommand{\D}{\mathpzc{D}}

\newcommand{\IC}{\mathbb{C}}
\newcommand{\IR}{\mathbb{R}}

\title{Instantons on multi-Taub-NUT Spaces III:\\
Down Transform, Completeness, and Isometry}

\author[*]{Sergey A. Cherkis}
\author[$\dagger$]{Andr\'es Larra\'in-Hubach}
\author[$\star$]{Mark Stern}
\affil[*]{\small Department of Mathematics, University of Arizona, Tucson, AZ 85721-0089}
\affil[ ]{\small cherkis@math.arizona.edu}
\affil[$\dagger$]{\small Department of Mathematics, University of Dayton, Dayton, OH 45469}
\affil[ ]{\small  alarrainhubach1@udayton.edu}
\affil[$\star$]{\small Department of Mathematics, Duke University, Durham, NC 27708-0320}
\affil[ ]{\small stern@math.duke.edu}

\begin{document}

\begin{titlepage}
\renewcommand{\thepage}{i}
\date{}  

\maketitle
\abstract{
The index bundle of a family of Dirac operators associated to an instanton on a multi-Taub-NUT space forms a bow representation.  We prove that the gauge equivalence classes of solutions of this bow representation are in one-to-one correspondence with the instantons.  

We also prove that this correspondence establishes an isometry of the bow and instanton moduli spaces.
}
\thispagestyle{empty}

\end{titlepage} 

\tableofcontents

\section{Introduction}
Atiyah, Drinfeld, Hitchin, Manin \cite{Atiyah:1978ri} and Nahm \cite{NahmADHM, Nahm1979-2} discovered a nonlinear generalization of the Fourier transform relating instantons and monopoles to, respectively, algebraic data and solutions of an ordinary differential equation.  Kronheimer and Nakajima \cite{KN} generalized this transform to instantons on Asymptotically Locally Euclidean (ALE) spaces. The result of the transform in this case is algebraic data associated to a quiver.  Instantons on Asymptotically Locally Flat (ALF) spaces, on the other hand, are associated to a bow \cite{Cherkis:2010bn}.  The focus of this paper is exactly such a  transformation mapping an instanton on the prototypical ALF space, the multi-Taub-NUT,  $\TN_k$, to (a gauge equivalence class of) a bow solution.

In each of these cases there are two hyperk\"ahler moduli spaces: the moduli space of instantons (or monopoles) and the moduli space of their related quiver or bow data. 
The transform not only provides a construction of all relevant instantons and monopoles but also defines an isometry between the two moduli spaces.  The analytic details of the proof of this fact  differ significantly case by case. 
For monopoles, the completeness of the Nahm construction was proved in \cite{Hitchin:1983ay} and  isometry in \cite{NakajimaMon}, while for instantons on ALE manifolds both completeness and isometry are proved in \cite{KN}.
There is also a relation between instantons on a four-torus $\mathbb{R}^4/\Lambda$ and its dual four-torus $\mathbb{R}^4/\Lambda^*.$  The isometry in this case is proved in \cite{Braam:1988qk} (see also \cite{DK}).  
{Doubly periodic instantons of rank 2 with quadratic curvature decay were related to Hitchin systems on a torus bijectively in \cite{jardim99,Jardim01,Jardim02}. 
These results were refined (by removing some assumptions) and the isometry of the corresponding moduli spaces proved in \cite{BiquardJardim}. 
They were generalized further in \cite{Mochizuki} to arbitrary rank and to \(L^2\) instanton curvature without any additional decay assumptions.\footnote{ 
Another case of isometry proof is the Nahm transform between two Hitchin systems \cite{Szabo}.  This, however, does not involve instantons or monopoles.
}
}

We have already established many of the key analytic results, including the index theorem, the curvature decay rate, and the asymptotic form of the instanton, in \cite{First}. We also  proved in \cite{Second} that the connection resulting from the Up transform \cite{Cherkis:2009jm,Cherkis:2010bn} of any bow solution is, indeed, an instanton.  Here, we complete the circle by formulating the Down transform 
{(Def.~\ref{DefDown},p.~\pageref{DefDown})}, 
in which the bow {representation} emerges on an index bundle of a family of Dirac operators associated to an instanton on the multi-Taub-NUT space.   We prove that the two transforms, Up and Down, are inverse of each other ({Thm.~\ref{template},p.~\pageref{template}}).  
We also prove that each acts as isometry between the moduli space of a bow {representation} and the moduli space of instantons ({Thm.~\ref{IsomThm},p.~\pageref{IsomThm}}).

Our approach is analytic, descending from \cite{Corrigan:1983sv,Corrigan:1978ce}. For an algebro-geometric approach to the  bow construction see \cite{Cherkis:2017pop,Cherkis:2020kpr}.
Both the instanton on  $\TN_k$  and the bow solution come equipped with a corresponding natural Dirac type operator (a partial differential operator in case of an instanton and an ordinary differential operator in the case of the bow). The Up transform \cite{Second} constructs the vector bundle over $\TN_k$ (together with its induced instanton connection) from the kernel of the bow Dirac operator. The Down transform constructs the bundle over the bow 
(with an induced bow solution) 
from the kernel of the instanton  Dirac operator ({Def.~\ref{DefDown},p.~\pageref{DefDown}}).  
As described in \cite{Cherkis:2010bn,Second}, in this way every bow solution can be mapped to an instanton via the Up transform. Here, after setting our conventions in Section~\ref{sec:Setup}, we formulate the Down transform  in Section~\ref{sec:DownNahmBif},  associating a bow solution to each instanton.  We prove that the two transforms, Up and Down, are inverse of each other in Section~\ref{sec:CompUni}.  After defining the moduli space of instantons ({Def.~\ref{ModSp} and Thm.~\ref{23andme},p.~\pageref{23andme}}) in Section~\ref{sec:Moduli}, we prove in Section~\ref{sec:Isometry} that each transform acts as isometry between the moduli space of  bow solutions and the moduli space of instantons ({Thm.~\ref{IsomThm},p.~\pageref{IsomThm}}). 

\section*{Acknowledgements}
SCh thanks the Berkeley Center for Theoretical Physics for hospitality during the final stages of this work. 
The work of MS is supported by the Simons Foundation Grant 3553857. 
We thank the referees for their helpful suggestions. 
\section{Setup}
\label{sec:Setup}
\subsection{The multi-Taub-NUT Space}
For the underlying base manifold we choose the prototypical  hyperk\"ahler ALF manifold: the multi-Taub-NUT space, $\TN_k^\nu$ (or simply $\TN_k$ with $ \nu$ implied). It is a circle fibration 
\begin{align}\label{fibr}S^1\rightarrow\TN_k^\nu\xrightarrow{\pi_k}\mathbb{R}^3
\end{align} over a Euclidean  three-space $\mathbb{R}^3\simeq\mathfrak{Im}\,\mathbb{H}$ with a triholomorphic isometry, generated by the vector field $\partial_\tau,$ rotating the fiber.  We denote the fixed points\footnote{
Since the orbit of a fixed point $\nu_\sigma$ consists of that single point, we denote both the fixed point $\nu_\sigma$ in $\TN_k$ and the corresponding point $\pi_k(\nu_\sigma)$ of the base $\mathbb{R}^3$ by $\nu_\sigma,$ and view $\nu_\sigma$ as an imaginary quaternion: $\nu_\sigma\in\mathfrak{Im}\, \mathbb{H}.$
} 
of this action by $\{\nu_\sigma\}_{\sigma=1}^k.$  On the complement of these $k$ points the space is a circle bundle.  Let $\hat{\omega}$ denote the connection one-form of this circle bundle.   
If we let $t\in\mathfrak{Im}\,\mathbb{H}$ and let $t_\sigma=|t-\nu_\sigma|,$ 
then the $\TN_k^\nu$ metric  has a Gibbons-Hawking form 
$$V |dt|^2+\frac{\hat{\omega}^2}{V},$$
with $V=l+\sum_{\sigma=1}^k  \frac{1}{2 t_\sigma}.$ In coordinates over a contractible open subset of $\mathbb{R}^3\setminus\{\nu_\sigma\}_\sigma,$
with $\tau\sim\tau+2\pi$ a periodic coordinate along the $S^1$ fiber, the one-form $\hat{\omega}$ is 
 $\hat{\omega}=d\tau+\pi_k^*(\omega),$ where the one-form $\omega$ on $\mathbb{R}^3$ has $d\omega$ Hodge dual to $dV$: 
\begin{align}\label{domega}d\omega=*_3 dV.
\end{align}  We choose the orientation with volume form $\mathrm{Vol}=Vdt^1\wedge dt^2\wedge dt^3\wedge d\tau$. Next, we consider on this space   a Hermitian bundle $\E$ {of rank n}, with an instanton connection $A$.  As  in \cite{First}, an {\em instanton} is a connection   with  anti-self-dual,  square integrable curvature $F_A$.  
{Throughout this paper we restrict ourselves to the generic case; namely, we assume that the instanton has {\em generic asymptotic holonomy} \cite[Sec.~4.1]{First}.  
}
{
This implies that the eigenvalues of the holonomy around the circle fiber \(S^1_t\) of the Taub-NUT space over \(t\in\mathbb{R}^3\), have distinct limits $\{\exp(2\pi\ii \lambda_j/l)\}_{j=1}^{n}$ as \(|t|\to\infty\).
 } We order $0\leq\lambda_1< \lambda_2<\ldots<\lambda_n<l.$

To discuss the Dirac operator, it is convenient to  introduce an associated orthonormal frame,
$\Theta_j=\frac{1}{\sqrt{V}}(\partial_j-\omega_j\partial_\tau), \Theta_4=\sqrt{V}\partial_\tau,$ and its dual coframe
$\theta^j=\sqrt{V}dt^j, \theta^4=\frac{d\tau+\omega}{\sqrt{V}}$.
  In this frame, the three symplectic forms are $w^i=\frac{1}{2}\theta^i\wedge\theta^4+\frac{1}{4}
  {\sum_{jk}\epsilon_{ijk}
  }
  \theta^j\wedge\theta^k.$ These are self-dual.  The corresponding Clifford algebra bundle is generated by $c^q:=Cl(\theta^q)$, satisfying $\{c^q,c^p\}=-2\delta^{pq}.$  The chirality operator $\gamma:=- Cl(\mathrm{Vol})=-   c^1c^2c^3c^4$ satisfies $\gamma\gamma=1$ and splits the spin bundle into positive and negative chirality eigenbundles $S^+\oplus S^-\!\!\rightarrow\TN_k$.  The sign in our definition of $\gamma$ is chosen so that its Clifford action is compatible with the Hodge  star action on two-forms: $\gamma Cl(\eta)= Cl(*\eta)$, for any two-form $\eta.$  In particular, Clifford multiplication by any self-dual $2$-form annihilates $S^-$. 
Since the Riemann curvature\footnote{
The $\TN_k$ space is hyperk\"ahler, which implies the anti-self-duality of its Riemann curvature.
}  
is also anti-self-dual, the bundle $S^+$ is trivial. 
The Clifford action of the symplectic forms  annihilates $S^-$ and generates a two-dimensional representation of the quaternions on (the space of covariantly constant sections of) $S^+$, with quaternionic units $I_j:=Cl(w^j)=\frac{1+\gamma}{2}c^jc^4$.  Importantly, since the base space is hyperk\"ahler, these complex structures are covariantly constant.  To complete our basis of quaternionic units, we set $I_4:=\frac{1+\gamma}{2}$, which acts as the identity on $S^+.$ Hence  $I_a=-\frac{1+\gamma}{2}c^4c^a$ for $a=1,2,3,4.$

\subsection{Bows}\label{subBows}
\begin{figure}[htb!]
\centering
\resizebox{0.7\textwidth}{!}{%
%
%
%
%
%
%
%
%
%
\begin{tikzpicture}[>=latex']
    \draw (0,0) ellipse (4 and 2);
           \node(p1) at ($(-155:4 and 2)$) [label=below:$p_1$]{\tikz\draw[red,fill=white] (0,0) circle (.5ex);};
           \node(p2) at ($(-55:4 and 2)$) [label=below:$p_2$] {\tikz\draw[red,fill=white] (0,0) circle (.5ex);};
           \node(p3) at ($(75:4 and 2)$)  [label=below:$p_3$] {\tikz\draw[red,fill=white] (0,0) circle (.5ex);};
    \draw[blue,thick] (p1)+(0,5)  arc (-155:-55:4 and 2);
     \draw[blue,thick] (p2)+(0,6) arc (-55:75:4 and 2);
     \draw[blue,thick] (p3)+(0,7) arc (75:205:4 and 2);
\draw[dotted] (p1) -- ++(0,7);
\draw[dotted] (p2) -- ++(0,6);
\draw[dotted] (p3) -- ++(0,7);
 \node(Int1) at ($(0,7)+(150:4 and 2)$) [label=above:$J_1$] {};
 \node(Int2) at ($(0,5)+(-100:4 and 2)$) [label=above:$J_2$] {};
 \node(Int3) at ($(0,6)+(15:4 and 2)$) [label=right:$J_3$] {};
\draw[->,thick] (p1)+(0,7) node {\tiny$\bullet$} node[anchor=south west] {$p_{1}-$} -- ++(0,6) node[left] {
} -- ++(0,-1) node {\tiny$\bullet$} node[anchor=south west] {$p_{1}+$};
\draw[->,thick] (p2)+(0,5) node {\tiny$\bullet$} node[anchor= north west] {$p_{2}-$} -- ++(0,5.5) node[right]{
} -- ++ (0,0.5) node {\tiny$\bullet$} node[inner sep=5pt, anchor=south] {$p_{2}+$};
\draw[->,thick] (p3)+(0,6) node {\tiny$\bullet$} node[anchor=north] {$p_{3}-$} -- ++(0,6.5) node[left]{
} -- ++(0,0.5) node {\tiny$\bullet$} node[anchor=south] {$p_{3}+$};
    \draw[opacity=0] ($(-130:6 and 4)$)  arc (-130:-85:6 and 4);
\end{tikzpicture}%
}
\caption{The $A_k$ bow (presented here with \(k=3\)) is obtained by cutting the circle at points $p_\sigma$ and connecting the resulting intervals $J_\sigma=[p_{\sigma-1}+,p_\sigma-]$ by edges.  The small representation of this bow has Hermitian line bundle $\Feu$ of rank one everywhere. The moduli space at level $i\nu$ of this representation is the multi-Taub-NUT space $\TN_k^\nu$ \cite{Cherkis:2008ip}.}
\label{fig:TNkBow}
\end{figure}
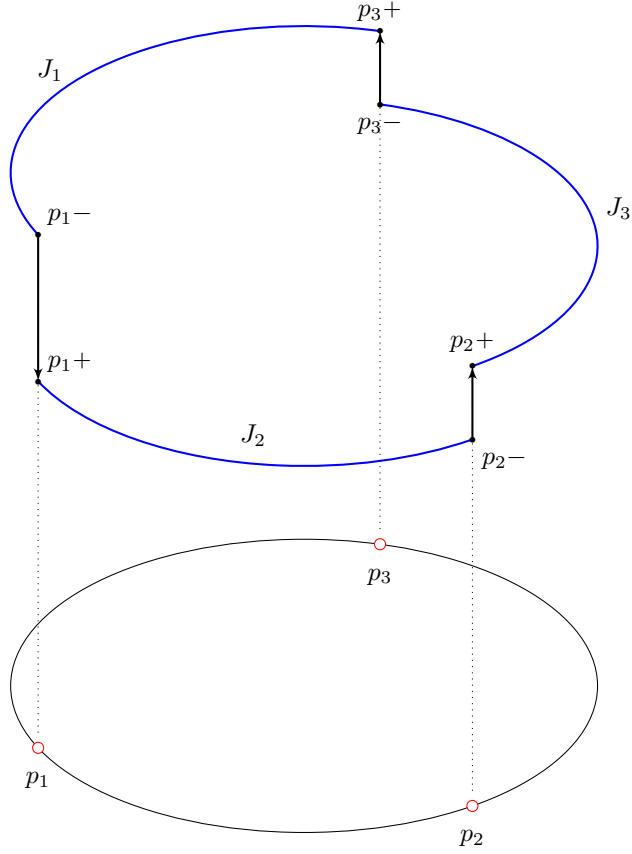

We give an abbreviated discussion of bows here. For more detail see \cite{Second}. An $A_{k-1}$ {\em bow}  consists of $k$ intervals $J_\sigma = [p_{\sigma-1}+,p_{\sigma}-]$ of length $l_\sigma$ and an edge from $p_\sigma-$ to $p_\sigma+$ for each $\sigma=1,\ldots,k$.  
(Note, that if one identifies $p_{\sigma}-$ with $p_{\sigma}+$, then one obtains a circle of length $l:= \sum_{\sigma=1}^kl_\sigma$, with $k$ marked points $P:=\{p_\sigma\}_{\sigma=1}^k$, as in Fig.~\ref{fig:TNkBow}. 
Let $s\sim s+l$ be the coordinate along this circle.)  
{Let $\calS$ be a trivial (with product connection) $\IC^2$ bundle over $\mathcal{J}:=\sqcup_\sigma J_\sigma$, carrying a representation $\e_j, j=1,2,3$ of quaternionic units.} (As we explain in Sec.~\ref{subBundles}, $\e_j$ action  will be contragredient to that of $I_j$.) 
A {\em representation} of this bow {\cite[Def.~4]{Second}
}
 consists of a set of points  $\Lambda =\{\lambda_i\}_{i=1}^n\subset \cup_\sigma J_\sigma$, and a collection  of  Hermitian vector bundles, collectively denoted $\mathpzc{E}$, with a bundle on {(the closure of)} each of the maximal subintervals of $\cup_\sigma J_\sigma$ that contain no $\Lambda$ or $P$ point in its interior.  
For simplicity we assume from now on that $\Lambda$ is disjoint from $P$.  
  Let $\Lambda^0\subset  \Lambda$ be the subset consisting of all $\lambda$ for which the Hermitian bundles on the two subintervals 
containing $\lambda$ as an endpoint have the same rank. One last piece of data in a bow representation is a rank one Hermitian vector space $W_\lambda$, for each such $\lambda\in \Lambda^0$. 

We associate {\em bow data} {\((T,B,Q)\)}  to each bow representation 
{\cite[Sec.~3.2]{First}}. 
This consists of 
\begin{itemize}
\item[(i)] Nahm data comprised of (a) a connection $\nabla_{\frac{d}{ds}}$ on each bundle together with  (b) three Hermitian sections (called Nahm matrices) $\{T^j\}_{j=1}^3$ of  $\End(\mathpzc{E})$,
\item[(ii)] bifundamental data comprised of an element $B_\sigma\in \Hom (\mathpzc{E}_{p_\sigma+},\calS\otimes\mathpzc{E}_{p_\sigma-}),$ for each  $\sigma$, and
\item[(iii)] fundamental data consisting of an element $Q_\lambda\in \Hom(W_\lambda,\calS\otimes \mathpzc{E}_\lambda), \forall\lambda\in \Lambda^0$. 
\end{itemize}
The bow data form a hyperk\"ahler affine space with the metric induced by the norm
{
\begin{align}\label{AffNorm}
	\left\|(\dot{T},\dot{B},\dot{Q})\right\|^2 = \sum_\sigma |\dot{B}_\sigma|^2 + \sum_{\lambda\in\Lambda^0} |\dot{Q}_\lambda|^2 
	+ \sum_\sigma \int_{J_\sigma} \left(|\dot{\nabla}_{\frac{d}{ds}}|^2 + \sum_{j=1}^3 \dot{T}^j\right) ds
.\end{align}
The linear spaces of fundamental and bifundamental data inherit a quaternionic action from \(\calS\).  
The tangent space to the affine space of the Nahm data acquires a quaternionic action as follows:  
consider \(\dot{T}^0:=\dot{\nabla}_{\frac{d}{ds}}\) and \(\dot{T}^1,\dot{T}^2,\dot{T}^3\) as four components of a quaternion, \(\dot{T}^0\otimes\e_0 + \sum_{j=1}^3\dot{T}^j\otimes \e_j\),  
with the left quaternionic action.  (We put, \(\e_0=1\).)
}

There is a natural {isometric} gauge group action on the bow data, with an associated hyperk\"ahler moment map. The resulting moment map equations (see \cite[Sec.~3.3]{Second}) for moment map level $i\nu$, with 
\begin{align}\label{nanosec}
	\nu =   \sum_\sigma \nu_\sigma (\delta_{p_\sigma-}-\delta_{p_\sigma+}):= \sum_\sigma {\e_j}\nu_\sigma^j(\delta_{p_\sigma-}-\delta_{p_\sigma+})
\end{align}
are as follows:\\ 
 the Nahm equations within each subinterval:
\begin{align}\label{nahmeq0}
[\ii\nabla_{\frac{d}{ds}},T_1]=[T_2,T_3], \text{ and cyclic permutations,}
\end{align} 
 the boundary conditions at $\lambda\in\Lambda_0:$
\begin{align}\label{Tdisc}
	{\e_j} T^j(\lambda+)-{\e_j} T^j(\lambda-)=\mathfrak{Im}\,\ii Q_\lambda Q_\lambda^\dagger,
\end{align}
the usual Nahm pole boundary condition \cite[Sec.3.2, Eqs.(5,6)]{Second} at $\lambda\notin\Lambda_0$, and
  the boundary conditions at the ends of the intervals
\begin{align}\label{pEnd}
	{\e_j} T^j(p_\sigma-)-\nu_\sigma&=-\mathfrak{Im}\, \ii B_\sigma B_\sigma^\dagger,&
	{\e_j} T^j(p_\sigma+)-\nu_\sigma&=\mathfrak{Im}\, \ii B_\sigma^c \left(B_\sigma^c\right)^\dagger,&
\end{align} 
where $B^c_\sigma$ denotes the charge conjugate of $B_\sigma$ (see Appendix~\ref{ApA}). 
Here {$ \mathfrak{Im}\, X:= X - \frac{1}{4}\sum_{a=0}^3\e_a X \e_a^\dagger$} is the quaternion imaginary part of $X$, c.f. Eq.~\eqref{fierz}. 

Any set of data satisfying these moment map equations is called a {\em bow solution}. 
{Thus the space of bow solutions is the level set of the moment map; it inherits a  metric from the ambient affine space. The 
{\em moduli space of the bow representation} is the quotient of this level set by the gauge group.  
This is the hyperkahler quotient of the space of bow data by the group of  gauge transformations. We refer to \cite[Sec.~3,p.446]{Second} for the detailed definition of this space. 
}
\subsection{Manifolds and Bundles from Bows}\label{subBundles}
  
A key part of our story is that the multi-Taub-NUT space $\TN_k^\nu$ is itself  the moduli space of a {\em small  $A_{k-1}$  bow representation}.   See \cite{Cherkis:2010bn} and \cite[Sec.~3.5]{Second} for a review.  We summarize the relevant details.  The small representation has   $\Lambda=\emptyset$,  and involves only a line bundle $\Feu\to \cup_\sigma J_\sigma$.  The level set $\bm{\mu}^{-1}(\ii \nu )$ of the hyperk\"ahler reduction of the data at level $  i\nu$ consists of a connection along the bow, the (commuting) Nahm endomorphisms $\{t^j\}_{j=1}^3$, and 
{$b_\sigma\in \Hom(\Feu_{p_\sigma+},\calS\otimes \Feu_{p_\sigma-})$} 
satisfying  
\begin{align}\label{bbnd}
	b_\sigma b_\sigma^{\dagger} = |t-\nu_\sigma|+\ii\sum_{j=1}^3{\e_j}(t^j-\nu^j_\sigma).
\end{align}
 The gauge quotient  of $\bm{\mu}^{-1}(\ii\nu)$ with its canonical hyperk\"ahler structure is $\TN_k^\nu$.  Consequently, $\bm{\mu}^{-1}(\ii\nu)$ is a principal bundle over $\TN_k^\nu$.  

For any point $s$ on the bow, consider the fiber $\Feu_s$ over that point and the trivial bundle $\Feu_s\times \bm{\mu}^{-1}(\ii\nu) \to \bm{\mu}^{-1}(\ii\nu)$.  Since the gauge group acts on the level set and on this fiber, we obtain an associated line bundle 
$\Fe_s\rightarrow\mathrm{TN}_k.$   
Doing this for each value of $s$, we obtain a bow-parameterized family of 
 {\em  tautological line bundles} \cite{Second}.  
 {These bundles may change topologically  from one interval to another, but not within the interval. That is, setting $K_\sigma := \Fe_{p_\sigma+}\otimes \Fe_{p_\sigma-}^{-1}$, we have $\Fe_s=\Fe_0\otimes \otimes_{p_\sigma<s} K_\sigma$. (See \cite[Sec.~3.5]{Second}.})   By this definition, a section of $\Fe_s$ over $\TN_k^\nu$ is an equivariant section of $\Feu_s\times \bm{\mu}^{-1}(\ii\nu) $ over the level set.  We will use this relation extensively in what follows.

 {We now relate the Taub-NUT side representation of quaternions $I_j$ acting on $S^+$ with the bow side representation of quaternions $\e_j$ acting on $\calS$.} 
 Note that the gauge group  acts trivially  on the trivial $\IC^2$ bundle $\calS$ on the bow. 
 {Hence, in the above hyperk\"ahler quotient construction}, $\calS$ descends to a trivial $\IC^2$ bundle on $\TN_k^\nu$.  
 It is natural to identify this bundle with $(S^+)^*$ --- the dual of the (trivial) bundle of positive chirality spinors, equipped with the contragredient action of the unit quaternions: for $z\in (S^+)^*$, {$\e_a z:= z\circ I_a^\dagger.$}  

\subsection{Dirac Operators and Tautological  Connections }
Let $(\E,A)$ be a Hermitian bundle on $\TN_k$ equipped with an instanton connection. Each of the bundles $\Fe_s$ carries an abelian instanton connection $d+\ii a_s $ computed in \cite{Cherkis:2009jm}, which we describe momentarily.  This allows us to associate to  $\E$  a family of bundles $\E\otimes \Fe_s,$ each equipped with  the induced  instanton connection $\nabla^s$ and an associated Dirac operator $D_s = \left(\begin{array}{cc}0&D_s^-\\D_s^+&0\end{array}\right)$, where $D_s^+:\Gamma(S^+\otimes \E\otimes \Fe_s)\to \Gamma(S^-\otimes \E\otimes \Fe_s)$ and 
$D_s^- = D_s^{+\dagger}$ is the formal adjoint of $D_s^+$. We will write $\nabla^{(s)}_X$ for the $\nabla^s$ covariant derivative in the $X$ direction when confusion may arise. (Note, however, 
that, as the explicit form of the connection \eqref{Eq:connfam} will make clear,  
for any given $\sigma$ and any horizontal vector $X\perp \Theta_4$, the covariant directional derivative is $s$-independent: $\nabla_X^{s_1}= \nabla_X^{s_2}$ for all $s_1, s_2\in\mathring{J}_\sigma$. Hence we will often omit the superscript in this case.)  Anti-self-duality of the connection is equivalent to  $D_s$ satisfying  
\begin{align}
\label{eq:LapDir}
D_s^{+\dagger}D_s^+ &=\nabla^{s *}\nabla^s
   =-V(\nabla_\tau^s)^2-\frac{1}{V}(\nabla_j-\omega_j \nabla_\tau)(\nabla_j-\omega_j\nabla_\tau),
\end{align} 
where $\nabla_j:= \nabla_{\frac{\p}{\p t^j}}^{(s)}$ and $\nabla_\tau^{s}:= \nabla_{\frac{\p}{\p \tau}}^{s}$. 
{Therefore, \(\langle \psi, D_s^{+\dagger}D_s^+ \psi\rangle = \| \nabla^s \psi \|^2,\) 
and positive chirality elements of $\mathpzc{E}_s:=\mathrm{Ker}_{L^2}\, (D_s):=  \mathrm{Ker} \, (D_s)\cap L^2$ are covariant constant and therefore vanish. 
This immediately implies 
\begin{lemma}
\label{trivKer}
For any instanton connection  
\(\mathrm{Ker}_{L^2} D_s^+ =0\).
\end{lemma}
} 
In particular, $\mathpzc{E}_s\subset \Gamma(S^-\otimes \E\otimes \Fe_s).$ 

Let us turn to the tautological connection on $\Fe_s$. As a subset of an affine hyperk\"ahler space, the level set $\bm{\mu}^{-1}(\ii\nu)$   inherits the induced  metric.  Thus, as a principal bundle over $\TN_k$, it has a natural connection, which, in turn is inherited by the $\Fe_s$ family of line bundles.   
In suitable local coordinates, it  can be written as $d+\ii a_s $ with one-form 
\begin{align}\label{Eq:connfam}
a_s=sa^{(0)}+\sum_{\sigma | p_\sigma<s} a^{(\sigma)},
\end{align}
 with 
$
a^{(0)}=\frac{d\tau+\omega}{V}=\frac{\theta^4}{\sqrt{V}}
$
 the connection one-form on $\Fe_0$ and
$
a^{(\sigma)}=\frac{1}{2 t_\sigma}\frac{d\tau+\omega}{V}-\eta_\sigma
$
the connection one-form on the line bundle $K_\sigma.$ 
Here  the $\eta_\sigma$ are one-forms on a local coordinate patch in $\mathbb{R}^3\setminus\{\nu_\sigma\}_\sigma$, satisfying $d\eta_\sigma=*d\frac{1}{2t_\sigma}.$  
In fact, $a^{(\sigma)}$ and any multiple of $a^{(0)}$ are all connection one-forms for abelian instantons on $\TN_k^\nu$. The one-form $a^{(0)}$ is globally defined, as it is a connection form on a trivial line bundle.  These connections were found by Ruback in \cite{Ruback}, and their curvature forms form a basis of $L^2$ cohomology of $\TN_k^\nu$, as proved in \cite[Sec.7.1.2]{HHM}. 
Since $a^{(0)}$ is globally defined, we can write  $D_s=D_0+\ii\mathrm{Cl}(sa^{(0)}),$ where $D_0$ is locally constant in $s$ for $s\neq p_\sigma. $ Hence in the interior of each bow interval,  
\begin{align}\label{Eq:Comm}
[D_s,\ii\partial_s]&=\frac{c^4}{\sqrt{V}},&
[D_s,t^i]&=\frac{c^i}{\sqrt{V}},\\
\label{Eq:Comm2}[\nabla^{s*}\nabla^s , \ii\partial_s]&=-\frac{2}{\sqrt{V}}\nabla^{(s)}_{\Theta_4},&
[\nabla^{s*}\nabla^s , t^i]&=-\frac{2}{\sqrt{V}}\nabla^{(s)}_{\Theta_i}.
\end{align}

The $b_\sigma$  bow data are  sections  of $(S^+)^\ast\otimes K_\sigma^{-1} $. Let $b_\sigma^c$ denote the charge conjugate of $b_\sigma$. (See \cite[Sec.~3.1]{Second}). These sections are $O(|t|^{\frac{1}{2}})$ by \eqref{bbnd} 
{(in fact \(|b_\sigma|^2 = |t-\nu_\sigma|\))
 }
 and satisfy 
$
(d+\ii\eta_\sigma)b_\sigma=\ii \e_j\frac{dt^j}{2t_\sigma} b_\sigma
$ and 
$
(d-\ii\eta_\sigma)b_\sigma^c=-\ii e_j\frac{dt^j}{2t_\sigma} b_\sigma^c
$ (see Appendix~\ref{sec:CovDer}),
which implies
\begin{align}\label{Eq:bCov}
\nabla_{\Theta_a}b_\sigma&=-\frac{\ii}{\sqrt{V}} \e_a^\dagger \frac{b_\sigma}{2t_\sigma},&
\nabla_{\Theta_a}b_\sigma^c&=\frac{\ii}{\sqrt{V}} \e_a^\dagger \frac{b_\sigma^c}{2t_\sigma}.
\end{align}
This means, in turn,  that $b_\sigma$ and $b_\sigma^c$ are harmonic:
 \begin{align}\nabla^{*}\nabla b_\sigma=0 = \nabla^*\nabla b_\sigma^c.
\end{align} 
In fact, they also satisfy a Dirac like equation:
\begin{align}
{\nabla_{\Theta_a}\frac{b_\sigma}{t_\sigma} I_a=}
	\e_a^\dagger  \nabla_{\Theta_a}\frac{b_\sigma}{t_\sigma} &= 0, & 
 \nabla_{\Theta_a}\frac{b_\sigma^c}{t_\sigma} I_a=
 \e_a^\dagger  \nabla_{\Theta_a}\frac{b_\sigma^c}{t_\sigma}
=0
.\end{align}

\section{Down Transform}
\label{sec:DownNahmBif}
The Down Transform assigns to an instanton $(\mathcal{E}, A)$  (and any chosen small bow representation, whose moduli space is the underlying  $\TN_k$) a large bow representation and its solution. First we introduce the large bow representation associated to $(\mathcal{E}, A)$.  Then we focus on the corresponding bow solution.  
{We recall that we denote bow data (reviewed in Sec.~\ref{subBows}) by \((T,B,Q)\), and a {\em bow solution} is  bow data that solves the bow moment map equations (\ref{nahmeq0},\ref{Tdisc},\ref{pEnd}). 
}

As proved in \cite[Thm.~B]{First} (under the assumption of  generic asymptotic holonomy),   any instanton connection   has the form
\begin{align}\label{connection}
A= \, \mathop{\oplus}_{\lambda\in\Lambda} \left(\pi_k^* \eta_\lambda- \ii (\lambda+\frac{m_\lambda}{2 |t|})\frac{\varpi }{V}\right)+O(|t|^{-2}),
\end{align}
where $\Lambda$ is a subset of $[0,\ell)$ with $|\Lambda|=\mathrm{rank}(\mathcal{E})$, 
each $m_\lambda$ is an integer, and $\eta_\lambda$ a connection  
defined on  degree $m_\lambda$ Hopf line bundles over $\mathbb{R}^3\setminus K$, for some compact set $K$. The bundle $\mathcal{E}$ decomposes, near infinity, as a direct sum of line bundles pulled back from $\mathbb{R}^3\setminus K$. For details see \cite[Sec.~6]{Second}. 

Let $G_{D_s}$ denote the Green's operator for $D_s^-D_s^+$. 
Note that, thanks to \eqref{eq:LapDir}, this Green's operator commutes with quaternionic units.  
{We define the bow representation to consist of $\Lambda$ (the set of asymptotic eigenvalues of \(\ii\ell A(\frac{\partial}{\partial \tau})\)  read off from \eqref{connection})
 } and the   vector bundle $ \mathpzc{E} $ over the bow to be  the   index bundle, with fiber $ \mathpzc{E}_s:=\Ker_{L^2}(D_s) $ at $s$. 
 By Lemma~\ref{trivKer}, the locally constant rank function is $R(s):= \mathrm{Index}(D_s^-) = \dim(\Ker_{L^2}(D_s^-))= \dim(\Ker_{L^2}(D_s))$. 

Let $\Lambda^0:=\{\lambda\in\Lambda\, |\, m_\lambda = 0\}.$ 
Recall the large bow data {(T,B,Q)} comprises 
\begin{enumerate}
\item[(i)] A unitary connection on $ \mathpzc{E} $,
\item[(ii)] (Hermitian) Nahm matrices $\{T^j(s)\}_{j=1}^3,$ 
\item[(iii)] bifundamental data $B_\sigma\in \Hom(\mathpzc{E}_{p_\sigma+},{\calS}\otimes \mathpzc{E}_{p_\sigma-})$, for $\calS$ a  complex 2 dimensional irreducible representation of the quaternions, 
\item[(iv)] fundamental data $Q_\lambda\in \Hom(\mathpzc{E}_\lambda, {\calS}\otimes \IC),$ for each $\lambda\in \Lambda^0$. 
\end{enumerate}
{
\begin{define}
\label{DefDown}
The {\em Down transform} of an instanton $(\mathcal{E}, A)$ is the bow solution \((T,B,Q)\) on the bow representation \((\Lambda, \mathpzc{E}_s:=\Ker_{L^2}(D_s), W_\lambda:=\Ker_{L^\infty} \nabla_\lambda^*\nabla_\lambda)\) with 
\begin{itemize}
	\item[]
\(T=(\nabla_{\frac{d}{ds}}, T^1,T^2,T^3)\) defined by Eqs.~\eqref{DownNab} and \eqref{DownT}, 
    \item[]
\(B=\{B_\sigma\}_\sigma\) defined by Eq.~\eqref{DownB}, and 
    \item[]
\(Q=\{Q_\lambda\}_{\lambda\in\Lambda^0}\)  defined by Eq.~\eqref{disQ1}.
\end{itemize}
\end{define}
}

\subsection{Nahm and Bifundamental Data}\label{nd}
Let $\Pi_{D_s}$ denote the $L^2$ orthogonal projection onto  $\Ker_{L^2}(D_s)=\Ker_{L^2}(D_s^-) $. We assign as connection for $\mathpzc{E}$, the  canonical connection on the index bundle:  
\begin{align}\label{DownNab}
	\nabla_{\frac{d}{d s}} := \Pi_{D_s}\frac{d}{d s}\Pi_{D_s}.
\end{align}
For $s\not\in \Lambda$, we define the Nahm matrices
\begin{align}\label{DownT}
	T^j(s):=   \Pi_{D_s}t^j\Pi_{D_s}.
\end{align}
Here $t^j$ is shorthand for multiplication by the $t^j$ coordinate.  
The bifundamental $B_\sigma  = \left(\begin{array}{c}B^1_\sigma\\B^2_\sigma \end{array} \right)$ is similarly defined by 
\begin{align}\label{DownB}
	B_\sigma^\alpha:=   \Pi_{D_{p_\sigma-}}b_\sigma^\alpha\Pi_{D_{p_\sigma+}}.
\end{align}
By \cite[Prop.~25]{First}, for $s\not\in \Lambda$, elements of $\Ker_{L^2}(D_s^-) $ are exponentially decaying. 
Hence $T^j$ and $B^j$ are well defined bounded operators.  
We next verify that they satisfy the moment map equations. 
We record {four} useful commutators which are immediate consequences of \eqref{Eq:Comm} and \eqref{Eq:Comm2}. For $s\not\in\Lambda$, 
\begin{align}\label{scomG}
\left[\frac{d}{ds},G_{D_s}\right] &=  -G_{D_s}\left(\frac{\ii c^4}{\sqrt{V}}D_s^++D_s^-\frac{\ii c^4}{\sqrt{V}}\right)G_{D_s}
=   G_{D_s} \frac{2\ii }{\sqrt{V}}\nabla_{\Theta_4}^sG_{D_s}.\\
\label{tcomG}
\left[t^j,G_{D_s}\right] &=   G_{D_s}\left(\frac{c^j}{\sqrt{V}}D_s^++D_s^-\frac{c^j}{\sqrt{V}}\right)G_{D_s}
=   -G_{D_s} \frac{2 }{\sqrt{V}}\nabla_{\Theta_j}G_{D_s}.\\
\label{compi}
\left[\frac{d}{ds},\Pi_{D_s}\right] &=  -\left[\frac{d}{ds},D_s^+G_{D_s}D_s^-\right]
=  -\Pi_{D_s}\frac{\ii c^4}{\sqrt{V}}G_{D_s}D_s^- - D_s^+G_{D_s}\frac{ic^4}{\sqrt{V}}\Pi_{D_s}.\\
\label{comt} 
\left[t^j,\Pi_{D_s}\right] &= \Pi_{D_s}\frac{ c^j}{\sqrt{V}}G_{D_s}D_s^- + D_s^+G_{D_s}\frac{c^j}{\sqrt{V}}\Pi_{D_s}\nonumber\\
&= \Pi_{D_s}\frac{ c^4}{\sqrt{V}}G_{D_s}I_jD_s^- + D_s^+I_j^\dagger G_{D_s}\frac{c^4}{\sqrt{V}}\Pi_{D_s}.
\end{align}
\begin{proposition}
$T^j(s)$ and $B_\sigma$ satisfy the moment map equations \eqref{nahmeq0} and \eqref{pEnd}. 
\end{proposition}
\begin{proof}
First we show that the $T^j$ satisfy the Nahm equations \eqref{nahmeq0}. 
We compute 
\begin{align}&[\ii\nabla_{\frac{d}{ds}},T^1] - [T^2,T^3]\nonumber\\
& = 
[\ii\Pi_{D_s}\frac{d}{ds}\Pi_{D_s},\Pi_{D_s}t^1\Pi_{D_s}] - [\Pi_{D_s}t^2\Pi_{D_s},\Pi_{D_s}t^3\Pi_{D_s}]\nonumber\\
 &=2\ii\Pi_{D_s}\left[\frac{d}{ds} ,\Pi_{D_s}\right]t^1\Pi_{D_s} 
  +2\Pi_{D_s}[t^3,\Pi_{D_s}]t^2\Pi_{D_s}\nonumber\\
&= 2\Pi_{D_s}\frac{ c^4}{\sqrt{V}}G_{D_s} \frac{c^1-c^4c^3c^2}{\sqrt{V}}\Pi_{D_s}=0,
\end{align}
since $c^1-c^4c^3c^2=c^1(1+\gamma)$ annihilates $S^-$.  
This establishes moment map equation \eqref{nahmeq0}.

We next consider moment map equation \eqref{pEnd}. Using \eqref{Eq:bCov}, we compute 
\begin{align}
B_\sigma^\alpha B_\sigma^{\beta\dagger} &= \Pi_{D_{p_\sigma-}}b_\sigma^\alpha\Pi_{D_{p_\sigma+}}\bar b_\sigma^\beta\Pi_{D_{p_\sigma-}}\nonumber\\
&= \Pi_{D_{p_\sigma-}}b_\sigma^\alpha(I-D_{p_\sigma+}^+G_{D_{p_\sigma+}}D_{p_\sigma+}^-)\bar b_\sigma^\beta\Pi_{D_{p_\sigma-}}\nonumber\\
&= \Pi_{D_{p_\sigma-}}b_\sigma^\alpha \bar b_\sigma^\beta\Pi_{D_{p_\sigma-}}
+\Pi_{D_{p_\sigma-}}c^4I_b\frac{(\e_b^\dagger b_\sigma)^\alpha}{2\sqrt{V}t_\sigma} G_{D_{p_\sigma+}}I_a^\dagger c^4\frac{\overline{(\e_a^\dagger b_\sigma)^\beta}}{2\sqrt{V}t_\sigma}\Pi_{D_{p_\sigma-}}.\label{TD}
\end{align}

The last summand in  \eqref{TD} has the form 
{$\Pi_{D_{p_\sigma-}}c^4I_b I_a^\dagger (I_b X I_a^\dagger)^{\alpha\beta} c^4 \Pi_{D_{p_\sigma-}}$} 
and thus it is proportional to $\delta^{\alpha\beta}$ by the quaternionic identity  \eqref{eq:tens}.   
Hence this term cannot contribute to $\mathfrak{Im}\, B_\sigma B_\sigma^\dagger.$ 
On the other hand, by \eqref{bbnd},  
\begin{align}\label{bbjk}b_\sigma^\alpha \bar b_\sigma^\beta = (b_\sigma b_\sigma^\dagger)^{\alpha\beta} = (\ii {\e_j }(t^j-\nu_\sigma^j) + t_\sigma 1_{S})^{\alpha\beta}.
\end{align} 
Substituting \eqref{bbjk} back into \eqref{TD}, we prove the first relation of \eqref{pEnd}. The proof of the second is essentially the same. 
\end{proof} 

Before we can check the remaining moment map equation \eqref{Tdisc}, 
we must, of course, first define the fundamental data $Q_\lambda, \lambda\in \Lambda^0$. This requires an extensive study of bounded harmonic sections of $\mathcal{E}$, which is the subject of the following section.

\subsection{Fundamental Data}\label{fud}
{In this subsection, we define the fundamental bow data  $Q$ produced by the Down Transform. Extensive analytic preliminaries are required to show that the data is well defined and has the requisite properties. For the convenience of the reader wishing to understand the overall structure of the transform before plunging into analytical detail, we will state the required analytical preliminaries in this subsection, but postpone their proofs until later subsections.}

For the remainder of this subsection, we consider $\lambda\in \Lambda^0$. It will be useful in this and many subsequent sections to replace the coordinate frame $\{\Theta_a\}_{a=1}^4$ with a frame $\{\hat \Theta_a\}_{a=1}^4$,   with $\hat \Theta_1$ the unit vector in the direction of the horizontal lift of the radial vectorfield  in $\IR^3$, and $\hat\Theta_4 = \Theta_4$. We will abuse notation slightly and write $\frac{\p}{\p r}$ for both the radial vectorfield in $\IR^3$ and its horizontal lift.  We let $\hat c^a$ denote Clifford multiplication by the dual of $\hat \Theta_a$, but will usually write $c^4$, since $c^4 = \hat c^4$.

Let $W_\lambda:=\mathrm{Ker}_{L^\infty} \nabla^{\lambda *}\nabla^\lambda $ denote the space of bounded harmonic sections of ${\mathcal E}\otimes \Fe_{\lambda}$. Since $S^+$ with the Levi-Civita spin connection is a trivial $\IC^2$ bundle with the product connection, the bounded harmonic sections of $S^+\otimes{\mathcal E}\otimes \Fe_{\lambda}$ can (and will) be identified with  $\IC^2\otimes W_\lambda$. We have the following propositions concerning  bounded harmonic spinors, {whose  proofs  can be found in Sec.~\ref{BHS},  
pages \pageref{Proof:bddh2} and \pageref{Proof:wdim1}.}

\begin{proposition}\label{bddh2}
Let $v$ be a covariant constant section of $S^+$. For every 
$H\in W_\lambda$,  $\|r^{2}\nabla^\lambda H\|_{L^\infty(M)}<\infty$, and $D_\lambda (v\otimes H)\in L^2$. Moreover, $W_\lambda$ is an inner product space endowed with the inner product
\begin{align}\label{normalize}
\langle H_1,H_2\rangle_\infty:=  \frac{1}{2} \lim_{R\to\infty}\int_{\pi_k^{-1}(S_R^2)}\langle H_1,H_2\rangle\,  d\sigma,
\end{align}
 where $d\sigma$ is the  volume form of the round unit sphere.   
\end{proposition}
\begin{proposition}\label{wdim1}
$W_\lambda$ is at most one-dimensional.
\end{proposition}

Let $\{\Psi_\alpha(\lambda)\}_\alpha$ be an $L^2$ unitary basis of $\Ker_{L^2}(D_\lambda) $. We show in the proof of Proposition \ref{residue} that for $s$ small, $\Pi_{D_{\lambda+s}}: \Ker_{L^2}(D_\lambda) \to  \Ker_{L^2}(D_{\lambda+s}) $ is an isomorphism, and $\{\Pi_{D_{\lambda+s}}\Psi_\alpha(\lambda)\}_\alpha$ therefore defines a natural local frame for $\mathpzc{E}$ in a neighborhood of $\lambda$. This frame can be somewhat awkward to work with; so, in \eqref{exts} we define a more computable modification $\{\Psi_\alpha(\lambda+s)\}_\alpha$ of this local frame which differs from  it by $O(|s|\ln\frac{1}{|s|})$ corrections. 

Define 
\begin{align}\label{eq:chi} 
	\chi_\alpha(\lambda+s): = G_{D_{\lambda+s}}\frac{c^4}{\sqrt{V}}\Psi_{\alpha}(\lambda+s).
\end{align}

\begin{proposition}\label{cool}
There is a bounded harmonic section 
$U_{0,\alpha} = i \hat c^1 \sqrt{V}r^2  \Psi_\alpha(\lambda)  +   O\left(\frac{1}{\sqrt{1+r^2}}\right)$ such that 
\begin{align}&\lim_{s\to 0^+}\left[G_{D_{\lambda+s}}\frac{c^4}{\sqrt{V}}\Pi_{D_{\lambda+s}}\Psi_\alpha(\lambda) - G_{D_{\lambda-s}}\frac{c^4}{\sqrt{V}}\Pi_{D_{\lambda-s}}\Psi_\alpha(\lambda)\right]\nonumber\\
  = &\lim_{s\to 0^+}\left[\chi_{\alpha}(\lambda+s) - \chi_{\alpha}(\lambda-s)\right] =  U_{0,\alpha}.
\end{align}
\end{proposition} 
See Sec.~\ref{flbh}, 
page \pageref{postcool} for the proof of this Proposition.

Let $ q_\lambda  :\Ker_{L^2}(D_\lambda) \to S^+\otimes W_\lambda$ denote the map we have just constructed: 
\begin{align}\label{eq:Psas}
	q_\lambda: \Psi(\lambda)\mapsto  \lim_{s\to 0^+}\left[G_{D_{\lambda+s}}\frac{c^4}{\sqrt{V}}\Pi_{D_{\lambda+s}}\Psi(\lambda) - G_{D_{\lambda-s}}\frac{c^4}{\sqrt{V}}\Pi_{D_{\lambda-s}}\Psi(\lambda)\right].
\end{align}
Here, we are implicitly using a trivialization of $S^+$ by unitary covariant constant sections to identify the codomain of this map with   $S^+\otimes W_\lambda$. Let $\{f_a\}_{a=1}^2$ be a covariant constant unitary frame for $(S^+)^*$ with coframe $\{f_a^\dagger\}_{a=1}^2$
for $S^+$.
There is also a natural map 
$R_\lambda:\Ker_{L^2}(D_\lambda) \to S^+\otimes W_\lambda$ given as the adjoint of the linear map defined by 
$$R_\lambda^\dagger(f_a^\dagger\otimes w) = -\ii D_\lambda (f_a^\dagger \otimes w).$$
Here the adjoint is taken with respect to the inner product on $\IC^2\otimes W_\lambda$ defined in Proposition \ref{bddh2}. 
\begin{proposition}\label{R2q}
$R_\lambda =  2q_\lambda.$
\end{proposition}
See  Sec.~\ref{flbh}, page~\pageref{proofR2q} for the proof.

We finally define the fundamental data  $Q_\lambda\in \Hom(W_\lambda,\IC^2\otimes (\Ker_{L^2}(D_\lambda) ))$. 
Recall that the trivial $\IC^2$ bundle on the bow descends to the trivial $\IC^2$ bundle {$\calS$} on  $\TN_k^{\nu}$, which we identify with $(S^+)^\ast.$  Define
\begin{align}\label{disQ1}
		Q_\lambda(w) := -\ii f_a\otimes\Pi_{D_\lambda}D_\lambda (f_a^\dagger\otimes w)=-\ii f_a\otimes D_\lambda (f_a^\dagger\otimes w).
\end{align} 
 To determine the adjoint of $Q_\lambda$, we compute for $\psi_a\in \Ker_{L^2}(D_\lambda),$ $a=1,2$, 
\begin{align}\label{disQ2}
	\langle Q_\lambda(w),f_a\otimes \psi_a\rangle_{L^2}&= \langle -\ii f_b\otimes D_\lambda (f_b^\dagger\otimes w),f_a\otimes \psi_a\rangle_{L^2} \nonumber\\ 
&= \langle R_\lambda^\dagger (f_b^\dagger\otimes w),  \psi_b\rangle_{L^2} \nonumber\\ 
&=  \langle  f_a^\dagger\otimes w , 2q_\lambda \psi_b\rangle_{L^2} \nonumber\\  
&=  \langle    w , 2q_\lambda (\psi_b)(f_b)\rangle_{L^2}   .
\end{align}
Hence we have 
\begin{align}Q_\lambda^\dagger(f_a\otimes \psi_a) = 2q_\lambda (\psi_b)(f_b).
\end{align} 
 
We now establish the remaining moment map equation \eqref{Tdisc}. Let $f_a\otimes \psi_a(\lambda) $ (implicit sum on $a$) be an element of $\IC^2\otimes \Ker_{L^2}(D_\lambda^-),$ thanks to Lemma~\ref{trivKer}. Let $\psi_a(\lambda+s) := \Pi_{D_{\lambda+s}}\psi_a(\lambda)$. 
We compute
\begin{align}&\mathfrak{Im}\, i(Q_\lambda Q_\lambda^\dagger) f_b\otimes \psi_b(\lambda)\nonumber\\
 &=    f_a\otimes D_\lambda    q_\lambda(\psi_b)\langle f_b\otimes f_a^\dagger,I_k^\dagger\rangle I_k^\dagger  \nonumber\\
&=    \e_kf_a\otimes D_\lambda    I_k q_\lambda(\psi_a)     \nonumber\\
& =   \e_kf_a\otimes  \lim_{s\to 0^+}D_\lambda I_k  \left(G_{\lambda+s}\frac{c^4}{\sqrt{V}} \psi_a(\lambda+s)-G_{\lambda-s}\frac{c^4}{\sqrt{V}} \psi_a(\lambda-s)\right)   \nonumber\\
& = -  \e_kf_a\otimes  \lim_{s\to 0^+}D_\lambda   \left(G_{\lambda+s} \frac{c^k}{\sqrt{V}} \psi_a(\lambda+s)-G_{\lambda-s}\frac{ c^k}{\sqrt{V}} \psi_a(\lambda-s)\right)   \nonumber\\
& = -  \e_kf_a\otimes  \lim_{s\to 0^+}D_\lambda   \left(G_{\lambda+s} D_{\lambda+s}t^k \psi_a(\lambda+s)-G_{\lambda-s}D_{\lambda-s}t^k \psi_a(\lambda-s)\right)   \nonumber\\
& = -  \e_kf_a\otimes  \lim_{s\to 0^+}    \left(D_{\lambda+s}G_{\lambda+s} D_{\lambda+s}t^k \psi_a(\lambda+s)-D_{\lambda-s}G_{\lambda-s}D_{\lambda-s}t^k \psi_a(\lambda-s)\right)   \nonumber\\
&\quad +  \e_kf_a\otimes  \lim_{s\to 0^+}\frac{isc^4}{\sqrt{V}}   \left(G_{\lambda+s} D_{\lambda+s}t^k \psi_a(\lambda+s)+G_{\lambda-s}D_{\lambda-s}t^k \psi_a(\lambda-s)\right)   \nonumber\\
& = \e_kf_a\otimes  \lim_{s\to 0^+}    \left(T^k(\lambda+s) \psi_a(\lambda+s)-T^k(\lambda-s) \psi_a(\lambda-s)\right)   \nonumber\\
&\quad +  \e_kf_a\otimes  \lim_{s\to 0^+}\frac{isc^4}{\sqrt{V}}   \left(G_{\lambda+s} \frac{ c^k}{\sqrt{V}}  \psi_a(\lambda+s)+G_{\lambda-s}\frac{ c^k}{\sqrt{V}} \psi_a(\lambda-s)\right)   \nonumber\\
& = \e_kf_a\otimes  \lim_{s\to 0^+}    \left(T^k(\lambda+s) \psi_a(\lambda+s)-T^k(\lambda-s) \psi_a(\lambda-s)\right)   \nonumber\\
& = \e_kf_a\otimes  \lim_{s\to 0^+}    \left(T^k(\lambda+s) -T^k(\lambda-s)\right) \psi_a(\lambda), 
\end{align}
and \eqref{Tdisc} follows.
 In the next {three} subsections we derive the properties of $L^2$ strongly harmonic spinors and of bounded harmonic sections that we used in this subsection. 
 \subsubsection{\texorpdfstring{$L^2$}{L2} Harmonic spinors\texorpdfstring{ near $\lambda\in \Lambda^0$}{ }}\label{LHS}
In this subsection we gather some properties of strongly harmonic spinors at $\lambda+s$, for $\lambda\in\Lambda^0$ and $|s|$ small.   
First we recall some useful notation and estimates from   \cite[Eq.~(45), Eq.~(130), and Sec.~6.3]{First}. 

Fix $R_0>\sum_\sigma |\nu_\sigma|$. Given  sections $\zeta_1$ and $\zeta_2$ of a Hermitian bundle over $\TN_k$, we define for $|x|\geq R_0$, 
\begin{align*}
\Phi(\zeta_1)(x)&:= \int_{\pi_k^{-1}(x)}|\zeta_1|^2d\tau,& 
&\text{ and }& 
Q(\zeta_1,\zeta_2)(x)&:= \mathfrak{Re}\,\int_{\pi_k^{-1}(x)}\langle \zeta_1,\zeta_2\rangle d\tau.
\end{align*}
For $|x|\geq R_0$,  we define projection operators on $L^2(\pi_k^{-1}(x))$ with coefficients in any Riemannian bundle tensored with $\mathcal{E}\otimes \Fe_{\lambda+s}$  by 
\begin{align}\label{pi0}
\Pi_0 = \frac{1}{2\pi i} \oint_{C}\left(z-\ii\nabla^{\lambda}_{\frac{\p}{\p \tau}}\right)^{-1}dz,
\end{align}
where $C$ is a small circle around $0$. Set $\Pi_1:= 1-\Pi_0.$ Then there is a constant $c_\pi>0$, so that for all sections $\sigma$, 
\begin{align}\label{p1ilam}
\Phi\left(\Pi_1\nabla_{\Theta_4}^\lambda\sigma \right)\geq c_\pi^2\Phi(\Pi_1\sigma).
\end{align}
Moreover, by \cite[Eq.~(130)]{First} there is a constant $c_m>0$, so that for all sections $\sigma$, 
\begin{align}\label{p0ilam}
\Phi\left(\Pi_0\nabla_{\Theta_4}^\lambda\sigma\right)\leq c_m^2r^{-4}\Phi(\Pi_0\sigma).
\end{align}
Consequently 
\begin{align}\label{p0is}\Phi\left(\Pi_0\nabla_{\Theta_4}^{\lambda+s}\sigma\right)\geq \left(\frac{s^2}{V}-\frac{2c_ms}{r^2\sqrt{V}}\right)\Phi(\Pi_0\sigma),
\end{align}
and  
\begin{align}\label{p1is}\Phi\left(\Pi_1\nabla_{\Theta_4}^{\lambda+s}\sigma\right)\geq \left(\frac{1}{2}c_\pi^2 - \frac{s^2}{V}\right) \Phi(\Pi_1\sigma) .
\end{align}
By  \cite[Eq.~(131)]{First} and the cubic decay of $R_{ab}$, there exists $C_F>0$ such that 
\begin{align}\left| \langle \Pi_02(F_{ab}^{\lambda+s}+R_{ab})\nabla^{\lambda+s}_{\Theta_a}\sigma,\Pi_0\nabla^{\lambda+s}_{\Theta_b}\sigma\rangle \right|\leq \frac{C_F}{(1+r^3)}\Phi(\nabla^{\lambda+s}\sigma).
\end{align}
We will frequently use these inequalities to control the respective contributions of  $\Pi_0$ and $\Pi_1$ to our estimates.

Define 
$$r_{s,\tau}:= \frac{r}{(1+\tau^2s^2r^2)^{\frac{1}{2}}},$$
and set $r_s:= r_{s,1}$. {(Here we do not combine $s$ and $\tau$ into a single parameter $s\tau$, because $s$ arises from the parameter defining the connection, and it is convenient to make the dependence on this parameter explicit.) }  Then 
$$\frac{\p r_{s,\tau}}{\p r}=  \frac{r_{s,\tau}}{r(1+\tau^2s^2r^2)} .$$
 Also, define
\begin{align}
	\eta_{\sigma,t}&:=
	\frac{r^te^{\sqrt{\sigma^2r^2+t^2}}}{\left(1+\sqrt{\frac{\sigma^2r^2}{t^2}+1}\right)^t}\  \text{for } t\geq 0&
	&\text{and}&
	\eta_{\sigma,t}&:=r^te^{ \sigma r}\  \text{for }  t<0.
\end{align}
Then for $t\geq 0$, 
\begin{align}\frac{|d\eta_{\sigma,t}|}{\eta_{\sigma,t}} =  V^{-\frac{1}{2}}\sqrt{\frac{t^2}{r^2}+\sigma^2},
\end{align}
and for $t<0$
\begin{align}\frac{|d\eta_{\sigma,t}|}{\eta_{\sigma,t}} \leq  V^{-\frac{1}{2}}\sqrt{\frac{t^2}{r^2}+\sigma^2}.
\end{align}
Observe that for  some constants $c_1(t)$, $c_2(t)$ (depending on $t$)   
\begin{align}
c_1(t)r_s^te^{|s|r}\leq \eta_{s,t}\leq c_2(t)r_s^te^{|s|r}.
\end{align} 
We will frequently use the relation: 
\begin{align}\nabla^{ \lambda+s *}\nabla^{\lambda+s} = \nabla^{\lambda *}\nabla^\lambda+\frac{s^2}{V}- 2\ii s\nabla^{\lambda}_{ \frac{\p}{\p\tau}}.
\end{align}

We now record a basic technical lemma we will use often. 
\begin{lemma}
\label{oldnag}
 Let $w\in L^2(S^+\otimes \mathcal{E})$ such that  for some  $t^2\leq \frac{1}{4}-2c_m|s|$ and some $0\leq \sigma \leq |s|$, 
$ \frac{\eta_{\sigma,t}w}{\sqrt{(\frac{1}{4}-t^2-2c_m|s|)r^{-2}+(s^2-\sigma^2)}}\in L^2$ and $r^{1+p}w\in L^2$. Set $K_s:= G_{D_{\lambda+s}}w$. Then for  $p\in(-\frac{1}{2},\frac{1}{2})$, $\exists C_1=C_1(R_0)>0, C_2 = C_2(R_0)\geq 1$ independent of $w$,  $s$, and $\sigma$ for $|s|$ small,  such that
\begin{align}\label{mus1}
&\left\|\sqrt{(\frac{1}{4}-t^2-2c_m|s|)r^{-2}+(s^2-\sigma^2)}\frac{\eta_{\sigma,t}K_s}{\sqrt{V}}\right\|^2_{L^2}\nonumber\\
&  \leq  \frac{C_1}{1-4p^2}\| \sqrt{V}r^{1+p}w\|_{L^2}^2
+ C_2\left\|\frac{\sqrt{V}\eta_{\sigma,t}w}{\sqrt{(\frac{1}{4}-t^2-2c_m|s|)r^{-2}+(s^2-\sigma^2)}} \right\|^2_{L^2}.
\end{align}
When  $|t|<\frac{\tau}{2}$, for some $\tau\in [0,1)$, $\exists \tilde C >0$ such that  
\begin{multline}\label{mus2} 
\| \eta_{\sigma,t}\nabla^{\lambda+s}K_s \|^2_{L^2}+\|\nabla^{\lambda+s}(\eta_{\sigma,t}K_s )\|^2_{L^2}\\
  \leq \frac{\tilde C }{1-\tau^2}\Bigg( \frac{C_1}{1-4p^2}\| \sqrt{V}r^{1+p}w\|_{L^2}^2\\
+ C_2\left\|\frac{\sqrt{V}\eta_{\sigma,t}w}{\sqrt{(\frac{1}{4}-t^2-2c_m|s|)r^{-2}+(s^2-\sigma^2)}} \right\|^2_{L^2}\Bigg).
\end{multline}
\end{lemma}
 \begin{proof} We treat the case $t\sigma\geq 0$. The estimates simplify in the case $t\sigma<0$, but otherwise the proof is the same. So, for the remainder of the proof, we assume $t\sigma\geq 0$. 
First we have the simple Hardy inequality estimate (see \cite[Lem.~13]{First}) for $p\in (-\frac{1}{2},\frac{1}{2})$, 
\begin{align}\langle w,r^{2p} K_s\rangle_{L^2}&= \left\|\nabla^{\lambda+s}(r^pK_s)\right\|^2_{L^2}-p^2\left\|\frac{r^{p-1}K_s}{\sqrt{V}}\right\|^2_{L^2}\nonumber\\
& \geq (\frac{1}{4}-p^2)\left\|\frac{r^{p-1}K_s}{\sqrt{V}}\right\|^2_{L^2}.
\end{align}
Hence 
\begin{align}\label{innest}& (\frac{1}{4}-p^2)\left\|\frac{r^{p-1}K_s}{\sqrt{V}}\right\|_{L^2}\leq  \left\| \sqrt{V}r^{p+1}w\right\|_{L^2}.
\end{align}
Let $\eta_{\sigma,t,N} := \min\{\eta_{\sigma,t},N\}$. Let ${\chi}_{d\eta}$ denote the characteristic function of the support of $N-\eta_{\sigma,t,N}$. Let ${\chi}_{R_0}$ denote the characteristic function of $\pi_k^{-1}(B_{R_0})$. Then we have 
\begin{align}\label{localmn}
&\langle w,\eta^2_{\sigma,t,N} K_s\rangle_{L^2}\nonumber\\
&= \left\|\nabla^{\lambda+s}(\eta_{\sigma,t,N}K_s)\right\|^2_{L^2}
-\left\| \sqrt{\frac{t^2}{r^2}+\sigma^2}\frac{\chi_{d\eta} \eta_{\sigma,t,N}K_s }{\sqrt{V}}\right\|^2_{L^2}\nonumber\\
&= \left\|(\nabla^{\lambda}+\frac{ \chi_{R_0} is\theta^4\otimes}{\sqrt{V}})(\eta_{\sigma,t,N}K_s)\right\|^2_{L^2}+s^2\left\| (1-\chi_{R_0})\frac{\eta_{\sigma,t,N}K_s}{\sqrt{V}}\right\|^2_{L^2}\nonumber\\
&\quad-2\left\langle is\nabla^{\lambda}_{\frac{\p}{\p\tau}}\eta_{\sigma,t,N}K_s,(1-{\chi}_{R_0})\eta_{\sigma,t,N}K_s\right\rangle_{L^2}
- \left\| \sqrt{\frac{t^2}{r^2}+\sigma^2} \frac{{\chi}_{d\eta} \eta_{\sigma,t,N}K_s }{\sqrt{V}}\right\|^2_{L^2}\nonumber\\
&\geq \left(\frac{1}{4}-t^2-2c_m|s|\right)\left \|\frac{\eta_{\sigma,t,N}K_s}{\sqrt{V}r}\right\|^2_{L^2}+(s^2-\sigma^2)\left\|  \frac{\eta_{\sigma,t,N}K_s}{\sqrt{V}}\right\|^2_{L^2}\nonumber\\
&\quad-C_0\int_{\pi_k^{-1}(B_{R_0})}| \eta_{\sigma,t,N}K_s|^2dv\nonumber\\
&\geq \left(\frac{1}{4}-t^2-2c_m|s|\right)\left\|\frac{\eta_{\sigma,t,N}K_s}{\sqrt{V}r}\right\|^2_{L^2}+(s^2-\sigma^2)\left\|  \frac{\eta_{\sigma,t,N}K_s}{\sqrt{V}}\right\|^2_{L^2}\nonumber\\
&
\quad -\frac{\tilde C_1}{(1-4p^2)^2} \left\| \sqrt{V}r^{p+1}w\right\|_{L^2}^2,
\end{align}
for some $C_0, \tilde C_1= \tilde C_1(R_0)>0$, independent of $w$ and $s$ (for $|s| $ small).  Hence there exist $C_1=   C_1(R_0)>0$ and $C_2=  C_2(R_0)\geq 1$ such that 
\begin{multline} 
\left\|\sqrt{(\frac{1}{4}-t^2-2c_m|s|)r^{-2}+(s^2-\sigma^2)}\frac{\eta_{\sigma,t,N}K_s}{\sqrt{V}}\right\|^2_{L^2}\\
  \leq \frac{C_1}{(1-4p^2)^2}\left\| \sqrt{V}r^{p+1}w\right\|_{L^2}^2 \qquad\\
  + C_2\left\|\frac{\sqrt{V}\eta_{\sigma,t,N}w}{\sqrt{(\frac{1}{4}-t^2-2c_m|s|)r^{-2}+(s^2-\sigma^2)}} \right\|^2_{L^2}.
\end{multline}
Take the limit as $N\to \infty$ to obtain   inequality \eqref{mus1}. 

It is now straightforward to obtain \eqref{mus2}. By \eqref{mus1} we have
\begin{align}\label{partway}
&\|\nabla^{\lambda+s}(\eta_{\sigma,t}K_s)\|^2_{L^2}\nonumber\\
&= \langle w,\eta^2_{\sigma,t } K_s\rangle_{L^2} +\|d\eta_{\sigma,t}\otimes K_s\|^2_{L^2}
\nonumber\\
&\leq   \frac{C_1}{(1-4p^2)^2}\left\| \sqrt{V}r^{p+1}w\right\|_{L^2}^2 +\left\|V^{-\frac{1}{2}}\eta_{\sigma,t}\sqrt{ t^2 r^{-2}+\sigma^2}  K_s\right\|^2_{L^2}\nonumber \\
&\qquad + C_2\left\|\frac{\sqrt{V}\eta_{\sigma,t}w}{\sqrt{(\frac{1}{4}-t^2-2c_m|s|)r^{-2}+(s^2-\sigma^2)}} \right\|^2_{L^2}.
\end{align}
Using the Leibniz rule, the triangle inequality, \eqref{mus1}, \eqref{innest}, and \eqref{partway},  the claim \eqref{mus2} follows easily.
\end{proof}
\begin{proposition}\label{25}Let 
$\Psi\in \mathrm{Ker}_{L^2} (D_s)$. Then $\exists C_K>0$ such that for $\tau= 6\delta^{-\frac{3}{4}},$ and  $\forall \delta>0$,  
\begin{align}\label{speaker}\sqrt{\delta}\left\|r^{-1} r_{s,\tau}^{\frac{3}{2}-\delta} e^{|s|r}\Psi \right\|_{L^2}  \leq C_K\|\Psi\|_{L^2}.
\end{align}
\end{proposition}
\begin{proof}  We first consider an improved Kato type inequality. Consider
\begin{align}&\left|\nabla_{\hat\Theta_1} \Psi-\hat c^1c^4\nabla^{\lambda+s}_{\hat\Theta_4}\Psi \right|^2= 
\left|\hat{c}^2\nabla_{\hat\Theta_2}\Psi + \hat{c}^3 \nabla_{\hat\Theta_3}\Psi \right|^2\nonumber\\
&\leq  
2\left(|\nabla^{\lambda+s}\Psi|^2-|\nabla_{\hat\Theta_1}\Psi|^2- |\nabla^{\lambda+s}_{\hat\Theta_4}\Psi|^2\right) .
\end{align}
Hence 
\begin{align}\label{ufo}
&\frac{1}{2}\left|\nabla_{\hat\Theta_1}\Psi- \hat{c}^1c^4\nabla^{\lambda+s}_{\hat\Theta_4}\Psi \right|^2+ \left|\nabla_{\hat\Theta_1}\Psi\right|^2+ \left|\nabla^{\lambda+s}_{\hat\Theta_4}\Psi\right|^2\leq  
|\nabla^{\lambda+s}\Psi|^2  .
\end{align}
Outside a compact set we use \eqref{p0is} and \eqref{p1is} to rewrite \eqref{ufo} :
\begin{align}\label{uf0}
&\Phi(\nabla^{\lambda+s}\Psi)-\Phi(\nabla_{\hat\Theta_1}\Psi)  \nonumber\\
&\geq \frac{1}{2}\Phi\left(\Pi_0(\nabla_{\hat\Theta_1}\Psi-\hat c^1c^4\nabla^{\lambda+s}_{\hat\Theta_4}\Psi)\right)+\Phi\left(\nabla^{\lambda+s}_{\hat\Theta_4}\Psi\right) \nonumber\\
&= \frac{1}{2}\Phi\left(\Pi_0(\nabla_{\hat\Theta_1}\Psi-\frac{i\hat c^1c^4 s}{\sqrt{V}}\Psi)\right)+Q\left(\Pi_0(\nabla_{\hat\Theta_1}\Psi-\frac{i\hat c^1c^4 s}{\sqrt{V}} \Psi),\Pi_0(  \nabla^{\lambda }_{\hat\Theta_4}\Psi)\right)\nonumber\\
&\quad + \frac{1}{2}\Phi(\Pi_0(  \nabla^{\lambda }_{\hat\Theta_4}\Psi)) +  \Phi(\Pi_0(  \nabla^{\lambda +s}_{\hat\Theta_4}\Psi))+\Phi(\Pi_1\nabla^{\lambda+s}_{\hat\Theta_4}\Psi) \nonumber\\
&\geq \left(\frac{1}{2}-\frac{1}{r}\right)\Phi\left(\Pi_0(\nabla_{\hat\Theta_1}\Psi-\frac{i\hat c^1c^4 s}{\sqrt{V}}\Psi)\right) +\left(\frac{s^2}{V}-\frac{2c_m|s|}{r^2\sqrt{V}}  -\frac{c_m^2}{r^3}\right)\Phi(\Pi_0(  \Psi))\nonumber\\
&\quad +\left(\frac{3}{4}c_\pi^2 - 3\frac{s^2}{V}\right) \Phi(\Pi_1\Psi)   .
\end{align} 
Let  $P_{\pm}$ denote orthogonal projection onto the $\pm 1$ eigenspaces of $i\hat c^1c^4$. Since $[\nabla_{\hat\Theta_1},i\hat c^1c^4] = O(r^{-2})$, we have 
\begin{align}\label{pcom1}[\nabla_{\hat\Theta_1},P_{\pm} ] =  O(r^{-2}) .
\end{align}
Similarly by \cite[Eq.~(131)]{First}, 
\begin{align}\label{pcom2}[\nabla_{\hat\Theta_1},\Pi_0 ] =  O(r^{-2}) .
\end{align}
Using \eqref{pcom1} and \eqref{pcom2}, we rewrite \eqref{uf0} for some $C_3>0$, as 
\begin{align}\label{uf1}&\Phi(\nabla^{\lambda+s}\Psi)-\Phi(\nabla_{\hat\Theta_1}\Psi)  \nonumber\\
&\geq \left(\frac{1}{2}-\frac{2}{r}\right)\Phi( \nabla_{\hat\Theta_1}\Pi_0  \Psi)+\left(\frac{1}{2}-\frac{2}{r}\right)\frac{s^2}{V}\Phi(\Pi_0 \Psi)\nonumber\\
&\quad  +\frac{s}{V}\frac{\p}{\p r}\left(\frac{1}{2}-\frac{2}{r}\right)(\Phi( \Pi_0P_- \Psi)-\Phi( \Pi_0P_+ \Psi))\nonumber\\
&\quad + \frac{s^2}{V} \Phi(\Pi_0(  \Psi))+\left(\frac{3}{4}c_\pi^2 - 3\frac{s^2}{V}\right) \Phi(\Pi_1\Psi)  -C_3(r^{-3}+|s|r^{-2})\Phi( \Psi) .
\end{align} 
 For $0\leq \sigma<|s|<<c_\pi$, $p<\frac{3}{2}$,  using \eqref{uf1} and then Hardy's inequality, we have, 
\begin{align}&\int_{\TN_k} \langle -c(F_s)r_{s,\tau}^{p} e^{\sigma r}\Psi,r_{s,\tau}^{p} e^{\sigma r}\Psi\rangle dv \nonumber\\
&= \int_{\TN_k} |\nabla^s(r_{s,\tau}^{p}e^{\sigma r}\Psi)|^2dv-\int_{TN_k} \left|(\frac{p}{r(1+\tau^2s^2r^2)}+\sigma)r_{s,\tau}^{p}e^{\sigma r}\Psi \right|^2V^{-1}dv
\nonumber\\
&\geq  \int\limits_{\TN_k}   |\nabla_{\frac{\p}{\p r}}(r_{s,\tau}^{p}e^{\sigma r}\Psi)|^2 V^{-1}dv
-\int_{\TN_k} \left|\left(\frac{p}{r(1+\tau^2s^2r^2)}+\sigma\right)r_{s,\tau}^{p}e^{\sigma r}\Psi \right|^2V^{-1}dv
\nonumber\\
&\quad +\int_{\TN_k}\biggl[\left(\frac{1}{2}-\frac{2}{r}\right)V^{-1}\Phi( r_{s,\tau}^{p}e^{\sigma r}\nabla_{\frac{\p}{\p r}}\Pi_0  \Psi)+ \frac{3}{2} \frac{s^2}{V}\Phi(r_{s,\tau}^{p}e^{\sigma r}\Pi_0 \Psi)\nonumber\\
&\quad  +\frac{s}{V}r_{s,\tau}^{2p}e^{2\sigma r}\frac{\p}{\p r}\left(\frac{1}{2}-\frac{2}{r}\right)(\Phi( \Pi_0P_- \Psi)-\Phi( \Pi_0P_+ \Psi))\nonumber\\
&\quad  + \frac{c_\pi^2}{2}   \Phi(\Pi_1\Psi)  -O(r^{-3}+s^2r^{-1}+sr^{-2})\Phi( \Psi)\biggr]dv
\nonumber\\
&\geq  \int\limits_{\TN_k}  \bigl[  \frac{(1+2p)(3-2p) +(6-4p)\tau^2s^2r^2+3\tau^4s^4r^4 -8pr(\sigma+|s|)(1+\tau^2s^2r^2) }{8r^2(1+\tau^2s^2r^2)^2}
 \nonumber\\
&\quad  +     \frac{\sigma-|s|}{r}    
+\frac{(3|s|+\sigma)(|s| -\sigma)}{2}    -O(r^{-3}+s^2r^{-1}+sr^{-2})\bigr]|r_{s,\tau}^{p}e^{\sigma r}\Psi |^2V^{-1}dv
.
\end{align}
Write $p= \frac{3}{2}-\delta$, with $0<\delta\leq 1$, and choose $\sigma = (1-\epsilon)|s|$, with $0<\epsilon<\frac{\delta}{32}$. Write $X=\tau|s|r$. Choose $\tau = 96\delta^{-\frac{3}{4}}$. Then
we have 
\begin{align}&\int_{\TN_k}O(r^{-3}+s^2r^{-1}+sr^{-2})|r_{s,\tau}^{p}e^{\sigma r}\Psi |^2V^{-1}dv\nonumber\\
&\geq  \delta\int_{\TN_k}  \biggl[  \frac{  1 +    X^2+\frac{3}{4\delta}X^4 - \frac{1}{16} \delta^{-\frac{1}{4}}   (X+X^3) }{2r^2(1+X^2)^2}
 \nonumber\\
&\quad  +\frac{\epsilon}{\delta}\frac{ (4-\epsilon)    \tau^{-2}X^2(1+X^2)^2- \tau^{-1}X(1+X^2)^2 }{2r^2(1+X^2)^2} \biggr]|r_{s,\tau}^{p}e^{\sigma r}\Psi |^2V^{-1}dv
\nonumber\\
&\geq \frac{\delta}{4}\|r^{-1}r_{s,\tau}^{p}e^{\sigma r}\Psi \|^2_{L^2}.
\end{align}

Consequently,    we have 
$\delta\|r^{-1}r_{s,\tau}^{\frac{3}{2}-\delta}e^{\sigma r}\Psi\|^2_{L^2}\leq C_K \|\Psi\|^2_{L^2},$
for some $C_K$ independent of $s,\delta,\epsilon$. Taking the limit as $\epsilon\to 0$ with $\tau = 96\delta^{-\frac{3}{4}}$ gives   
\begin{align}\label{point}\delta\|r^{-1}r_{s,\tau}^{\frac{3}{2}-\delta}e^{|s| r}\Psi\|^2_{L^2} \leq C_K \|\Psi\|^2_{L^2}.
\end{align}
\end{proof}
\begin{coro}\label{water}For some $B_e>0$ independent of $s$, 
\begin{align}|\Psi|^2\leq   B_e\left(\frac{1}{r^4}+\frac{|s|^3}{r}\right)\ln(r)(1+\ln(r)^{\frac{3}{2}}s^2r^2)^{3}e^{-2|s|r} \|\Psi\|^2_{L^2}.
\end{align}
\end{coro}
\begin{proof}Applying the elliptic estimate \cite[Prop.~2]{First} to $f=\Phi(\Psi)$, with $W=O(1)$, and 
$R=\min\{\frac{1}{2}|x|,\frac{1}{|s|}\}$, we deduce, using \eqref{point} that for some $\tilde B_e>0$, independent of $s$, 
\begin{align*}\Phi(\Psi)(x)\leq \tilde B_e\left(\frac{1}{|x|^3}+|s|^3\right)|x|^{2\delta-1}\delta^{-1}(1+96^2\delta^{-\frac{3}{2}}s^2|x|^2)^{3-2\delta}e^{-2|s||x|} \|\Psi\|^2_{L^2}.
\end{align*}
Choosing $\delta = \frac{1}{\ln(|x|)}$ gives 
\begin{align*}\Phi(\Psi)(x)\leq 8\tilde B_e\left(\frac{1}{|x|^4}+\frac{|s|^3}{|x|}\right) \ln(|x|)(1+96^2\ln(|x|)^{\frac{3}{2}}s^2|x|^2)^{3}e^{-2|s||x|} \|\Psi\|^2_{L^2}.
\end{align*}
The result now follows from   \cite[Lem.~14]{First} (which is stated for powers of $r$, but clearly extends to more general functions of $r$). 
\end{proof}
 
\begin{proposition}\label{residue}
$$\lim_{s\to 0}\Ker_{L^2}(D_{\lambda+s})  = \Ker_{L^2}(D_{\lambda}).$$
\end{proposition}
\begin{proof}For 
$|s|>0$, $D_{\lambda+s}^+$ is Fredholm with zero kernel (by Lemma~\ref{trivKer}), and the index is constant on connected Fredholm families; so,   
 the subspace $\Ker_{L^2}(D_{\lambda+s}) $ has constant rank.  
 We recall that $\Pi_{D_{\lambda+s}}$ denotes $L^2$ orthogonal projection onto 
 $\Ker_{L^2}(D_{\lambda+s}) .$ We will show  that $\Pi_{D_{\lambda+s}}$ is injective on the image of $\Pi_{D_{\lambda}}$ for small $s$, and that $\|(I
-\Pi_{D_{\lambda}})\Pi_{D_{\lambda+s}}\|_{sup}= O(|s|^{\frac{1}{2}-\delta})$,   $\forall\delta>0$, which implies that $\Pi_{D_{\lambda}}$ is injective on the image of $\Pi_{D_{\lambda+s}}$ for small $s$. The proposition follows immediately from these two statements.     

Let $\psi_0\in \Ker_{L^2}(D_\lambda) \cap \Ker(\Pi_{D_{\lambda+s}})$. Then 
$\psi_0 = sD_{\lambda+s}G_{D_{\lambda+s}}  \frac{ic^4}{\sqrt{V}}\psi_0.$  Using \cite[Thm.~29]{First} and Lemma \ref{oldnag}, with $t=0=\sigma=p$, we have 
$$\|\frac{1}{r}\psi_0\|_{L^2}^2=\|\frac{1}{r}sD_{\lambda+s}G_{D_{\lambda+s}}  \frac{\ii c^4}{\sqrt{V}}\psi_0\|^2_{L^2}= O(s\|\psi_0\|_{L^2}^2).$$ Hence $\psi_0 = 0$ and the first injectivity statement follows. 

For the second injectivity statement, consider  $\psi_s\in \Ker_{L^2}(D_{\lambda+s}) .$ Then 
\begin{align}\Pi_{D_{\lambda}}\psi_s = \psi_s - D_\lambda G_\lambda D_{\lambda}\psi_s = \psi_s - \ii sD_\lambda G_\lambda \frac{c^4}{\sqrt{V}}\psi_s. 
\end{align}
Using Lemma \ref{oldnag}, with $w= \frac{c^4}{\sqrt{V}}\psi_s$, and $\sigma = t =p= 0$, and then Proposition \ref{25}, we have 
\begin{align} \left\|\nabla_\lambda G_\lambda \frac{c^4}{\sqrt{V}}\psi_s\right\|_{L^2}&= O(\|r\psi_s\|_{L^2})\nonumber\\
&
= O(\delta^{-9/4}|s|^{-\frac{1}{2}-\delta}\|  r^{-1}r_{s,\tau}^{\frac{3}{2}-\delta}e^{|s|r}\psi_s\|_{L^2})\nonumber\\
&
= O(\delta^{-13/4}|s|^{-\frac{1}{2}-\delta}\|  \psi_s\|_{L^2}).
\end{align}
Consequently 
$\|\psi_s-\Pi_{D_{\lambda}}\psi_s\|_{L^2} = O(|s|^{\frac{1}{2}-\delta})$, for all $\delta>0$, and the proposition follows. 
\end{proof}
We introduce a frame for $\Ker_{L^2}(D_{\lambda+s}) $ near $\lambda$. To orient the reader, we first consider a model problem in $\IR^3\times S^1$. Let $|0\rangle$ be a covariant constant spinor on 
$\IR^3\times S^1$. Let $\tilde D_s = \tilde D + isCl(d\tau)$, where $\tilde D$ is the (untwisted) Dirac operator  on $\IR^3\times S^1$. Then $\Delta \frac{1}{r}|0\rangle = 0=(\Delta +s^2)\frac{e^{-|s|r}}{r}|0\rangle$, away from the origin.  Hence $y_0:= -\tilde D(\frac{1}{r}|0\rangle) =  \frac{ Cl(dr)}{r}|0\rangle$ is in the kernel of $\tilde D$ (away from $r=0$), and 
\begin{align*}
y_s:&=   -\tilde D_s\left(\frac{e^{-|s|r}}{r}|0\rangle\right) =  \left(\frac{Cl(dr)}{r^2}+\frac{Cl(dr)|s|-\ii c^4s}{r}\right)e^{-|s|r}|0\rangle\nonumber\\
& = (1+|s|r+\ii s c^4 Cl(dr))e^{-|s|r}y_0
\end{align*}
is in the kernel of $\tilde D_s$.
 
Returning to $\TN_k$, let $\{\Psi_\alpha(\lambda)\}_\alpha$ be an $L^2$-unitary frame for $\Ker_{L^2}(D_{\lambda}) $. Define the approximate harmonic {spinor}
\begin{align}\label{aph}\tilde \Psi_\alpha(\lambda+s):= (1+|s|r+\ii s r c^4\hat c^1)e^{-|s|r}\Psi_\alpha(\lambda).
\end{align}
We compute
\begin{align}\label{parks}
D_{\lambda+s}\tilde \Psi_\alpha(\lambda+s)&= \left(\frac{2\ii s c^4}{\sqrt{V}}+\ii s r\hat c^j[\nabla_{ \hat\Theta_j},c^4\hat c^1]\right)e^{-|s|r}\Psi_\alpha(\lambda)\nonumber\\
&\quad + 2\ii s r c^4e^{-|s|r}\nabla_{e_1}\Psi_\alpha(\lambda)-2\ii s r\hat c^1\nabla^\lambda_{\hat\Theta_4}e^{-|s|r}\Psi_\alpha(\lambda)\nonumber\\
&= 2\ii s r e^{-|s|r}\left( \frac{ c^4}{ \sqrt{V}}(\frac{2 }{r }  +  \nabla_{\frac{\p}{\p r}})-\hat c^1\nabla^\lambda_{\hat\Theta_4}+ O(\frac{1}{r^2})\right)\Psi_\alpha(\lambda).
\end{align}

\begin{lemma}
$|D_{\lambda+s}\tilde \Psi_\alpha(\lambda+s)| = O( |s|r^{-3} e^{-|s|r}\|\Psi_\alpha(\lambda)\|_{L^2})$. 
\end{lemma}
\begin{proof}
Consider the tangential Dirac operator acting on sections over $\pi_k^{-1}(S_r)$, with $S_r\subset \IR^3$ the radius $r-$sphere about the origin. 
$$\mathcal{D}_S(r):= \frac{1}{2}\left(\sum_{a>1} \hat c^1\hat c^a\nabla^\lambda_{\hat\Theta_a}+ (\sum_{a>1}   \hat c^1\hat c^a\nabla^\lambda_{\hat\Theta_a})^*\right)= \sum_{a>1} \hat c^1\hat c^a\nabla^\lambda_{\hat\Theta_a}+O(r^{-2}).$$
This operator satisfies 
\begin{align}
\mathcal{D}_S^2  =\nabla^{S^*}\nabla^S+c(F_S)-\frac{1}{r}\mathcal{D}_S-\frac{1}{r}\hat c^1c^4\nabla^\lambda_{\hat\Theta_4}+ O\left(\frac{1}{r^2}\nabla\right) + O(r^{-3}),
\end{align} 
where $\nabla^S$ denotes the component of $\nabla^\lambda$ tangential to $\pi_k^{-1}(S_r)$. {Let $u$ be an eigenvector of $\nabla_S^*\nabla_S$ with eigenvalue $\alpha$. Then 
\begin{align}
\alpha \|u\|^2_{L^2(\pi_k^{-1}(S_r))}&=\|\nabla_Su\|^2_{L^2(\pi_k^{-1}(S_r))}\geq \|\nabla_{\Theta_4}^\lambda u\|^2_{L^2(\pi_k^{-1}(S_r))}\nonumber\\
&\geq c_\pi\|\Pi_1 u\|^2_{L^2(\pi_k^{-1}(S_r))},
\end{align}
where we have used  \eqref{p1ilam} for the last inequality. Hence, if $\alpha = O(r^{-2})$ and we normalize $u$ so that $\| u\|^2_{L^2(\pi_k^{-1}(S_r))}=1$,  then $\| \Pi_1u\|^2_{L^2(\pi_k^{-1}(S_r))}= O(r^{-2})$ and 
$\|\Pi_0u\|^2_{L^2(\pi_k^{-1}(S_r))}=1-O(r^{-2}).$ Thus the}  $O(r^{-2})$ eigenvalues of $\nabla_S^*\nabla_S$ have eigenspaces dominated by their $\Pi_0$ projection. {By \eqref{p0ilam},  $(\nabla_{\Theta_4}^\lambda)^*\nabla_{\Theta_4}^\lambda= O(r^{-4})$ on the image of $\Pi_0$. Moreover, the angular connection acting on the image of $\Pi_0$ differs from the Euclidean connection by $O(r^{-3})$. Hence the $O(r^{-2})$ eigenvalues (but not their multiplicities) of $\nabla_S^*\nabla_S$  lie in bands of width $O(r^{-3})$ about the Euclidean eigenvalues. }  In particular, the first three eigenvalue bands are 
$0+O(r^{-3}), \frac{2}{r^2V}+O(r^{-3}),\frac{6}{r^2V}+O(r^{-3}).$ Hence, the lowest (norm) eigenbands for $D_S$ are 
$0+O(r^{-3}), \frac{1}{r\sqrt{V}}+O(r^{-3}),-\frac{2}{r\sqrt{V}}+O(r^{-3}), \frac{2}{r\sqrt{V}}+O(r^{-3}), -\frac{3}{r\sqrt{V}}+O(r^{-3}).$ 
Let $C_\mu$ be a small circle around $\mu$. Writing the projection operator onto the $\frac{\mu}{r\sqrt{V}}+O(r^{-3})$ eigenband as 
\begin{align}q_\mu= \frac{1}{2\pi i}\int_{C_\mu}(z-r\sqrt{V}\mathcal{D}_S)^{-1}dz,
\end{align}
we have 
\begin{align}[\nabla_{\frac{\p}{\p r}},q_\mu]= \frac{1}{2\pi i}\int\displaylimits_{C_\mu}(z-r\sqrt{V}\mathcal{D}_S)^{-1}[\nabla_{\frac{\p}{\p r}},r\sqrt{V}\mathcal{D}_S](z-r\sqrt{V}\mathcal{D}_S)^{-1}dz= O(r^{-2}).
\end{align}
With these preliminaries, we use separation of variables to refine our understanding of $\Psi_\alpha(\lambda)$. Write 
$$\Psi_\alpha(\lambda) = \sum_\mu q_\mu\Psi_\alpha(\lambda)= : \sum_\mu  \Psi_\alpha^\mu,$$where the sum runs over the spectrum of the Dirac operator of the Euclidean  sphere. 
Then 
\begin{align}\label{expression}
0 &= -\sqrt{V}\hat c^1D_\lambda\Psi_\alpha(\lambda) =  ( \nabla_{\frac{\p}{\p r}}-\sqrt{V}\mathcal{D}_S+O(r^{-2}))\Psi_\alpha(\lambda)\nonumber\\
& =  \sum_\mu  \left(\nabla_{\frac{\p}{\p r}}-\frac{\mu}{r} +O(r^{-2})\right)\Psi_\alpha^\mu.
\end{align}
Here the $O(r^{-2})$ (rather than the $O(r^{-3})$ discrepancy for the eigenvalue) arises from differentiating the projection onto the eigenspace. This term is therefore an $O(r^{-2})$ map between distinct eigenbands.

Taking the inner product of \eqref{expression} with $\Psi_\alpha^\beta$, integrating from $r$ to $\infty$, and applying \cite[Prop.~26]{First} gives for $\beta> -2$:
\begin{align}\fint_{\pi_k^{-1}(S_r)}|\Psi_\alpha^\beta|^2d\sigma = O\left(  \frac{1}{\beta^2+1}r^{-6}  \|\Psi_\alpha(\lambda)\|^2_{L^2}\right),
\end{align}
and for $\beta \leq - 2$  and some fixed $R_0$, integrating from $R_0$ to $r$ gives: 
\begin{align}\fint_{\pi_k^{-1}(S_r)}|\Psi_\alpha^\beta|^2d\sigma = 
O\left(\frac{1}{\beta^2+1} r^{-6}\ln^2(r) +r^{2\beta}R_0^{-2\beta}\right) \|\Psi_\alpha(\lambda)\|^2_{L^2}).
\end{align}
In particular, the $|\Psi_\alpha^{-2}|$ mode decays like $r^{-2}$; all other modes are $O(r^{-3}\ln(r))$. Since the $O(r^{-2})$ term in \eqref{expression} exchanges bands, we have 
\begin{align} \label{quaddecay}\left|\left(\nabla_{\frac{\p}{\p r}}+\frac{2}{r}\right) \Psi_\alpha(\lambda) \right|^2  = O(r^{-8} \ln^2(r)\|\Psi_\alpha(\lambda)\|^2_{L^2}).
\end{align}
Using \eqref{quaddecay} to control the $\nabla_{\hat\Theta_1}$ term in \eqref{parks} and  \cite[Cor.~27]{First} and \eqref{p0ilam} to control the $\nabla_{\hat\Theta_4}$ term gives 
\begin{align}\label{screen}|D_{\lambda+s}\tilde \Psi_\alpha(\lambda+s)|&= 
 e^{-|s|r}  O\left(\frac{|s|}{r} |\Psi_\alpha(\lambda)|\right)= 
 e^{-|s|r} O\left(\frac{|s|}{r^{3}} \|\Psi_\alpha(\lambda)\|_{L^2}\right),
\end{align}
as claimed. In the last equality, we  use the quadratic decay of  $|\Psi_\alpha(\lambda)|$ proved in \cite[Thm.~29]{First}. 
\end{proof}
 
Passing from approximate harmonics to actual harmonics,  henceforth we work in the following unnormalized frame: 
\begin{align}
\label{exts}\Psi_\alpha(\lambda + s): = \tilde \Psi_\alpha(\lambda+s) - D_{\lambda+s}G_{D_{\lambda+s}}D_{\lambda+s}\tilde \Psi_\alpha(\lambda+s) .
\end{align}
The following proposition tells us that orthonormalization of this frame only introduces $O(|s|^{\frac{1}{2}})$ corrections. Moreover, this frame differs from the frame $\{\Pi_{D_{\lambda+s}}\Psi_\alpha(\lambda)\}_\alpha$ by at most $O\left(|s|\ln \frac{1}{|s|}\right). $
\begin{proposition}\label{ct}
\begin{align}\label{ct1}
\|\Psi_\alpha(\lambda + s)-\tilde \Psi_\alpha(\lambda +s)\|_{L^2} &= O( |s| ),\\\label{ct3}
\|\Psi_\alpha(\lambda + s)-\Pi_{D_{\lambda+s}}\Psi_\alpha(\lambda)\|_{L^2} &= O\left( |s|\ln \frac{1}{|s|} \right),\end{align}
 and \begin{align}\label{ct2}\|\Psi_\alpha(\lambda + s)-\Psi_\alpha(\lambda )\|_{L^2} = O(\sqrt{|s|}).\end{align}
\end{proposition}\label{moose}
\begin{proof}
Set 
\begin{align}\label{beta}
\beta_s:=  G_{D_{\lambda+s}}D_{\lambda+s}\tilde \Psi(\lambda+s).
\end{align}   
 By Lemma \ref{oldnag} with $K_s=\beta_s$ and $|w|=|D_{\lambda+s}\tilde\Psi(\lambda+s)|= O( |s|r^{-3}e^{-|s|r})$, 
we have for $|s|$ small, 
\begin{align}\left\| \frac{1}{r}  \eta_{ |s|,\frac{1}{2}-\delta}\beta_s\right\|^2_{L^2}= O( \delta^{-1}|s|^{2} ),
\end{align}
and
\begin{align}\left\| \eta_{ |s|,\frac{1}{2}-\delta}\nabla_{\lambda+s}\beta_s\right\|^2_{L^2}= O( \delta^{-1}|s|^{2} ).
\end{align}
In particular, 
\begin{align}\label{steve}\|  \eta_{ |s|,\frac{1}{2}-\delta}D_{\lambda+s}\beta_s\|^2_{L^2}= O(  \delta^{-1}|s|^{2} ),
\end{align}
and \eqref{ct1} follows. To obtain \eqref{ct2}, we are left to estimate
\begin{align}\|\Psi_{\alpha}(\lambda+s)-\Psi_\alpha(\lambda)\|_{L^2}^2&\leq 2\|\Psi_{\alpha}(\lambda+s)-\tilde \Psi_\alpha(\lambda+s)\|_{L^2}^2
+O(|s|)\nonumber\\
& \leq O(|s| ).
\end{align}
The relation \eqref{ct3} follows immediately from \eqref{ct1} with the use of the simple estimate 
$\|\tilde \Psi_\alpha(\lambda + s)-\Pi_{D_{\lambda+s}}\Psi_\alpha(\lambda)\|_{L^2} = O( |s|\ln\frac{1}{|s|} ).$
\end{proof}

\subsubsection{Bounded harmonic sections\texorpdfstring{ at $\lambda\in \Lambda^0$}{ }}\label{BHS}
In this subsection, we study bounded weakly harmonic spinors at  $\lambda\in \Lambda^0$, but first we need some geometric preliminaries.

Let $h$ be a section of $S\otimes {\mathcal E}\otimes \Fe_{\lambda+s}$. Let $F^{\lambda+s}$ denote the curvature of ${\mathcal E}\otimes \Fe_{\lambda+s}$. Then  
\begin{align}
[D_{\lambda+s},\nabla^{\lambda+s}]h = \theta^b\otimes c^a(F_{ab}^{\lambda+s}+R_{ab})h, \end{align}
and 
\begin{align}\label{harmgrad}
[D_{\lambda+s}^2,\nabla^{\lambda+s}]h =&  
\frac{1}{2} \theta^b\otimes c^c c^a (F_{ac;b}^{\lambda+s}+R_{ac;b})h -\theta^b\otimes (F_{ab;a}^{\lambda+s}+R_{ab;a})h \nonumber\\
&
-2\theta^b\otimes (F_{ab}^{\lambda+s}+R_{ab}) \nabla_{\Theta_a}^{\lambda+s} h.\end{align}
{(Here, as usual, the \(;a\) subscript signifies the covariant derivative by $\Theta_a$.)
} 
Since $F^{\lambda+s}$ and $R$ are Yang-Mills, the second term on the right-hand side of \eqref{harmgrad} vanishes. Since $F^{\lambda+s}$ and $R$ are anti-self-dual, the first term drops for $h$ a section of $S^+\otimes {\mathcal E}\otimes \Fe_{\lambda+s}$.   Hence, for positive chirality $h$, we have 
\begin{align}\label{stripped0}[\nabla^{\lambda+s*}\nabla^{\lambda+s},\nabla^{\lambda+s}]h = [D_{\lambda+s}^2,\nabla^{\lambda+s}]h =  -2 \theta^b\otimes (F_{ab}^{\lambda+s} +R_{ab})\nabla_{\Theta_a}^{\lambda+s}h.\end{align}

Let $W_\lambda:=\mathrm{Ker}_{L^\infty} \nabla^{\lambda *}\nabla^\lambda $ denote the space of bounded harmonic sections of ${\mathcal E}\otimes \Fe_{\lambda}$. Since $S^+$ with the spin connection is a trivial $\IC^2$ bundle with the product connection, the space of bounded harmonic sections of $S^+\otimes{\mathcal E}\otimes \Fe_{\lambda}$ is simply $\IC^2\otimes W_\lambda$. 

Now we are ready for the
\begin{proof}[Proof of Prop.~\ref{bddh2}]
\label{Proof:bddh2}
Let $\{\eta_n\}_{n=1}^\infty$ be a sequence of smooth compactly supported cutoff functions converging to $1$ pointwise and satisfying $r|\nabla\eta_n|+r^2|\nabla^2\eta_n|\leq c_\eta$, for some $c_\eta$ independent of $n$. Let $H$ be a bounded harmonic {section}. 
By \cite[Eq.~(50)]{First}, 
{which, when specialized to harmonic sections, states \(V^{-1}\Delta_{\IR^3} \frac{1}{2}\Phi(H)  = -\Phi(\nabla H)\)
}, 
we have 
\begin{align}\label{revpoinc}\int_{\IR^3}\Delta_{\IR^3}((r^2+1)^{-p}\eta_n)\frac{1}{2}\Phi(H)dv&=\int_{\IR^3}(r^2+1)^{-p}\eta_n  \Delta_{\IR^3}\frac{1}{2}\Phi(H)dv\nonumber\\
& = -\int_{\IR^3}(r^2+1)^{-p}\eta_n  \Phi(\nabla^\lambda H)Vdv.
\end{align}
For each $p<-\frac{1}{2}$, the left-hand-side of \eqref{revpoinc} is bounded independent of $n$. Taking the limit as $n\to\infty$, we have 
\begin{align}\label{revpoinc2}\int_{\IR^3}(r^2+1)^{-p}\Phi(\nabla^\lambda H)Vdv=\int_{\IR^3}-\Delta_{\IR^3} (r^2+1)^{-p}\frac{1}{2} \Phi(H)dv.
\end{align}
The right-hand-side of \eqref{revpoinc2} is uniformly bounded as $p\to \frac{1}{2}$. (The potentially problematic term behaves like $p(1-2p)\int_1^\infty r^{-2p}dr$.)  Take the limit as $p\to \frac{1}{2}$ to deduce  
\begin{align} \label{finite} \int_{\TN_k} \frac{1}{r}|\nabla^\lambda H|dv <\infty.
\end{align}
 
Applying the elliptic estimate \cite[Prop.~2]{First}  to $f:= \frac{1}{2}\Phi(\nabla^\lambda H),$ $w=O(R^{-1})$, on $B_{R}(x)$ for $|x|\geq 2R$, we deduce that
 there exists $ C_{e}>0$, independent of $R$ such that 
\begin{align}\Phi( \nabla^\lambda H)\leq  C_eR^{-2}\left\|\frac{\nabla^\lambda H}{r^{\frac{1}{2}}}\right\|_{L^2\left(B^c_{\frac{R}{2}}(o)\right)}.
\end{align}
By \cite[Lem.~14]{First},  $|\nabla^\lambda H|=O(\frac{1}{r}),$  as desired. 
By \eqref{stripped0}, 
\begin{align}| \nabla^{\lambda*}\nabla^\lambda\nabla^\lambda H| = O(r^{-2}|\nabla^\lambda H|) = O(r^{-3}).\end{align} 
Hence, taking the $L^2$ inner product of $ \nabla^{\lambda*}\nabla^\lambda\nabla^\lambda H$ with $\eta_n^2(r)r^{2p}\nabla^{\lambda}H$, and sending $n\to \infty$, we see that  
$r^p(\nabla^\lambda)^2H\in L^2$ for $p<\frac{1}{2}$.   

Let  $\zeta_v:= D_\lambda (v\otimes H)$. Then $D_\lambda\zeta_v = 0$. \cite[Prop.~29]{First} states that if $\zeta_v\in  \Ker_{L^2}(D_\lambda)$, then $\|r^2\zeta_v\|_{L^\infty(M)}<\infty $. The proof of that proposition, however, actually only uses the weaker condition that   $r^{-\frac{1}{2}}\zeta_v\in L^2$. This weaker hypothesis follows in our case from \eqref{finite}. Hence the proof of \cite[Prop.~29]{First} extends to our situation and implies $\zeta_v\in L^2$, and $\|r^2\zeta_v\|_{L^\infty(M)}<\infty $. We now extend this decay estimate to $\nabla^\lambda H$. 
Write 
\begin{align} \sum_{c=1}^4|\zeta_{I_cv}|^2 &= \sum_{a,b,c}\langle c^aI_cv,c^bI_cv\rangle \langle \nabla^\lambda_{\Theta_a}H,\nabla^\lambda_{\Theta_b}H\rangle\nonumber\\
&
= \sum_{a,b}\langle I_c^\dagger I_b^\dagger I_aI_cv, v\rangle \langle \nabla_{\Theta_a}^\lambda H,\nabla^\lambda_{\Theta_b}H\rangle 
= 4 |v|^2 |\nabla^\lambda H|^2.
\end{align}
Hence $\|r^2\nabla^\lambda H\|_{L^\infty(M)}<\infty $.

Consequently,  $|\nabla_\lambda H|$ restricted to any radial ray is $L^1$, and therefore $H$ has a well defined limit as $r\to\infty$ along any ray. If $H$ limits to $0$ at $\infty$, then it vanishes, by the maximum principle. This  implies the remaining claims. 
\end{proof}
Recall $\Pi_1 = I-\Pi_0$, where $\Pi_0$ is the projection defined in \eqref{pi0}.
\begin{coro}\label{coro27}Let $H$ be a bounded harmonic {section}. Then
 $|\Pi_1H| = O(r^{-N}),$ $\forall N$ as $r\to\infty$.
\end{coro}
\begin{proof}The assertion follows by an obvious modification of the proof of \cite[Prop.~26 and Cor.~27]{First}. 
\end{proof}
Now, to proving that $W_\lambda$ is one-dimensional:
\begin{proof}[Proof of Prop.~\ref{wdim1}]
\label{Proof:wdim1}
Let $L_\lambda$ denote the line bundle defined implicitly by the block decomposition of $\nabla^\lambda$ in \eqref{connection}. Then by \cite[Thm.~23]{First} 
applied to $\lambda\in\Lambda^0$, the curvature restricted to $L_\lambda$ is $O(r^{-3})$. Hence it is easy to construct orthonormal frames for  $L_\lambda$ with $O(r^{-2})$ connection matrices and $O(r^{-3})$ second covariant derivative on  $\TN_k\setminus K$, $K$ a large compact set. Denote this frame    $\{1_\lambda\}$. Let $\eta$ be a cutoff function supported in $\TN_k\setminus K$ and identically one outside a compact set. 
Define a section  
 $$\sigma := \eta   1_\lambda - G_{D_\lambda}\nabla^{\lambda\ast}\nabla^\lambda(\eta  1_\lambda).$$
 Then Hardy's inequality implies 
 \begin{align}\|r^{-1}G_{D_\lambda}\nabla^{\lambda\ast}\nabla^\lambda(\eta   1_\lambda)\|^2_{L^2}\leq 16\|r \nabla^{\lambda\ast}\nabla^\lambda(\eta   1_\lambda)\|^2_{L^2}.
 \end{align}
 Elliptic estimates then imply $ G_{D_\lambda}\nabla^{\lambda\ast}\nabla^\lambda(\eta   1_\lambda)$ decays pointwise as $r\to\infty$. Hence, there exist at least one element  $\sigma \in W_\lambda$. 
 
The image of $\Pi_0$ takes values in a line bundle. 
By Corollary \ref{coro27},  elements of $W_\lambda$ are asymptotically covariant constant and asymptotically take values in the image of $\Pi_0$. Thus ${\dim(W_\lambda)= 1}$.  
\end{proof}
\subsubsection{From \texorpdfstring{$L^2$}{L2}-harmonic spinors to bounded harmonic sections: the discontinuity of \texorpdfstring{$\chi_\alpha$}{chi}.}\label{flbh} 
Earlier, in  Eq.~\eqref{eq:chi}, we defined 
$ \chi_\alpha(\lambda+s): = G_{D_{\lambda+s}}\frac{c^4}{\sqrt{V}}\Psi_{\alpha}(\lambda+s).$   We now study the discontinuity of $\chi_\alpha$ at $\lambda  $ in order to prove Propositions \ref{cool} and \ref{R2q}.
Write
\begin{align}\label{eq:split}
\chi_\alpha(\lambda+s)-\chi_\alpha(\lambda-s)
  = &G_{D_{\lambda+s}}\frac{c^4}{\sqrt{V}}\Psi_{\alpha}(\lambda+s)- G_{D_{\lambda-s}}\frac{c^4}{\sqrt{V}}\Psi_{\alpha}(\lambda-s)\nonumber\\
= &\frac{1}{2}(G_{D_{\lambda+s}}+G_{D_{\lambda-s}})\frac{c^4}{\sqrt{V}}(\Psi_{\alpha}(\lambda+s)-\Psi_{\alpha}(\lambda-s))\nonumber\\
&+ \frac{1}{2}(G_{D_{\lambda+s}}-G_{D_{\lambda-s}})\frac{c^4}{\sqrt{V}}(\Psi_{\alpha}(\lambda+s)+\Psi_{\alpha}(\lambda-s)).
\end{align}
Lemma~\ref{Z} below shows that 
$\frac{1}{2}(G_{D_{\lambda+s}}-G_{D_{\lambda-s}})\frac{c^4}{\sqrt{V}}(\Psi_{\alpha}(\lambda+s)+\Psi_{\alpha}(\lambda-s))$ vanishes at $s=0$. 
The discontinuity of $\chi_\alpha$ is contained in the other term 
$\frac{1}{2}(G_{D_{\lambda+s}}+G_{D_{\lambda-s}})\frac{c^4}{\sqrt{V}}(\Psi_{\alpha}(\lambda+s)-\Psi_{\alpha}(\lambda-s))$. 
Observe 
\begin{multline}
\frac{c^4}{\sqrt{V}}(\Psi_\alpha(\lambda + s)- \Psi_\alpha(\lambda - s)) \\ 
= -  \frac{2isr\hat c^1}{\sqrt{V}}e^{-|s|r}\Psi_\alpha(\lambda) - \frac{c^4}{\sqrt{V}}(D_{\lambda+s}\beta_s -D_{\lambda-s}\beta_{-s}),
\end{multline}
where we recall $\beta_s:=  G_{D_{\lambda+s}}D_{\lambda+s}\tilde \Psi(\lambda+s) $ was defined in \eqref{beta}.
We first study the contribution of the $\beta_{\pm s}$ terms to the discontinuity. 
\begin{lemma}\label{X}For all $\delta >0$
\begin{align}\left\|G_{D_{\lambda+s}}\frac{c^4}{\sqrt{V}}D_{\lambda+ s}\beta_{ s}\right\|^2_{L^2}+\left\|G_{D_{\lambda-s}}\frac{c^4}{\sqrt{V}}D_{\lambda+ s}\beta_{ s}\right\|^2_{L^2}=O(\delta^{-3}|s|^{1-2\delta}), 
\end{align}
\begin{align}\left\|G_{D_{\lambda+s}}\frac{c^4}{\sqrt{V}}D_{\lambda+ s}\beta_{ s}\right\|_{L^1(M)}+\left\|G_{D_{\lambda-s}}\frac{c^4}{\sqrt{V}}D_{\lambda+ s}\beta_{ s}\right\|_{L^1(M)}=O(\delta^{-\frac{3}{2}}|s|^{-\frac{5}{2}}), 
\end{align}
and 
\begin{multline}
\left\| e^{-\frac{|s|r}{2}}\nabla^{\lambda+s}(e^{|s|r}G_{D_{\lambda+s}}\frac{c^4}{\sqrt{V}}D_{\lambda+  s}\beta_{   s})\right\|^2_{L^2}\\
+\left\| e^{-\frac{|s|r}{2}}\nabla^{\lambda+s}(e^{|s|r}G_{D_{\lambda-s}}\frac{c^4}{\sqrt{V}}D_{\lambda+  s}\beta_{  s})\right\|^2_{L^2}= O(\delta^{-3}|s|^{1-2\delta}).
\end{multline}
\end{lemma}
\begin{proof}
Recall equation \eqref{steve}: $\| e^{\frac{3}{4}|s|r}D_{\lambda+s}\beta_s\|^2_{L^2}= O( |s|^{2} ).$ 
 Applying Lemma \ref{oldnag} to $G_{D_{\lambda+s}}\frac{c^4}{\sqrt{V}}D_{\lambda+\epsilon s}\beta_{\epsilon s}=K_s$, $\epsilon = \pm 1$, we have 
\begin{align}&\left\|\sqrt{ \frac{1}{4r^2}  + |s|^2}\eta_{\frac{|s|}{2}, \delta-\frac{1}{2}}G_{D_{\lambda+s}}\frac{c^4}{\sqrt{V}}D_{\lambda+\epsilon s}\beta_{\epsilon s}\right\|^2_{L^2}\nonumber\\
&= O\left(\delta^{-2}\|\frac{\eta_{\frac{|s|}{2},\delta-\frac{1}{2}}D_{\lambda+\epsilon s}\beta_{\epsilon s}}{\sqrt{ \frac{1}{4r^2}  + |s|^2}} \|^2_{L^2}\right)+
O(\delta^{-3}\|r^{\frac{1}{2}+\delta}D_{\lambda+\epsilon s}\beta_{\epsilon s} \|^2_{L^2})\nonumber\\
&= O\Big(\delta^{-2}\|e^{-\frac{1}{4}|s|r}r_s^{\delta+\frac{1}{2}} \eta_{\frac{3|s|}{4},0}D_{\lambda+\epsilon s}\beta_{\epsilon s} \|^2_{L^2}\nonumber\\
&\phantom{= O(} +\delta^{-3}\|r^{\frac{1}{2}+\delta}e^{- \frac{3}{4}|s|r}\eta_{\frac{3}{4}|s|,0}D_{\lambda+\epsilon s}\beta_{\epsilon s} \|^2_{L^2}\Big)
= O(\delta^{-3}|s|^{1-2\delta}),
\end{align}
and
\begin{align}\left\| e^{-\frac{|s|r}{2}}\nabla^{\lambda+s}(e^{|s|r}G_{D_{\lambda+s}}\frac{c^4}{\sqrt{V}}D_{\lambda+\epsilon s}\beta_{\epsilon s})\right\|^2_{L^2}= O(\delta^{-3}|s|^{1-2\delta}).
\end{align}
We use Cauchy-Schwartz to estimate  the $L^1$ norm:
\begin{align}&\left\|G_{D_{\lambda+s}}\frac{c^4}{\sqrt{V}}D_{\lambda+ \epsilon s}\beta_{\epsilon s}\right\|_{L^1(M)}\nonumber\\
& \leq \left\|\frac{r^{\frac{3}{2}-\delta}e^{-\frac{|s|}{2}}}{\sqrt{ 1  + |s|^2r^2}}\right\|_{L^2}
\left\|\sqrt{ \frac{1}{4r^2}  + |s|^2}\eta_{\frac{|s|}{2}, \delta-\frac{1}{2}}G_{D_{\lambda+s}}\frac{c^4}{\sqrt{V}}D_{\lambda+ \epsilon s}\beta_{\epsilon s}\right\|_{L^2}\nonumber\\
& \leq  O(\delta^{-\frac{3}{2}}|s|^{-\frac{5}{2}}),
\end{align}
as claimed. 
\end{proof}
Hence the contribution of the $\beta_s$ terms to $\chi_{\alpha}(s)-\chi_{\alpha}(-s)$ vanishes in $H_1$ as $s\to 0$.

We turn to the dominant term. Set 
$$U_s:= -2\ii s G_{D_{\lambda+s}}\frac{\hat c^1}{\sqrt{V}}r e^{-|s|r} \Psi_\alpha(\lambda).$$
\begin{proposition}\label{Y}
$|U_s|$ is uniformly bounded as $s\to 0$,  and as $s\to 0$, 
$U_{s }$ converges in $L^2_{loc}$ to a bounded harmonic {spinor} of the form  $\ii\frac{s}{|s|}\hat c^1 \sqrt{V}r^2  \Psi_\alpha(\lambda)+O(\frac{1}{\sqrt{1+r^2}}).$ 
\end{proposition}
\begin{proof}
To estimate $U_s$, first compute
\begin{align}
&D_{\lambda+s}^2( \hat c^1 \sqrt{V}r^2 e^{-|s|r} \Psi_\alpha(\lambda))\nonumber\\
&= r^2 e^{-|s|r}\hat c^1(-D_{\lambda }+2\hat c^1\nabla_{e_1} +\frac{\ii s c^4}{\sqrt{V}}+\frac{\hat c^1}{\sqrt{V}}\left(\frac{4}{r}-|s|\right)+O(r^{-2}))^2(\sqrt{V} \Psi_\alpha(\lambda))\nonumber\\
&=  -2r|s| e^{-|s|r}\frac{\hat c^1}{\sqrt{V}}    \Psi_\alpha(\lambda)+    O\left((r^{-1}+ |s|)e^{-|s|r}|\Psi_\alpha(\lambda)|\right).
\end{align}
Here we have used \eqref{quaddecay} and $\nabla^*\nabla\hat c^1 = \frac{2}{r^2}\hat c^1+O(r^{-3})$.
Consequently
\begin{align}\label{luna}
\nabla^{\lambda+s\ast} \nabla^{\lambda+s}\left(U_s-\ii\frac{s}{|s|}\hat c^1 \sqrt{V}r^2 e^{-|s|r} \Psi_\alpha(\lambda)\right) 
=     O((r^{-1}+ |s|)e^{-|s|r}|\Psi_\alpha(\lambda)|).
\end{align} 
By Lemma \ref{oldnag}, with $K_s= U_s-\ii\frac{s}{|s|}\hat c^1 \sqrt{V}r^2 e^{-|s|r} \Psi_\alpha(\lambda)  $ and $|w|= O((r^{-3}+|s|r^{-2})e^{-|s|r})\|\Psi_\alpha(\lambda)\|_{L^2}$,  we have (after 
slight rearrangement)
\begin{multline}\label{anodyne}
\left\|e^{-\frac{1}{4}|s|r}\nabla^{\lambda+s}([e^{|s|r}U_s-\ii\frac{s}{|s|}\hat c^1 \sqrt{V}r^2 \Psi_\alpha(\lambda)])\right\|^2_{L^2}\\
+ \left\|e^{\frac{3}{4}|s|r} r_s^{-1}(U_s-\ii\frac{s}{|s|}\hat c^1 \sqrt{V}r^2 e^{-|s|r} \Psi_\alpha(\lambda)\right\|^2_{L^2} = O(\| \Psi_\alpha(\lambda)\|^2_{L^2}). 
\end{multline}
By \eqref{quaddecay}, $|\nabla_{\frac{\p}{\p r}}(r^2 \Psi) |^2  = O(r^{-4}\ln^2(r) \|\Psi_\alpha(\lambda)\|^2_{L^2}).$ 
 Hence \eqref{anodyne} implies 
\begin{align}\label{anodyne2}&\|e^{-\frac{1}{4}|s|r}\nabla_{\frac{\p}{\p r}}(e^{|s|r}U_s)\|^2_{L^2} 
 = O(\| \Psi_\alpha(\lambda)\|^2_{L^2}). 
\end{align}
Next, applying elliptic estimate obtained in  \cite[Prop.~2]{First} to the function 
$f:= \Phi(U_s-\ii\frac{s}{|s|}\hat c^1 \sqrt{V}r^2 e^{-|s|r} \Psi_\alpha(\lambda))+\frac{1}{r^2+1}$, using \eqref{luna} and \eqref{anodyne}
 yields,
\begin{align} U_s=\ii\frac{s}{|s|}\hat c^1 \sqrt{V}r^2 e^{-|s|r} \Psi_\alpha(\lambda)+ O\left(\frac{1}{\sqrt{1+r^2}}\right). 
\end{align}
In particular, $U_s$ is uniformly bounded as $s\to 0$. For every sequence $s_n\to 0$, Rellich compactness gives a subsequence converging in $L^2_{loc}$ (and in fact in any local Sobolev norm) to a bounded harmonic {section}. Any two such limits differ by an $O\left(\frac{1}{\sqrt{1+r^2}}\right)$ harmonic 
{section}, 
$Z$. By the maximum principle, $Z  = 0$. Therefore $U_s$ converges as $s\to 0\pm$ to a unique limit of the form $\ii\frac{s}{|s|}\hat c^1 \sqrt{V}r^2  \Psi_\alpha(\lambda)+ O\left(\frac{1}{\sqrt{1+r^2}}\right).$ 
\end{proof}

Set $v_s:= (G_{D_{\lambda+s}}-G_{D_{\lambda-s}})\frac{c^4}{\sqrt{V}}\frac{1}{2}(\Psi_\alpha(\lambda+s)+\Psi_\alpha(\lambda-s)).$
\begin{lemma}\label{Z}
$\lim_{s\to 0}v_s= 0$ in $H_1$.
\end{lemma} 
\begin{proof}
Compute
\begin{align*}\left(\nabla^{\lambda\ast } \nabla^\lambda + \frac{s^2}{V}\right)v_s= 2\ii s\nabla^{\lambda}_{\frac{\p}{\p\tau}}(G_{D_{\lambda+s}}+G_{D_{\lambda-s}})\frac{c^4}{\sqrt{V}}\frac{1}{2}(\Psi_\alpha(\lambda+s)+\Psi_\alpha(\lambda-s)).
\end{align*}
Set $\mu_s: = (G_{D_{\lambda+s}}+G_{D_{\lambda-s}})\frac{c^4}{\sqrt{V}}\frac{1}{2}(\Psi_\alpha(\lambda+s)+\Psi_\alpha(\lambda-s))$. 
By Lemma \ref{oldnag}, 
\begin{align}\label{xiaoma}\|  e^{\frac{3}{4}|s|r} \nabla^\lambda \mu_s\|_{L^2}^2+\left\| \frac{1}{r}e^{\frac{3}{4}|s|r} \mu_s\right\|_{L^2}^2= O(|s|^{-1}\| \Psi_\alpha(\lambda))\|^2_{L^2}).
\end{align}
We slightly modify the argument of Lemma \ref{oldnag} to estimate $v_s$. We compute
\begin{align}
&\|\nabla^{\lambda}(e^{\frac{3}{4}|s|r} v_s)\|^2_{L^2}+ \frac{7}{16}s^2\|V^{-\frac{1}{2}}e^{\frac{3}{4}|s|r} v_s\|^2_{L^2}\nonumber\\
&= 2s \mathfrak{Re}\, \langle e^{\frac{3}{4}|s|r} \mu_s, e^{\frac{3}{4}|s|r} i\nabla^{\lambda}_{\frac{\p}{\p\tau}}v_s\rangle 
\nonumber\\
&\leq   2\left(c_m+\frac{R_0^2}{4}\right)|s|\left\| \frac{1}{r}e^{\frac{3}{4}|s|r} \mu_s\right\|_{L^2}\left\| \frac{1}{r}e^{\frac{3}{4}|s|r}   v_s\right\|_{L^2} \nonumber\\
&\quad +\frac{2}{c_\pi}|s|\|  e^{\frac{3}{4}|s|r} \nabla^\lambda \mu_s\|_{L^2}\| e^{\frac{3}{4}|s|r} \nabla^\lambda v_s\|_{L^2} .
\end{align}
Applying Hardy's inequality  and \eqref{xiaoma} we deduce 
\begin{align}&\|e^{ \frac{3}{4}|s|r}\nabla_{\lambda}  v_s \|^2_{L^2}+  \| e^{\frac{3}{4}|s|r} v_s\|^2_{L^2}    = O(|s| \| \Psi_\alpha(\lambda)\|^2_{L^2}).
\end{align}
\end{proof}
\begin{proof}[Proof of Prop.~\ref{cool}]
\label{postcool}
Prop.~\ref{cool} is a direct corollary of Lemma  \ref{X}, Lemma \ref{Z}, and Proposition \ref{Y}.
\end{proof}

Recall from  Eq.~\eqref{eq:Psas}, $ q_\lambda  :\Ker_{L^2}(D_\lambda) \to S^+\otimes W_\lambda$ denotes the map: 
$$q_\lambda: \Psi(\lambda)\mapsto  \lim_{s\to 0^+}\left[G_{D_{\lambda+s}}\frac{c^4}{\sqrt{V}}\Pi_{D_{\lambda+s}}\Psi(\lambda) - G_{D_{\lambda-s}}\frac{c^4}{\sqrt{V}}\Pi_{D_{\lambda-s}}\Psi(\lambda)\right],$$
and  
$R_\lambda:\Ker_{L^2}(D_\lambda) \to S^+\otimes W_\lambda$ is the adjoint of the linear map defined by 
$$R_\lambda^\dagger(f_j^\dagger\otimes w) = -\ii D_\lambda (f_j^\dagger \otimes w).$$
 
We now prove that $R_\lambda =  2q_\lambda$:
\begin{proof}[Proof of Prop.~\ref{R2q}]
\label{proofR2q}
Recall $\nabla^{\lambda+s\ast} \nabla^{\lambda+s}=\nabla^{\lambda\ast} \nabla^{\lambda}-2\ii s\nabla^{\lambda}_{\frac{\p}{\p\tau}}+\frac{s^2}{V}.$  
Let $f\in S^+\otimes W_\lambda$ and $\Psi\in \Ker_{L^2}(D_\lambda) .$ Let $\Psi(\lambda+s)$ denote the extension of $\Psi$ to an element of $\Ker_{L^2}(D_{\lambda+s}) $ as constructed for an entire frame in \eqref{aph} and  \eqref{exts}. Let $\chi(\lambda+s):= G_{D_{\lambda+s}}\frac{c^4}{\sqrt{V}}\Psi.$ Then 
\begin{align}\nonumber
\langle   D_\lambda f ,&\Psi  \rangle_{L^2} =\lim_{s\rightarrow 0}\langle   D_\lambda f ,\Psi (\lambda+s)\rangle_{L^2}\nonumber\\
&= \lim_{s\rightarrow 0}\left\langle    f ,\left(D_{\lambda+s}-\frac{\ii  s c^4}{\sqrt{V}}\right)\Psi(\lambda +s)\right\rangle_{L^2}=\lim_{s\rightarrow 0} \ii s\left\langle f ,\frac{c^4}{\sqrt{V}}\Psi (\lambda +s)\right\rangle_{L^2}\nonumber\\
\nonumber
&= \lim_{s\rightarrow 0} \ii s\langle f ,\nabla^{\lambda+s\ast}\nabla^{\lambda+s}\chi (\lambda +s)\rangle_{L^2}\\
\nonumber
&=\lim_{s\rightarrow 0}  \left\langle \left(2s^2 \nabla^{\lambda}_{\frac{\p}{\p\tau}}+\frac{is^3}{V}\right)f ,\chi(\lambda +s)\right\rangle_{L^2}\\
&=   \lim_{s\rightarrow 0}  \frac{\ii s^3}{2}\left\langle  \frac{f}{V}   ,\chi(\lambda +s)-\chi(\lambda -s)\right\rangle_{L^2}- \lim_{s\rightarrow 0}  2s^2\langle \nabla^{\lambda}_{\frac{\p}{\p\tau}} f , \chi(\lambda +s)\rangle_{L^2}. 
\end{align}
Using Corollary \ref{coro27} and Lemma \ref{oldnag}, it is easy to check that \\
$  \lim_{s\rightarrow 0}  2s^2\langle \Pi_1\nabla^{\lambda}_{\frac{\p}{\p\tau}} f , \chi(\lambda +s)\rangle_{L^2} =0.$ On the other hand, we have  
\begin{multline} 
 \lim_{s\rightarrow 0}  2s^2\langle \Pi_0\nabla^{\lambda}_{\frac{\p}{\p\tau}} f , \chi(\lambda +s)\rangle_{L^2} 
\leq \lim_{s\rightarrow 0}  2c_ms^2\|f\|_{L^\infty}\|r^{-2}  \chi(\lambda +s)\|_{L^1(M)}\\
\leq \lim_{s\rightarrow 0}  2c_ms^2\|f\|_{L^\infty}\|r^{-1}e^{-\frac{|s|}{2}r} \|_{L^2}\|r^{-1} e^{\frac{|s|}{2}r} \chi(\lambda +s)\|_{L^2}
= O(|s|),
\end{multline} 
by Lemma \ref{oldnag},  with  $w= \frac{c^4}{\sqrt{V}}\Psi(\lambda+s).$ 
Hence  
\begin{align} 
&\langle   D_\lambda f ,\Psi \rangle_{L^2} =  \lim_{s\rightarrow 0}  \frac{\ii s^3}{2}\left\langle  \frac{f}{V}   ,\chi(\lambda +s)-\chi(\lambda -s)\right\rangle_{L^2}. 
\end{align}
Using \eqref{anodyne} to show $\lim_{s\rightarrow 0}|s|^3\|U_s-\ii\frac{s}{|s|}\hat c^1\sqrt{V}r^2e^{-|s|r}\Psi \|_{L^1(M)} = 0$ and using Lemma~\ref{X} to control the remaining terms, we have 
\begin{align}\label{almost} 
\langle    R_\lambda^\dagger f ,\Psi \rangle_{L^2} &=  \lim_{s\rightarrow 0^+}  \frac{s^3}{2}\left\langle  \frac{f}{\sqrt{V}}   , \ii\hat c^1 r^2e^{-|s|r}\Psi\right\rangle_{L^2}\nonumber\\
& =      2\langle  f  ,\chi(\lambda^+)-\chi(\lambda^-)\rangle_{ \infty } =    \langle  f  ,2q_\lambda\Psi\rangle_{ \infty }, 
\end{align}
where $\langle\cdot,\cdot\rangle_\infty$ is the inner product defined in \eqref{normalize}. 
Hence $  R_\lambda = 2q_\lambda .$  
 \end{proof}

\section{Completeness and Uniqueness}
\label{sec:CompUni}
In this section, we develop and modify Nakajima's treatment of completeness and uniqueness of the Nahm transform in \cite[Sec.~4 and 5]{NakajimaMon}. 
We first record several elementary identities that we will frequently employ. Let  $\D_{(t,b)}$ denote the bow Dirac operator determined by small bow data $(t,b)$. We recall \cite[Sec.~4.1]{Second} that for $\psi$ a section of  $S\otimes \FE\otimes \Feu^*$,
\begin{align}
\D_{(t,b)}\psi:=
\begin{pmatrix}
\left(-\nabla^0_{\frac{d}{ds}}+\ii t^0-\ii T^0 +{ \ii {\e_j}(T^j-t^j)}\right)\psi
\\ -Q_\lambda^\dagger  \psi(\lambda) \\ B_\sigma^\dagger \psi(p_{\sigma}-)-b_\sigma^\dagger \psi(p_{\sigma}+)\\  - B_\sigma^{c \dagger} \psi(p_{\sigma}+)+ b_\sigma^{c\dagger} \psi(p_{\sigma}-)
\end{pmatrix} .
\end{align}  
Here the bundles  $S$ and  $\FE$  are defined in  Sec.~\ref{subBows} and $\Feu$ is introduced in  Sec.~\ref{subBundles}.
   Let $G_{\D_{(t,b)}}$ and $G_{D_s}$ denote the Green's functions for
$ \D_{(t,b)}^\dagger \D_{(t,b)} $ and $D_{s}^{+\dagger}D_s^+$ respectively. Let $\Pi_{\D_{(t,b)}}$ and $\Pi_{D_s}$ denote the $L^2-$unitary projections onto the kernels of $\D_{(t,b)}^\dagger$ and $D_{s}^{+\dagger}$ respectively.  We will drop the operator subscripts when recording identities that are true for both cases, using in that case  $D$ to represent $\D_{(t,b)}$ or $D_s$. 
We will do the same for the respective Green's operators and the projections.

For $B$ a linear operator, elementary computations yield the following identities, {whenever the compositions are well defined}:  
\begin{align}[B,G] = -G([B,D^\dagger]D+D^\dagger[B,D])G,
\end{align}
and
\begin{align}\Pi[B,\Pi] = -\Pi[B,D]GD^\dagger.
\end{align} 

We will also have frequent need of the following identities. If $\psi$ is a section of $S^+\otimes \mathcal{E}$,  
\begin{align}\label{geomprelim} 
D^-\frac{c^4}{\sqrt{V}} \psi=  -\frac{1}{\sqrt{V}}I_c^\dagger\nabla_{\Theta_c}\psi.
\end{align}
Consequently, using the quaternionic relation \eqref{fierz}, 
\begin{align}\label{geomprelim2} \frac{1}{4}I_aI_bD^-\frac{c^4}{\sqrt{V}}I_a^\dagger =  -\frac{1}{\sqrt{V}} \nabla_{\Theta_b}.
\end{align}
The identity \eqref{geomprelim}  follows from 3 observations. First, 
if $\phi$ is a one-form, we have  
\begin{align}\{D,c(\phi)\} &=  d^*\phi+c(d\phi)-2g^{ij}\phi_i\nabla_j.
\end{align}
The second fact we utilize is that $d(V^{-\frac{1}{2}}\theta^4)   =  V^{-\frac{3}{2}}[\ast_4(dV\wedge \theta^4) -  dV\wedge \theta^4] $ 
  is anti-self-dual. Hence, its Clifford action  $Cl(d(V^{-\frac{1}{2}}\theta^4))$ annihilates $S^+$.  Finally, $d^*(V^{-\frac{1}{2}}\theta^4)=0$, and \eqref{geomprelim} follows. 
   
\subsection{Completeness}
The multi-centered Taub-NUT space  $\TN_k^\nu$ arises as a hyperk\"ahler quotient of the small bow data space  by its gauge group at level $i\nu$, with $\nu$ as given in \eqref{nanosec}. Let $\bm{\mu}_\mathfrak{s}$ denote the hyperk\"ahler moment map of the gauge action and  let $\mathcal{P}:\bm{\mu}_\mathfrak{s}^{-1}(\ii\nu)\to \text{TN}_k^\nu$ denote the quotient map.

In this subsection we construct  a map $\kappa :\mathcal{P}^*\E \to \Ker(\D^\dagger)$, the kernel bundle  of the Bow Dirac operator $\D^\dagger$ over the level set $\bm{\mu}_\mathfrak{s}^{-1}(\ii\nu).$ The map is  equivariant and descends to the quotients.  The domain of $\D_{(t,b)}^\dagger$ is 
 $L^2(S\otimes \FE\otimes\Feu^*)\oplus \mathop{\oplus}_{\lambda\in\Lambda^0}W_\lambda\otimes \Feu^*_\lambda\oplus\mathop{\oplus}_\sigma N_\sigma(\FE,\Feu^*).$ We construct this map one summand at a time. First we consider the summand $L^2(S\otimes \FE\otimes\Feu^*).$

\subsubsection{Nahm}
The Green's operator $G_{D_s}$ maps $\Gamma(S^+\otimes \E\otimes \Fe_s)$ to $\Gamma(S^+\otimes \E\otimes \Fe_s)$. For $y\in \TN_k^\nu$, and $\delta_y$ the $\delta$ function supported on $y$, $v\in \E_y,$ $G_{D_s}v\delta_y$ therefore defines a section of $\Gamma(S^+\otimes \E\otimes \Fe_s\otimes (S^+\otimes  \Fe_s)^*_y),$ which lies in the Sobolev space\footnote{
	{$H_{loc}^{-\epsilon}$ denotes  distributions which, after multiplication by a smooth compactly supported function, lie in  $H^{-\epsilon}$.
	Here $H^{-\epsilon}$ is the dual space of  $H^{\epsilon}$ - the Sobolev space of functions with ``$\epsilon$ derivatives in $L^2$''.}
	}
 $H^{-\epsilon}_{loc}$, for all $\epsilon >0$. In particular, for $f\in S^+_y$ and $\beta\in (\Fe_s)_y$, $G_{D_s}v\delta_y(f\otimes\beta) := G_{D_s}f\otimes v\otimes\beta \delta_y\in \Gamma(S^+\otimes \E\otimes \Fe_s).$ 

{\bf Convention:} Henceforth, for $A\in \Hom(X\otimes Y,Z)$, for $x\in X$, we let $Ax\in \Hom(Y,Z)$  be defined by 
$$(Ax)(z):= A(x\otimes z).$$
In particular, our notation will allow frequent mismatch between the nominal domain of a map or inner product and the given argument.

Define $\kappa_N  : \mathcal{P}^*\E  \to  \Gamma(S\otimes \mathpzc{E}\otimes \Feu^*)$ as follows. 
For $v\in (\mathcal{P}^*\E)_{(t,b)}$ and $s\in \text{Bow}
{:=\sqcup_\sigma J_\sigma}$, set 
\begin{align}
\kappa_N( v)(s):=   \Pi_{D_s}\frac{ \ii c^4}{\sqrt{V}}G_{D_s}v\delta_{\mathcal{P}((t,b))}.
\end{align}
 
\begin{lemma}
$\kappa_N( v)\in L^2(ds).$ 
\end{lemma}
\begin{proof}Let 
 $\{\psi_a\}_a$ be an $L^2-$unitary basis for $\Ker_{L^2}(D_s^-)$. Then 
\begin{align} \int|\kappa_N(v)(s)|^2ds   &=  \int\left|  \Pi_{D_s}\frac{ \ii c^4}{\sqrt{V}}G_{D_s} v\delta_{\mathcal{P}((t,b))}\right|^2 ds \nonumber\\
&= \int\sum_a\left|\langle     v ,\left(G_{D_s}\frac{ c^4}{\sqrt{V}}\psi_a\right)\mathcal{P}((t,b))\rangle  \right|^2 ds.
\end{align}
By Lemma \ref{oldnag} (and the analogous estimate for $\lambda\in \Lambda\setminus \Lambda^0$), $\|G_{D_s}\frac{ c^4}{\sqrt{V}}\psi_a\|_{L^2}\in L^1(ds)$ (since it has at worst $O(\frac{1}{|s-\lambda |^{\frac{1}{2}}})$ singularities for $\lambda \in \Lambda $ and  is bounded away from $\Lambda$).  
Applying the elliptic estimate \cite[Prop.~2]{First} to $f=|G_{D_s}\frac{ c^4}{\sqrt{V}}\psi_a|^2, $ we find 
$|G_{D_s}\frac{ c^4}{\sqrt{V}}\psi_a(\mathcal{P}((t,b)))|\in L^2(ds),$ and the claimed estimate follows. 
\end{proof} 

\begin{theo}Let
 $v\in (\mathcal{P}^*\mathcal{E})_{(t_0,b_0)}$. Then $\kappa_N(v)$ satisfies \\
$\left(\nabla_{\frac{d}{ds}}+\ii t^0_0+\ii \e_j (T^j-t^j_0)\right)\kappa_N(v) =0$. 
\end{theo}
\begin{proof}
 We break the computation into smaller pieces.  Recall the action of the unit quaternions on $S$ is contragredient to the action on $(S^+)^*$. 
\begin{align}\label{potential}&\ii \e_m (T^m-t^m_0)\kappa_N( v)(s ) = \ii  \Pi_{D_s}\left[t^m,\Pi_{D_s}\frac{ ic^4}{\sqrt{V}}  G_{D_s}\right]I_m^\dagger v \delta_{\mathcal{P}((t_0,b_0))}  \nonumber\\
& =  \Pi_{D_s}\left([t^m,D_s^+]G_{D_s}D_s^-\frac{c^4}{\sqrt{V}}     +\frac{c^4}{\sqrt{V}}G_{D_s}([t^m,\nabla^*\nabla])\right)I_m^\dagger G_{D_s} v \delta_{\mathcal{P}((t_0,b_0))}\nonumber\\
& =  -\Pi_{D_s}\left(\frac{c^m}{\sqrt{V}}G_{D_s}D_s^-\frac{c^4}{\sqrt{V}}     +\frac{c^4}{\sqrt{V}}G_{D_s}\frac{2}{\sqrt{V}}\nabla_{\Theta_m}\right)I_m^\dagger G_{D_s} v \delta_{\mathcal{P}((t_0,b_0))}\nonumber\\
& =  -\Pi_{D_s}\frac{c^4}{\sqrt{V}}G_{D_s}\left(I_m  D_s^-\frac{c^4I_m }{\sqrt{V}}     +\frac{2}{\sqrt{V}}I_m^\dagger\nabla_{\Theta_m}\right) G_{D_s} v \delta_{\mathcal{P}((t_0,b_0))}\nonumber\\
& =  -\Pi_{D_s}\frac{c^4}{\sqrt{V}}G_{D_s}\left(3\nabla_{\Theta_4}      +\frac{1}{\sqrt{V}}I_m^\dagger\nabla_{\Theta_m}\right) G_{D_s} v \delta_{\mathcal{P}((t_0,b_0))},
\end{align}
where we used \eqref{geomprelim} to obtain the last equality. 

Similarly, 
\begin{align}\label{kinetic}& \left(\nabla_{\frac{d}{ds}}+\ii t^0_0 \right)\kappa_N(  v )(s)
= \ii  \Pi_{D_s}\left[\nabla_{\frac{d}{ds}},\Pi_{D_s}\frac{c^4}{\sqrt{V}}  G_{D_s}\right]  v \delta_{\mathcal{P}((t_0,b_0))}  \nonumber\\\
& =  \ii\Pi_{D_s}\left([\nabla_{\frac{d}{ds}},D_s^+]G_{D_s}D_s^-\frac{c^4}{\sqrt{V}}     +\frac{c^4}{\sqrt{V}}G_{D_s}([\nabla_{\frac{d}{ds}},\nabla^*\nabla])\right)  G_{D_s} v \delta_{\mathcal{P}((t_0,b_0))}\nonumber\\
& =  -\Pi_{D_s}\frac{c^4}{\sqrt{V}}G_{D_s}\left(D_s^-\frac{c^4}{\sqrt{V}} -2\frac{1}{\sqrt{V}}\nabla_{\Theta_4}\right)  G_{D_s} v \delta_{\mathcal{P}((t_0,b_0))}\nonumber\\
& =   \Pi_{D_s}\frac{c^4}{\sqrt{V}}G_{D_s}\left(\frac{I_m^\dagger}{\sqrt{V}}\nabla_{\Theta_m}+3\frac{1}{\sqrt{V}}\nabla_{\Theta_4}\right)  G_{D_s} v \delta_{\mathcal{P}((t_0,b_0))}.
\end{align}
Combining \eqref{potential} and \eqref{kinetic}  yields
\begin{align}& (\nabla_{\frac{d}{ds}}+\ii t^0_0+  \ii \e_m (T^m-t^m_0))\kappa_N( v)(s) =0,
\end{align}
as claimed. 
\end{proof}

\subsubsection{Bifundamental }
 Next we turn to constructing  a map  from $\mathcal{P}^*\E$ to   $ \mathop{\oplus}_\sigma N_\sigma(\FE,\Feu^*).$ Elements of $N_\sigma(\FE,\Feu^*)$ 
 relate boundary values of solutions to $(\frac{d}{ds}+\ii t^0_0+\ii {\e_j}(T^j-t^j_0))y=0$ at $p_{\sigma-}$ to values at $p_{\sigma_+}$. In our context, that means we need to relate $\Pi_{D_{p_\sigma-}}$ to $\Pi_{D_{p_\sigma+}}$.  Since $D_{p_\sigma\pm}$ acts on $(S^+\oplus S^-)\otimes \E\otimes \Fe_{p_\sigma\pm}$, we require maps between 
$\Fe_{p_\sigma-}$ and $\Fe_{p_\sigma+}$. These are provided by the bifundamentals. In the following discussion, we will frequently replace the vector bundles $S^\pm\otimes \E\otimes \Fe_s$ by $S^\pm\otimes S^{+\ast}\otimes \E\otimes \Fe_s$. We will use the same notation $D_{s}$, $D_{s}^{\pm}$, etc. to represent the natural extension of the previously defined operators to these spaces. We  write $G_s $ for the Green's function for  $\nabla^*\nabla$ acting on sections of $\E\otimes \Fe_s$. We will simply write $G$ when the particular Green's function is clear from context.   

Recall $b_\sigma^\dagger\in \Gamma( S^+ \otimes \Hom(\Fe_{p_{\sigma-}},\Fe_{p_{\sigma+}})) =\Gamma(  S^* \otimes \Hom(\Fe_{p_{\sigma-}},\Fe_{p_{\sigma+}}))$.
We defined the bifundamental data $B_\sigma \in \Hom(  \Ker_{L^2}(D_{p_\sigma+}),  S\otimes\Ker_{L^2}( D_{p_\sigma-}^-))$, 
as well as its charge conjugate $B_\sigma^c \in \Hom( \Ker_{L^2}(D_{p_\sigma-}), S\otimes\Ker_{L^2}( D_{p_\sigma+}^-))$, by 
\begin{align}
	B_\sigma &:= \Pi_{D_{p_\sigma-}} b_\sigma \Pi_{D_{p_\sigma+}},& 
		 &\text{ and }&
	B_\sigma^c&:= \Pi_{D_{p_\sigma+}} b_\sigma^c \Pi_{D_{p_\sigma-}},
\end{align}
with the understanding that, for any covariant constant section $a$ of $\calS=(S^+)^*$, we have 
$a^\dagger B_\sigma := \Pi_{D_{p_\sigma-}} (a^\dagger b_\sigma)\Pi_{D_{p_\sigma+}}$ and 
$a^\dagger B_\sigma^c := \Pi_{D_{p_\sigma+}} (a^\dagger b_\sigma^c) \Pi_{D_{p_\sigma-}}.$ 
This convention resolves potential confusion in interpreting expressions such as  $\Pi b\Pi$, and will be used liberally below.

Recall $N_\sigma(  \mathpzc{E}, \Feu^*):= \mathpzc{E}_{p_\sigma+}\otimes\Feu^*_{p_\sigma-}\oplus \mathpzc{E}_{p_\sigma-}\otimes  \Feu^*_{p_\sigma+}$. 
Define $\kappa_B  : \mathcal{P}^*\mathcal{E} \to  N_\sigma(  \mathpzc{E}, \Feu^*)$ as follows.  For $v\in (\mathcal{P}^*\mathcal{E})_{(t_0,b_0)}$, 
\begin{align} 
\kappa_B(v)(\sigma-) 
:=  \Pi_{D_{p_\sigma+}}\frac{\ii c^4}{\sqrt{V}} \frac{  b_\sigma^\dagger }{2t_\sigma}\otimes G_{D_{p_\sigma-}} v \delta_{\mathcal{P}((t_0,b_0))} ,
\end{align}
 and
\begin{align} \kappa_B(v)(\sigma+ )  &:= \Pi_{D_{p_\sigma-}}\frac{\ii c^4}{\sqrt{V}}\frac{b_\sigma^{c\dagger}}{2t_\sigma}\otimes G_{D_{p_\sigma+}} v \delta_{\mathcal{P}((t_0,b_0))}  .
\end{align}
   In the following, we suppress the  subscripts on $D$ and $G$. 
We make use of the following relation: 
\begin{align}\label{Eq:charm}
	\Pi_{D_{p_\sigma-}} (a^\dagger b) D^+ G D^- \frac{c^4}{\sqrt{V} }\, \frac{b^\dagger}{2t}G
= \Pi_{D_{p_\sigma-}} \frac{c^4}{\sqrt{V} } \frac{(b^c)^\dagger}{2t} \otimes  a^\dagger [G,b^c]
.\end{align}
Using \eqref{Eq:commbG}, the right-hand-side is equal to 
\begin{align}\label{Eq:Gb}
\ii \Pi_{D_{p_\sigma-}} \frac{c^4}{\sqrt{V} }\, \frac{(b^c)^\dagger}{2t}\otimes  G \frac{1}{\sqrt{V} } \frac{a^\dagger \e_a^\dagger b^c}{t} \nabla_{\Theta_a} G 
,\end{align}
while, applying \eqref{ssa}, \eqref{diracadj}, and \eqref{chargeconj1}, the left-hand-side is equal to 
\begin{align}
	\Pi_{D_{p_\sigma-}} c^b  (-(a^\dagger \nabla_{\Theta_b} b)) & (G \frac{1}{\sqrt{V}}\nabla_{\Theta_a} c^a c^4 \frac{b^\dagger }{2t}  G) \\
	= \Pi_{D_{p_\sigma-}} c^4 I_b a^\dagger 
	\left(\frac{\ii}{2\sqrt{V}} \frac{b}{t}  I_b\right) &
\left(G \frac{1}{\sqrt{V}} \nabla_a (-I_a^\dagger) \frac{b^\dagger }{2t}\otimes  G\right)\\
	=-\frac{\ii}{2} \Pi_{D_{p_\sigma-}} c^4 a^\dagger \left( \frac{1}{\sqrt{V}}  \frac{b}{2t} I_b\right) &
(I_b G \frac{1}{\sqrt{V}} \frac{b^\dagger \e_a}{t}  \nabla_{\Theta_a} G)\\
	=2\frac{\ii}{2} \Pi_{D_{p_\sigma-}} \frac{c^4}{\sqrt{V} }\,\left( \frac{b}{2t}\right)^{c\dagger} &\otimes 
a^\dagger \left( G \frac{1}{\sqrt{V}} \frac{b^\dagger \e_a}{t}\right)^{c\dagger}  \nabla_a G
.\end{align}
Which, indeed, equals to \eqref{Eq:Gb}.

\begin{proposition}
\begin{align}\label{bif1}
\kappa_N(\cdot)(  p_\sigma-,\cdot) = B_\sigma \kappa_B(\cdot)(\sigma-, \cdot)+ \kappa_B(\cdot)(\sigma+, \cdot)b_{0\sigma}^c,
\end{align}
and
\begin{align}\label{bif2}\kappa_N(\cdot)( p_\sigma+,\cdot ) = B_\sigma^c \kappa_B(\cdot)(\sigma+, \cdot)+ \kappa_B(\cdot)(\sigma-, \cdot)b_{0\sigma}.
\end{align}
\end{proposition} 
\begin{proof}
	The two proofs are very similar. 
We prove \eqref{bif1}, which amounts to 
\begin{align}
	a^\dagger \Pi_{D_{p_\sigma-}} \frac{c^4}{\sqrt{V} }\, G =  \Pi_{D_{p_\sigma-}} a^\dagger  b \Pi_{D_{p_\sigma+}} \frac{c^4}{\sqrt{V} }\, \frac{b^\dagger}{2t}G   
	  +  \Pi_{D_{p_\sigma-}} \frac{c^4}{\sqrt{V} }\, \frac{(b^c)^\dagger}{2t} G  a^\dagger b^c
.\end{align}
Its right-hand-side is equal to
\begin{align}
	&\Pi_{D_{p_\sigma-}} a^\dagger  b \left((1-D^+ G D^-) \frac{c^4}{\sqrt{V} }\, \frac{b^\dagger}{2t}\otimes  G\right)\nonumber\\
	&\qquad\qquad + \Pi_{D_{p_\sigma-}} \frac{c^4}{\sqrt{V} } \frac{(b^c)^\dagger}{2t}   a^\dagger b^c G 
	+ \Pi \frac{c^4}{\sqrt{V} }\, \frac{(b^c)^\dagger}{2t} a^\dagger [G,b^c]\\
	&\quad =\Pi_{D_{p_\sigma-}} \frac{c^4}{\sqrt{V} }\, a^\dagger \frac{b\otimes b^\dagger}{2t} G 
	+ \Pi_{D_{p_\sigma-}} \frac{c^4}{\sqrt{V} }\, a^\dagger \frac{b^c\otimes (b^c)^\dagger}{2t} G \label{Terms:in} \\
	&\qquad\qquad - \Pi_{D_{p_\sigma-}} a^\dagger  b \left(D^+ G D^- \frac{c^4}{\sqrt{V} }\, \frac{b^\dagger}{2t}\otimes  G\right) 
	+ \Pi_{D_{p_\sigma-}} \frac{c^4}{\sqrt{V} } \frac{(b^c)^\dagger}{2t} a^\dagger [G,b^c]\label{Terms:out}
.\end{align}
The two terms on line \eqref{Terms:in} add up to $a^\dagger \Pi_{D_{p_\sigma-}} \frac{c^4}{\sqrt{V} }\, G$ and the two terms on line \eqref{Terms:out} cancel thanks to \eqref{Eq:charm}. 
\end{proof}

\subsubsection{Fundamental}
For $\lambda\in \Lambda^0$ and $v\in (\mathcal{P}^*\mathcal{E})_{(t,b)}$, define 
$$\kappa_F( v)(\lambda):= \frac{1}{2\ii}\langle v,\phi_\lambda(\mathcal{P}((t,b)))\rangle\phi_\lambda,$$ where $\phi_\lambda$ is a unit vector in $W_\lambda$. 
Set $$  \kappa:= (\kappa_N,\kappa_B,\kappa_F).$$ 
\begin{proposition}\label{isodown}
For $v\in (\mathcal{P}^*\mathcal{E})_{(t,b)}$,  
{we have} 
$\D^\dagger_{(t,b)}  \kappa(v) =0.$   
\end{proposition}
\begin{proof}
We are left to show the desired behavior at points $\lambda\in \Lambda^0$ :  
\begin{align}\label{desired}\lim_{s\to 0^+}(\kappa_N(v)(\lambda+s)-\kappa_N(v)(\lambda-s)) =  Q_\lambda(\kappa_{W}(v) (\lambda)).
\end{align}  
We expand the left-hand side of \eqref{desired}
\begin{align}\label{d}
&\lim_{s\to 0^+}  \psi_a(\lambda+s)\left\langle   v  ,G\frac{\ii c^4}{\sqrt{V}}\psi_a(\lambda+s)(\mathcal{P}((t,b)))\right\rangle \nonumber\\
&\quad -\lim_{s\to 0^+}  \psi_a(\lambda-s)\left\langle   v ,G \frac{ic^4}{\sqrt{V}}\psi_a(\lambda-s)(\mathcal{P}((t,b)))\right\rangle \nonumber\\
&= -\ii  \psi_a(\lambda )\langle   v ,  q_\lambda(\psi_a(\lambda ))(\mathcal{P}((t,b)))\rangle \nonumber\\
&\quad +\lim_{s\to 0^+}  (\psi_a(\lambda+s)-\psi_a(\lambda))\left\langle  v ,G\frac{\ii c^4}{\sqrt{V}}\psi_a(\lambda+s)(\mathcal{P}((t,b)))\right\rangle \nonumber\\
&\quad -\lim_{s\to 0^+} (\psi_a(\lambda-s)-\psi_a(\lambda))\left\langle   v ,G\frac{\ii c^4}{\sqrt{V}}\psi_a(\lambda-s)(\mathcal{P}((t,b)))\right\rangle \nonumber\\
&= -\ii  \psi_a(\lambda )\langle  v  ,  q_\lambda(\psi_a(\lambda ))(\mathcal{P}((t,b)))\rangle ,
\end{align} 
by the Lebesgue dominated convergence theorem. 
On the other hand, 
\begin{align} Q_\lambda(\kappa_{W}(v))  &=   2q_\lambda^\dagger( \kappa_F(v)(\lambda ))\nonumber\\
&= -\ii   \langle v ,\phi_\lambda(\mathcal{P}((t,b)))\rangle q_\lambda^\dagger( \phi_\lambda )\nonumber\\
&= -\ii  \psi_a \langle   v ,q_\lambda(\psi_a(\lambda))(\mathcal{P}((t,b)))\rangle,
\end{align}
 yielding the claimed result. 
\end{proof}
\subsubsection{Up after Down}
The preceding proposition shows that $ \kappa $ defines a map  
$$\kappa :(\mathcal{E},A)\to Up\circ Down (\mathcal{E},A).$$ 
\begin{theo}\label{template} Let 
$(\mathcal{E}',A'):=Up\circ Down (\mathcal{E},A)$. Then $ \kappa :(\mathcal{E},A)\to (\mathcal{E}',A')$  is covariant constant with respect to the induced connection on $\Hom(\mathcal{E},\mathcal{E}')$. If $\mathcal{E}$ is irreducible then $Up\circ Down$ is bijective and the metrics on $\mathcal{E}$ and $\mathcal{E}'$ agree up to scale. 
\end{theo}
\begin{proof} 
 Let $v$ be a differentiable section of $\mathcal{E}$.  Let $X\in T_{z_0}\TN_k^\nu$. We are to show 
\begin{align}\label{goal1}\nabla_X \kappa(v) =  \kappa(\nabla_X v).
\end{align}
 By linearity, it suffices to consider $X=\frac{1}{\sqrt{V}}\Theta_a$. By the definition of the connection induced by the Up transform, \eqref{goal1} becomes 
\begin{align}\label{goal2}\frac{1}{V}W_a \kappa(v) -  \kappa(\nabla_{\frac{1}{\sqrt{V}}\Theta_a} v)\in \text{Im}\D,
\end{align}
where $W_a$ denotes the horizontal lift of $V^{\frac{1}{2}}\Theta_a$ to the moment map level set (see \cite[Eq.~(54)]{Second}). 
 We first compute  
\begin{align}
&\frac{1}{V}W_a \kappa_N(v )(s) -  \kappa_N( \nabla_{\frac{1}{\sqrt{V}}\Theta_a} v)(s)= - \Pi_{D_s}\frac{ \ii c^4}{\sqrt{V}}G\nabla_{\frac{1}{\sqrt{V}}\Theta_a}  v\delta_{\mathcal{P}((t,b))} \nonumber\\
&=  \frac{1}{4}\Pi_{D_s}\frac{ ic^4}{\sqrt{V}}GI_b I_aD_s^- \frac{c^4}{\sqrt{V}} I_b^\dagger v\delta_{\mathcal{P}((t,b))} \nonumber\\
&=  \frac{1}{4} \Pi_{D_s}\left(\frac{ ic^4}{\sqrt{V}}G  D_s^- \frac{c^4}{\sqrt{V}} I_a  
+ \frac{ \ii c^m}{\sqrt{V}}G D_s^- \frac{c^4}{\sqrt{V}} I_aI_m^\dagger \right) v\delta_{\mathcal{P}((t,b))} \nonumber\\
&
=  \frac{1}{4}\Big(\Pi_{D_s}\frac{ d}{ds}D_s^+G  D_s^- \frac{c^4}{\sqrt{V}} - \e_m  \Pi_{D_s}t^mD_s^+G D_s^- \frac{ic^4}{\sqrt{V}} \Big)I_a  v\delta_{\mathcal{P}((t,b))}\nonumber\\
&=  -\frac{1}{4} \Pi_{D_s}\frac{ d}{ds}\Pi_{D_s}   \frac{c^4}{\sqrt{V}} I_a   v\delta_{\mathcal{P}((t,b))}   -\frac{1}{4}\e_m \Pi_{D_s}t^m(I-\Pi_{D_s}) \frac{\ii c^4}{\sqrt{V}} I_a   v \delta_{\mathcal{P}((t,b))} \nonumber\\
&=  \left(- \nabla_{\frac{ d}{ds}}+ \ii I_m(T^m-t^m_0)\right) \Pi_{D_s}  \frac{c^4}{4\sqrt{V}} I_a v\delta_{\mathcal{P}((t,b))}.
\end{align}
Set $$y_a:=  \Pi_{D_s}  \frac{c^4}{4\sqrt{V}} I_a v\delta_{\mathcal{P}((t,b))} .$$
We are left to show that for all $\sigma$, and for all $\lambda\in \Lambda^0$, 
\begin{align}\label{vitb1} 
(B_\sigma^\dagger y_a (p_\sigma-) - b_\sigma^{\dagger}y_a (p_\sigma+))& = \frac{1}{V}W_a\kappa_B(v)(\sigma^-)- \kappa_B(\nabla_{\frac{1}{\sqrt{V}}\Theta_a} v)(\sigma^-)\nonumber\\
& - \kappa_B( v)(\sigma^-,\nabla_{\frac{1}{\sqrt{V}}\Theta_a}),\\
 \label{vitb2}-(B_\sigma^{c\dagger}y_a (p_\sigma+)-b_\sigma^{ c\dagger}y_a (p_\sigma-)) &= \frac{1}{V}W_a\kappa_B(v)(\sigma^+,\cdot) - \kappa_B(\nabla_{\frac{1}{\sqrt{V}}\Theta_a} v)(\sigma^+)\nonumber\\
&
- \kappa_B( v)(\sigma^+,\nabla_{\frac{1}{\sqrt{V}}\Theta_a}\beta),\text{  and}\\
 \label{choir}-Q_\lambda^\dagger y_a (\lambda) &=\frac{1}{V}W_a\kappa_F(v)(\lambda) - \kappa_F(\nabla_{\frac{1}{\sqrt{V}}\Theta_a} v)(\lambda).
\end{align}
We will demonstrate \eqref{vitb1} and \eqref{choir}. The proof of \eqref{vitb2} is similar to that of \eqref{vitb1}.

Consider  \eqref{vitb1}, and compute, with $\{f_a\}_a$ a covariant constant unitary frame of $S$: 
\begin{align}
& B_\sigma^\dagger y_a (p_\sigma-)-b_\sigma^{\dagger}y_a (p_\sigma+)\nonumber\\
& =\begin{aligned}[t]  \Pi_{D_{p_\sigma+}}b_\sigma^\dagger(f_j)  \Pi_{D_{p_\sigma-}}\frac{c^4I_a f_j^\dagger}{4\sqrt{V}}  \otimes v \delta_{\mathcal{P}((t,b))} \\
{} - b_\sigma^{ \dagger}(f_j)\Pi_{D_{p_\sigma+}}  \frac{c^4I_a f_j^\dagger}{4\sqrt{V}} \otimes v\delta_{\mathcal{P}((t,b))}
\end{aligned} 
\nonumber\\
& =\begin{aligned}[t]  \Pi_{D_{p_\sigma+}}b_\sigma^\dagger(f_j) (I- D_{p_\sigma-}^+GD_{p_\sigma-}^-)\frac{c^4I_a f_j^\dagger}{4\sqrt{V}}  \otimes v \delta_{\mathcal{P}((t,b))} \\
- \Pi_{D_{p_\sigma+}}  \frac{c^4I_a b_\sigma^{ \dagger}}{4\sqrt{V}} \otimes v \delta_{\mathcal{P}((t,b))} 
\end{aligned} 
\nonumber\\
& =   -\Pi_{D_{p_\sigma+}}\frac{c^4I_b}{\sqrt{V}}\frac{\ii b_\sigma^\dagger I_b(f_j)}{2t_\sigma }   G \frac{1}{\sqrt{V}}I_c^\dagger\nabla_{\Theta_c}\frac{ I_a f_j^\dagger}{4 }  \otimes v\delta_{\mathcal{P}((t,b))} 
  \nonumber\\
& =   -\Pi_{D_{p_\sigma+}}\frac{c^4 }{\sqrt{V}}\frac{\ii b_\sigma^\dagger  (f_j)}{2t_\sigma }   G \frac{1}{\sqrt{V}}\nabla_{\Theta_a}  f_j^\dagger \otimes v\delta_{\mathcal{P}((t,b))} .
\end{align}
On the other hand, 
\begin{align}&\frac{1}{V}W_a\kappa_B(v)(\sigma^-) - \kappa_B(\nabla_{\frac{1}{\sqrt{V}}\Theta_a} v)(\sigma^-)
\nonumber\\
&= -\Pi_{D_{p_\sigma+}}\frac{\ii c^4}{\sqrt{V}}\frac{b_\sigma^\dagger(f_j)}{2t_\sigma} G\frac{1}{\sqrt{V}}\nabla_{\Theta_a}  f_j^\dagger\otimes v \delta_{\mathcal{P}((t,b))}.
\end{align}
Finally we check the $\kappa_F$ term: 
\begin{align}& \langle Q_\lambda^\dagger y_a,\phi_\lambda\rangle_\infty 
=  \left\langle \Pi_{D_\lambda}  \frac{c^4}{4\sqrt{V}} I_a  v\delta_{\mathcal{P}((t,b))} ,Q_\lambda\phi_\lambda\right\rangle\nonumber\\
&=  \ii\left\langle    \frac{c^4}{4\sqrt{V}} I_a  v\delta_{\mathcal{P}((t,b))} ,   c^\mu f_j^\dagger \nabla_{\Theta_\mu} \phi_\lambda\right\rangle\nonumber\\
&=   \ii \mathrm{tr}_S\, I_a^\dagger I_\mu\left\langle      \frac{1}{4\sqrt{V}} v\delta_{\mathcal{P}((t,b))} ,   \nabla_{\Theta_\mu} \phi_\lambda\right\rangle\nonumber\\
&=   \ii  \frac{1}{2 } \left\langle        v\delta_{\mathcal{P}((t,b))} ,    \frac{1}{ \sqrt{V}}\nabla_{\Theta_a} \phi_\lambda\right\rangle.
\end{align}
On the other hand, we have  
\begin{align} 
&V^{-1}W_a\kappa_F(v)(\lambda)-\kappa_F(\nabla_{V^{-1/2}\Theta_a}v)(\lambda)\nonumber\\
& = \frac{1}{2\ii}\left\langle v, \frac{1}{\sqrt{V}}\nabla_{\Theta_a}\phi_\lambda(\mathcal{P}((t,b)))\right\rangle\phi_\lambda
=-Q_\lambda^\dagger y_a,
\end{align}
and   we have $\frac{1}{V}W_a \kappa(\frac{1}{ \sqrt{V}}\nabla_{\Theta_a}v) -  \kappa(\nabla_{\frac{1}{\sqrt{V}}\Theta_a} v) = \D_{(t,b)} y_a,$ as desired. Hence $ \kappa$ is covariant constant. 
 
Consequently  $  \kappa^\dagger\circ \kappa$ and $  \kappa\circ \kappa^\dagger$ are  covariant constant sections of $\End(\mathcal{E})$ and $\End(\mathcal{E}')$ respectively. If $\mathcal{E} $, respectively $\mathcal{E}'$ is irreducible, then $  \kappa^\dagger\circ \kappa=\lambda I_{\mathcal{E}},$ respectively $  \kappa\circ\kappa^\dagger=\lambda I_{\mathcal{E}'}$ for some scalar $\lambda\geq 0$. If $\E$ is reducible, then $\kappa^\dagger\circ \kappa = \sum_j\lambda_jP_j$ for some hermitian commuting projection operators $P_j$, with $\sum_j P_j =1$. We claim each $\lambda_j>0$. To see this, suppose some $\lambda_j=0$. Then $ \kappa(P_jv) = 0,\forall v$. This implies $\Pi_{P_jD_sP_j}$ annihilates the image of $\frac{c^4}{\sqrt{V}}G$. The cokernel of this operator is zero and therefore $\Pi_{P_jD_sP_j}$ must vanish for all $s$. This contradicts the index computation \cite[Thm.~44]{First}. Hence $\lambda_j\not = 0$ and $\kappa P_j$ is an isometry onto its image, up to scale.  By \cite[Thm.~1]{Second}, the rank of $\E'$ is $|\Lambda|$. Using \cite[Thm.~44]{First} to compute the rank jumps in $\mathpzc{E}$, we see that  $|\Lambda| = \text{rank}\, \E$. Hence $\kappa$ is surjective and $Up\circ Down$ is bijective.
\end{proof}
 \begin{coro}[Completeness]
	\label{halfway}The Down transform is injective and the $Up$ transformation is surjective. 
\end{coro}
\begin{proof} $ Up\circ  Down $ is  bijective. Therefore $Up$ is surjective and $Down$ is injective. 
\end{proof} 

\subsection{Uniqueness} 

The proof of uniqueness is similar to that of completeness.
Given a bow solution $\mathcal{B}=(T,Q,B)$, with bow Dirac operator family  $\D_{(t,b)}$ and an element $s$ of the bow, we now construct a   map $\Phi_s :   \mathpzc{E}_s\to \Gamma(\mathcal{P}^*\left(S^-\otimes \mathcal{E}\otimes \Fe_s\right))$ as follows.   For $v\in \mathpzc{E}_s$,   set
\begin{align}\Phi_s(v)((t,b)):= c^4\Pi_{\D_{ (t,b)}}[
\nabla^s_{\Theta_4}, \D_{ (t,b)}]G_{\D_{ (t,b)}}  v \delta_s.
\end{align}
\begin{theo}
$D_{s}\Phi_s(v)((t,b)) = 0.$
\end{theo}
\begin{proof}
The proof is a direct computation, using identities developed in the appendix.  We recall that the connection on the index bundle of the bow is given by $\nabla_{\Theta_a}  = \Pi_{\D_{ (t,b)}}\Theta_a^H\Pi_{\D_{ (t,b)}}$, and the connection on $\Fe_s$ is given simply by the action of $\Theta_a^H$ on the corresponding equivariant section. Here $\Theta_a^H$ denotes the horizontal lift of $\Theta_a$ to the given level set 
$\bm{\mu}^{-1}(\ii \nu )$ of the small bow moment map. 
First we observe that for $F$ a section of $S^+\otimes \E\otimes \Fe_s$, 
\begin{align}\label{red1}D_{s}V^{-\frac{1}{2}}c^4 F &= V^{-\frac{1}{2}}c^ac^4  \nabla_{\Theta_a}^sF + Cl(d(V^{-\frac{1}{2}}\theta^4))  F\nonumber\\
&= -V^{-\frac{1}{2}}I_a^\dagger  \nabla_{\Theta_a}^sF,
\end{align}
where the last equality follows from observing that $d(V^{-\frac{1}{2}}\theta^4)$ is anti-self-dual, and therefore its Clifford action $Cl(d(V^{-\frac{1}{2}}\theta^4))$ annihilates sections of $S^+$, and $I_a^\dagger = c^ac^4$ when acting on sections of $S^+$. 
Applying \eqref{red1}, we have 
\begin{align}
&D_s\Phi_s(v)((t,b))\nonumber\\
&= -V^{-\frac{1}{2}}  \nabla_{\Theta_a}^s\left(\Pi_{\D_{ (t,b)}}[V^{\frac{1}{2}} \Theta_4^H, \D_{ (t,b)}]G_{\D_{ (t,b)}}\e_a v \delta_s\right)\nonumber\\
&=\label{fr1} -V^{-\frac{1}{2}} [\nabla^s_{\Theta_a},\Pi_{\D_{ (t,b)}}][V^{\frac{1}{2}} \Theta_4^H, \D_{ (t,b)}]G_{\D_{ (t,b)}}\e_a  v \delta_s \\
&\label{fr2} 
\quad -V^{-\frac{1}{2}}   \Pi_{\D_{ (t,b)}}[ \Theta_a^H,[V^{\frac{1}{2}}\Theta_4^H, \D_{ (t,b)}]]G_{\D_{ (t,b)}}\e_a  v \delta_s \\
&\label{fr3}
\quad  -V^{-\frac{1}{2}}  \Pi_{\D_{ (t,b)}}[V^{\frac{1}{2}}\Theta_4^H, \D_{ (t,b)}][ \Theta_a^H,G_{\D_{ (t,b)}}]\e_a  v \delta_s.
\end{align}
To simplify, we use (the adjoint of) \eqref{ssa} to compute that 
\begin{align}[ \Theta_4^H, \D_{ (t,b)}] = \left(\begin{array}{c}-\frac{\ii}{ \sqrt{V}} \\
0\\
  - \frac{\ii}{2t\sqrt{V}}b_\sigma^{ \dagger} ev_{p_\sigma+}\\
- \frac{\ii}{2t\sqrt{V}}b_\sigma^{c\dagger} ev_{p_\sigma-}\end{array}\right).
\end{align}
Consequently 
\begin{align}[ \Theta_a^H,[V^{\frac{1}{2}}\Theta_4^H, \D_{ (t,b)}]]G_{\D_{ (t,b)}}\e_a = 
 \left(\begin{array}{c}0 \\
0\\
  - \nabla_{\Theta_a}\frac{\ii}{2t\sqrt{V}}b_\sigma^{ \dagger}\e_a ev_{p_\sigma+}\\
- \nabla_{\Theta_a}\frac{\ii}{2t\sqrt{V}}b_\sigma^{c\dagger} \e_a ev_{p_\sigma-}\end{array}\right)G_{\D_{ (t,b)}} =0,
\end{align}
by \eqref{diracadj} and the corresponding equality for $b^{c\dagger}$. Hence the summand \eqref{fr2} vanishes. 

Using \eqref{ssa} and its adjoint, we compute 
\begin{align}\label{bike1}[\Theta_a^H,\D_{ (t,b)}] = [\Theta_4^H,\D_{ (t,b)}]\e_a\end{align}
and 
\begin{align}\label{bike2}[\Theta_a^H,\D_{ (t,b)}^\dagger] = \e_a^\dagger[\Theta_4^H,\D_{ (t,b)}^\dagger].\end{align}
Formal manipulations and \eqref{bike1} and \eqref{bike2} yield 
\begin{align}\label{pablum1}[\Theta_a^H,G_{\D_{ (t,b)}}] &=  -G_{\D_{ (t,b)}}(\e_a^\dagger[\Theta_4^H,\D_{ (t,b)}^\dagger]\D_{ (t,b)} +\D_{ (t,b)}^\dagger[\Theta_4^H,\D_{ (t,b)}]\e_a)G_{\D_{ (t,b)}},\end{align}
and
\begin{multline}\label{pablum2}
[\Theta_a^H,\Pi_{\D_{ (t,b)}}] 
=   -\Pi_{\D_{ (t,b)}}[\Theta_4^H,\D_{ (t,b)}]\e_aG_{\D_{ (t,b)}}\D_{ (t,b)}^\dagger \\
-\D_{ (t,b)} G_{\D_{ (t,b)}} \e_a^\dagger[\Theta_4^H,\D_{ (t,b)}^\dagger]\Pi_{\D_{ (t,b)}}
\end{multline}

Substituting \eqref{pablum1} and \eqref{pablum2} into \eqref{fr1} and \eqref{fr3} yields 
\begin{align}
&D_{s}\Phi_s(v)((t,b))\nonumber\\
&=   V^{-\frac{1}{2}}   \Pi_{\D_{ (t,b)}}\left([ \Theta_4^H,\D_{ (t,b)}]\e_aG_{\D_{ (t,b)}}\D_{ (t,b)}^\dagger[V^{\frac{1}{2}} \Theta_4^H, \D_{ (t,b)}]G_{\D_{ (t,b)}}\e_a\right.  \nonumber\\
&\quad  + [V^{\frac{1}{2}} \Theta_4^H, \D_{ (t,b)}]G_{\D_{ (t,b)}} \e_a^\dagger[\Theta_4^H,\D_{ (t,b)}^\dagger]\D_{ (t,b)}   G_{\D_{ (t,b)}}\e_a \nonumber\\
&\quad +\left.[V^{\frac{1}{2}} \Theta_4^H, \D_{ (t,b)}]G_{\D_{ (t,b)}}\D_{ (t,b)}^\dagger[\Theta_4^H,\D_{ (t,b)}]\e_aG_{\D_{ (t,b)}}\e_a\right)  v \delta_s\nonumber\\
&=   V^{-\frac{1}{2}}  \Pi_{\D_{ (t,b)}}\left(-[\Theta_4^H,\D_{ (t,b)}] G_{\D_{ (t,b)}}\e_a^\dagger\D_{ (t,b)}^\dagger[V^{\frac{1}{2}} \Theta_4^H, \D_{ (t,b)}]G_{\D_{ (t,b)}}\e_a\right.  \nonumber\\
&\quad +[\Theta_4^H,\D_{ (t,b)}] G_{\D_{ (t,b)}}2\D_{ (t,b)}^\dagger[V^{\frac{1}{2}} \Theta_4^H, \D_{ (t,b)}]G_{\D_{ (t,b)}}   \nonumber\\
&\quad   + [V^{\frac{1}{2}}\Theta_4^H, \D_{ (t,b)}]G_{\D_{ (t,b)}} \e_a^\dagger[\Theta_4^H,\D_{ (t,b)}^\dagger]\D_{ (t,b)}   G_{\D_{ (t,b)}}\e_a \nonumber\\
&\quad -2\left.[V^{\frac{1}{2}}\Theta_4^H, \D_{ (t,b)}]G_{\D_{ (t,b)}}\D_{ (t,b)}^\dagger[\Theta_4^H,\D_{ (t,b)}] G_{\D_{ (t,b)}} \right) v \delta_s\nonumber\\
&=   V^{-\frac{1}{2}}   \Pi_{\D_{ (t,b)}}[\Theta_4^H,\D_{ (t,b)}] G_{\D_{ (t,b)}}I_a^\dagger\left(-\D_{ (t,b)}^\dagger[V^{\frac{1}{2}} \Theta_4^H, \D_{ (t,b)}] \right.  \nonumber\\
&\quad   \left.+    [\Theta_4^H,\D_{ (t,b)}^\dagger]\D_{ (t,b)}   \right)\e_a  G_{\D_{ (t,b)}}v \delta_s\nonumber\\
&=   -\frac{2}{V^{\frac{1}{2}}}  \Pi_{\D_{ (t,b)}}[\Theta_4^H,\D_{ (t,b)}] G_{\D_{ (t,b)}}\e_a^\dagger \e_k(T^k-t^k)  I_a  G_{\D_{ (t,b)}}  v \delta_s\nonumber\\
&=0,
\end{align}
as claimed.
\end{proof}
\begin{proposition}
$\Phi_s( (T^m-t^m)v) \in \mathfrak{Im}\,(D_{s})$, and $\frac{\p\Phi_s}{\p s}( v ) = \Phi_s(  \ii(T^0-t^0)v).$
\end{proposition}
\begin{proof}
 
\begin{align}&D_{s}  \left(  \Pi_{\D_{ (t,b)}}\e_m^\dagger   v  \delta_s \right)\nonumber\\
&
= \Pi_{\D_{ (t,b)}} c^a   [ \Theta_a^H,\Pi_{\D_{ (t,b)}}]\e_m^\dagger  v  \delta_s 
\nonumber\\
&
= c^4 I_a   \Pi_{\D_{ (t,b)}}[\Theta_4^H,\D_{ (t,b)}]\e_aG_{\D_{ (t,b)}}\D_{ (t,b)}^\dagger \e_m^\dagger   v \delta_s \nonumber\\
&
=   c^4 \Pi_{\D_{ (t,b)}}[\Theta_4^H,\D_{ (t,b)}]G_{\D_{ (t,b)}}\e_a\D_{ (t,b)}^\dagger \e_m^\dagger  \e_a^\dagger  v \delta_s \nonumber\\
&
=   c^4  \Pi_{\D_{ (t,b)}}[\Theta_4^H,\D_{ (t,b)}]G_{\D_{ (t,b)}}4\ii(T^m-t^m)   v  \delta_s \nonumber\\
&
=   4\ii\Phi_s( (T^m-t^m)v).
\end{align}
Similarly, if $v= v(s)$ for some section $v(\cdot)$, then 
\begin{align}&D_{s}    \left(f_j^\dagger\otimes \Pi_{\D_{ (t,b)}} f_j\otimes v \delta_s \right)\nonumber\\
&
=   c^4 f_j^\dagger\otimes  \Pi_{\D_{ (t,b)}}[\Theta_4^H,\D_{ (t,b)}]G_{\D_{ (t,b)}}4 (\frac{d}{ds}+\ii(T^0-t^0))  f_j\otimes v  \delta_s \nonumber\\
&
=  - 4 \frac{\p\Phi_s}{\p s}(  v)  +4\ii \Phi_s( (T^0-t^0)v).
\end{align}
\end{proof}
Let $ \mathpzc{E} $ be the vector bundle of the large bow representation, and let $ \mathpzc{E}' $ denote $Down\circ Up  ( \mathpzc{E}). $ Arguing as in Proposition \eqref{template} and Corollary \ref{halfway}, we obtain the following corollaries. 
\begin{coro}\label{coro25}
$\Phi$ is a covariant constant element of $\Hom(  \mathpzc{E}, \mathpzc{E}' ).$ $Down\circ Up$ is bijective, and if $ \mathpzc{E} $ is irreducible, then $\Phi$ is an isometry up to scale. 
\end{coro}
\begin{coro}\label{rest}$Up $ is injective and $Down$ is surjectice.
\end{coro}
\begin{proof}
$   Down\circ Up$ is {bijective.}
\end{proof}
\begin{theo}The $Up$ and $Down$ transforms are bijective. 
\end{theo}
\begin{proof} This follows from Corollaries \ref{halfway} and \ref{rest}.
\end{proof}

\section{Moduli}
\label{sec:Moduli}
In this section, following \cite{FU}, we show that the moduli space of irreducible instantons on $\TN_k^\nu$ is a smooth  (albeit  not necessarily complete)  manifold. 
We first recall some features of instantons on $\TN_k^\nu$. Let $(\mathcal{E},A)$ be a Hermitian vector bundle over $\TN_k^\nu$, with finite energy anti-self-dual connection $A$. Outside a  compact set, $\mathcal{E}$ splits into a sum of holonomy eigen bundles $\mathcal{E} = \oplus_a\mathcal{E}_a$, where the holonomy of $\mathcal{E}_a$ around the circle fiber is $e^{2\pi i\frac{\lambda_a+O(r^{-1})}{l}}$. 
We will say such a connection has asymptotic holonomy $\vec{\lambda}$, where $\vec{\lambda }$ is the vector with components $\lambda_a$ with  $0\leq\lambda_1<\lambda_2\cdots<\lambda_r<l.$ 
We further restrict to the case where each eigenbundle has rank 1. 
So, each $\lambda_a$ has multiplicity 1.  
Connections satisfying these hypotheses are said to have {\em generic asymptotic holonomy}. 
In  \cite[Thm.~23]{First} we show that every instanton connection with generic asymptotic holonomy has the form 
\begin{align}\label{model}A = \oplus_a\left(-\ii(\lambda_a+\frac{m_a}{2r})\frac{d\tau+\omega}{V}+\pi_k^*\eta_a\right)+ O(r^{-2}),
\end{align}
where in \cite[Sec.~6.3]{Second},
we show that $\eta_a$ can be chosen to be the pullback of the standard connection on the  Hopf bundle of degree $m_a$. 

Given a Hermitian vector bundle $\mathcal{V}$ equipped with a connection, we let $L^2_{j}(\mathcal{V})$ denote the closure of the smooth compactly supported sections of $V$ with respect to the norm $\|f\|^2_{L^2_j}:= \sum_{0\leq i\leq j} \|\nabla^i f\|^2_{L^2}.$ Introduce the weighted Sobolev spaces $L^2_{j,w}(\mathcal{V})$ which are the closure of the smooth compactly supported sections of $\mathcal{V}$ with respect to the norm $\|f\|^2_{L^2_{j,w}}:= \sum_{1\leq i\leq j} \| (1+r^{2i})^{\frac{1}{2}}\nabla^i f\|^2_{L^2} + L\|  f\|^2_{L^2}.$  Here $L$ is a positive constant which will be constrained later. 

 Fix a smooth  irreducible connection $A_0$ of the form \eqref{model}, with $r^2 F_{A_0}\in L^\infty$, and $F_{A_0}^+\in L^2_{1,w}$, and let  $\mathcal{A}$ denote the space of connections of the form $A= A_0 +a$, where $a\in L^2_{2,w}(T^*\text{TN}_k\otimes ad(\mathcal{E}))$.  Let $\mathcal{G}$ denote the  unitary automorphisms of $E$ which satisfy $g-I\in L^2_{3,w}$. By an obvious extension of the Sobolev multiplication theorem, $\mathcal{G}$ is a   Hilbert Lie group, with Lie algebra $\mathfrak{g}:=L^2_{3,w}(ad(E))$. 
We check that the verification of the hypotheses  of the Slice theorem given in \cite[Thm.~3.2 and its Cor.]{FU}  for compact manifolds applies to our particular noncompact setting.  Let $\delta_A$ denote the $L^2$ adjoint to $d_A$. Let $\nabla_A^*$ denote the $L^2$ adjoint of $\nabla_A$. We recall that for $d_A:L^2_{j,w}\to L^2_{j-1,w}$, the adjoint operator is given by 
\begin{align}d_A^{*_{j,w}} = \left(L+\sum_{m=1}^{j}\nabla_A^{\ast m}(1+r^{2m})\nabla^m\right)^{-1}\delta_A\left(L+\sum_{l=0}^{j-1}\nabla_A^{\ast l}(1+r^{2l})\nabla_A^l\right).
\end{align}
 
\begin{proposition}\label{27}Let 
$A\in \mathcal{A}$. There exists a neighborhood of $A$ diffeomorphic to an open subset of $(\Ker(d_A^{*_{3,w}})\cap \mathcal{A})\times \mathcal{G}.$ 
\end{proposition}
\begin{proof}Define 
$\Phi: \Ker(d_A^{*_{3,w}}) \times \mathcal{G}\to \mathcal{A}$ by 
$\Phi(b,s) = s^{-1}(d_A+b)s.$ Then we have 
$D\Phi_{(0,1)}(\dot b,\dot s) = \dot b + d_A\dot s.$ 
To show $\Phi$ is a local diffeomorphism, it suffices to show that $d_A$ is a boundedly invertible isomorphism from $L^2_{3,w}(T^*\TN_k^\nu\otimes ad(\mathcal{E}))$ to the  orthogonal complement of $\Ker(d_A^{*_{3,w}}) $. 
To show surjectivity, it suffices to show that the image of $d_A$ is closed, since the orthogonal complement of the kernel of $d_A^*$ is the closure of the image of $d_A$. To show injectivity and closure of the image, it suffices to bound from below $\frac{\|d_Ab\|^2_{L^2_{2}}}{\| b\|^2_{L^2_{3,w}}}.$ To bound this from below, it suffices to find some $c>0$ such that 
$\| b\|^2_{L^2}\leq c \|(1+r^2)^{\frac{1}{2}}d_A b\|^2_{L^2}.$ By the weighted Hardy's inequality, 
\begin{align}\label{benignwt}
\frac{9}{4}\|V^{-\frac{1}{2}}  b\|^2_{L^2 }\leq  \|rd_A b\|^2_{L^2}.
\end{align}
By the unweighted Hardy's inequality, we have for each choice of origin $o$ and associated radial coordinate $r_o$, 
\begin{align}\label{benign}
\frac{1}{4}\|V^{-\frac{1}{2}} r_o^{-1}b\|^2_{L^2 }\leq  \|d_A b\|^2_{L^2}.
\end{align}   
Summing the estimate \eqref{benign} over origins 
$\nu_\sigma$ (given by the singularities of $V$), 
and combining with \eqref{benignwt} gives the desired estimate. 
\end{proof}
\begin{proposition}
$\mathcal{A}/\mathcal{G}$ is Hausdorff, and the action of $\mathcal{G}$ is free. 
\end{proposition}
\begin{proof} It suffices (see \cite[Sec.~2.9]{Varadarajan})  to show that the set 
$$\Gamma:= \{(A, s^{-1}As): (A,s)\in \mathcal{A}\times \mathcal{G}\}$$is a closed subset of  $\mathcal{A}\times \mathcal{A}.$ \\
So,  consider a  Cauchy sequence 
$\{(d_{A_0}+a_n,  s_n^{-1}(d_{A_0}+a_n)s_n)\}_{n=1}^\infty\subset \Gamma.$ 
Then we have $(a,b)\in L^2_{2,w}(T^*M\otimes ad(\mathcal{E}))\times L^2_{2,w}(T^*M\otimes ad(\mathcal{E}))$ such that $a_n\stackrel{L^2_{2}}{\to} a$ and $s_n^{-1}a_ns_n + s_n^{-1}[d_{A_0},s_n]\stackrel{L^2_{2}}{\to} b.$ The proof that $\{s_n\}_n$ is Cauchy in $L^2_{3,w}$ is now exactly the same as in \cite[Cor., p. 50]{FU}, except that we again use Hardy's inequality to control $\|r^{-1}s_n\|_{L^2}$. 

Unlike the compact case, the freedom of the action of $\mathcal{G}$ requires no additional irreducibility hypotheses in our case. Any element $g$ of the stabilizer of a connection is covariant constant and therefore necessarily the identity, since $\frac{(I-g)}{r}\in L^2$. 
\end{proof}
\begin{coro}
$\mathcal{A}/\mathcal{G}$ is a (Hilbert) manifold. 
\end{coro}
Now we can use the implicit function theorem to give the moduli space of anti-self-dual connections a smooth manifold structure. 
Taking the self-dual component of the curvature, define the smooth map\\ 
$F^+:\mathcal{A}\to L^2_{1,w}(\Lambda^2_+\otimes ad(E)),$ by 
$F^+:A\to F_A^+$. 
\begin{lemma}
$0$ is a regular value of $F^+$.
\end{lemma}
\begin{proof}
Suppose $F^+(A)=0$.  The derivative of $F^+$ at $A$ is $(DF^+)_A = d_A^+.$  Consider $\psi\in  L^2_{1,w}(\Lambda^2_+\otimes ad(E))$ such that  $d_A^+\psi = \delta_A\psi = 0.$ Then we have the Bochner formula 
\begin{align}\label{bff}0 &= \|\nabla_A\psi\|^2_{L^2} + \langle \theta^i\wedge \theta^lR_{ijkl}i_{\Theta_j}i_{\Theta_k}\psi,\psi\rangle_{L^2}
 - \langle  \theta^i\wedge i_{\Theta_j}F_{ij}^+\psi,\psi\rangle_{L^2}\nonumber\\
&= \|\nabla_A\psi\|^2_{L^2}  ,
\end{align}
and $\psi = 0$. 
Here we have used the fact that $R$ is anti-self-dual,  which yields
\begin{align}\langle  \theta^i\wedge\theta^l R_{ijkl}i_{\Theta_j}i_{\Theta_k} \psi,\psi\rangle_{L^2} = \langle  \theta^i\wedge \theta^l R_{ikjl}i_{\Theta_j}i_{\Theta_k}\psi,\psi\rangle_{L^2}. 
\end{align}
Hence, this term is zero by symmetry.  

Suppose now that $\psi$ is $L^2_{1,w}(\Lambda^2_+\otimes ad(E))$ orthogonal to the image of $d_A^+$. Then we have for all smooth compactly supported $b$, 
\begin{align}\label{chalk0}
0 &= \langle \nabla_Ad_A^+b,(1+r^2)\nabla_A\psi\rangle_{L^2 }
 +L\langle d_A^+b,\psi\rangle_{L^2 }.
\end{align}
By elliptic regularity, \eqref{chalk0} implies $\psi$ is smooth. Let $\eta_n$ satisfy $|d\eta_n|\leq \frac{2}{r}$ and $\lim_{n\to\infty }\eta_n(x) = 1$, for every $x$. Manipulating \eqref{chalk0} and then substituting {$b=\eta^2_n\delta_A\psi$} yields 
\begin{align}\label{chalk}
0&= \langle \nabla_Ab,(1+r^2)\delta_A\nabla_A\psi\rangle_{L^2 }+L\langle b,\delta_A\psi\rangle_{L^2 }\nonumber\\
&\quad + \langle (1+r^2)(F_{ma}+R_{ma})b,
i_{\Theta_a}\nabla_{\Theta_m}\psi\rangle_{L^2 }-2\langle r\nabla_Ab,i_{\nabla r}\nabla_A\psi\rangle_{L^2 }\nonumber\\
&= \langle \nabla_Ab,(1+r^2)\nabla_A\delta_A\psi\rangle_{L^2 }+\langle \nabla_{\Theta_m}b,(1+r^2)i_{\Theta_a}(F_{ma}+R_{ma})\psi\rangle_{L^2 }\nonumber\\
&\quad + \langle (1+r^2)(F_{ma}+R_{ma})b,i_{\Theta_a}\nabla_{\Theta_m}\psi\rangle_{L^2 }\nonumber\\
&\quad -2\langle r\nabla_Ab,i_{\nabla r}\nabla_A\psi\rangle_{L^2 }+L\langle b,\delta_A\psi\rangle_{L^2 }\nonumber\\
&= \|(1+r^2)^{\frac{1}{2}}\nabla_A(\eta_n\delta_A\psi)\|^2_{L^2 }+L\|\eta_n\delta_A\psi\|^2_{L^2 }\nonumber\\
&\quad +\langle \nabla_{\Theta_m}(\eta_n^2\delta_A\psi),(1+r^2)i_{\Theta_a}(F_{ma}+R_{ma})\psi\rangle_{L^2 }\nonumber\\
&\quad + \langle (1+r^2)(F_{ma}+R_{ma})(\eta_n^2\delta_A\psi),i_{\Theta_a}\nabla_{\Theta_m}\psi\rangle_{L^2 }\nonumber\\
&\quad -2\langle r\nabla_A(\eta_n^2\delta_A\psi),i_{\nabla r}\nabla_A\psi\rangle_{L^2 }-\|(1+r^2)^{\frac{1}{2}}d\eta_n\otimes \delta_A\psi\|^2_{L^2 }.
\end{align}
Since $\psi\in  L^2_{1,w},$ we can take the limit as $n\to\infty$ in \eqref{chalk} to deduce 
\begin{align}\label{chalk2}
0=&   \|(1+r^2)^{\frac{1}{2}}\nabla_A( \delta_A\psi)\|^2_{L^2 }+L\| \delta_A\psi\|^2_{L^2 }\nonumber\\
&+\langle \nabla_{\Theta_m}( \delta_A\psi),(1+r^2)i_{\Theta_a}(F_{ma}+R_{ma})\psi\rangle_{L^2 }\nonumber\\
&+ \langle (1+r^2)(F_{ma}+R_{ma})( \delta_A\psi),i_{\Theta_a}\nabla_{\Theta_m}\psi\rangle_{L^2 }\nonumber\\
&-2\langle r\nabla_A( \delta_A\psi),i_{\nabla r}\nabla_A\psi\rangle_{L^2 } .
\end{align}
For $L$ sufficiently large (depending on $A$), the sum of the final 3 lines in \eqref{chalk2}  can be absorbed into the first line, 
yielding $\delta_A\psi = 0$, which in turn implies $\psi = 0$ by the Bochner formula \eqref{bff}. Hence the orthogonal complement to the image of $d_A^+$ is zero. The proof that the image of $d_A^+$ is closed is essentially the same as in Proposition \ref{27}. 
Hence $d_A^+$ is surjective.
\end{proof}

\begin{coro}
$(F^+)^{-1}(0)$ is a smooth Hilbert submanifold of $\mathcal{A}.$ 
\end{coro}
Since  $\frac{1}{r}\not\in L^2$,  the data $\vec{\lambda}$ and $\vec{m}$ in \eqref{model} are constant on  $ \mathcal{A}$.
\begin{define}\label{ModSp}
Set $ \mathcal{M}(\vec{\lambda},\vec{m}):= (F^+)^{-1}(0)/\mathcal{G}.$ 
\end{define}
\begin{theo}\label{23andme}
$\mathcal{M}(\vec{\lambda},\vec{m})$ is a smooth manifold. For $A\in (F^+)^{-1}(0)$, there is a neighborhood of  [A] (the orbit of $A$) in the quotient $(F^+)^{-1}(0)/\mathcal{G}$ diffeomorphic to $(B_\epsilon(0)\cap \Ker(d_A^{\ast_{3,w}})\cap [(F^+)^{-1}(0)-d_A])$, for some $\epsilon >0$. The tangent space at $[A]$ is  isomorphic to $ \Ker(d_A^+)\cap \Ker(\delta_A)$.  
\end{theo}
\begin{proof}The proof that the quotient is a manifold and a neighborhood of the orbit of $A$ in $\mathcal{M}(\vec{\lambda},\vec{m})$ has the claimed form is exactly the same as \cite[Thm.~3.16]{FU}, except that we need not carry the metric factor required there. 
The tangent space to $(F^+)^{-1}(0)/\mathcal{G}$ at the orbit of $A$  is naturally isomorphic to $\Ker(d_A^+)/\mathrm{Im}(d_A)$, with $\mathrm{Im}(d_A)$ closed. If we replace the inner product on $L^2_{2,w}\cap \Ker(d_A^+)$ by the equivalent   
$\langle P,Q\rangle_{new} 
=  \langle (T^2+Td_Ar^2\delta_A + d_A\delta_Ar^4d_A\delta_A)P, Q\rangle_{L^2}$, with $T$ a positive constant, then we can identify $\Ker(d_A^+)/\mathrm{Im}(d_A)$ with 
\begin{align}
&\Ker(d_A^+)\cap \Ker(  \delta_A(T^2+Td_Ar^2\delta_A + d_A\delta_Ar^4d_A\delta_A))\nonumber\\
&= \Ker(d_A^+)\cap \Ker(  (T^2+T\delta_Ar^2d_A + \delta_Ad_Ar^4\delta_Ad_A \nonumber\\
&+2T\delta_Ardr\wedge - 4\delta_Ad_Ar^2i_{\nabla_r} rd_A)\delta_A)\nonumber\\
&=  \Ker(d_A^+)\cap \Ker ( \delta_A) ,
\end{align}  for 
$T$ so large that $(T^2+T\delta_Ar^2d_A + \delta_Ad_Ar^4\delta_Ad_A  +2T\delta_Ardr\wedge - 4\delta_Ad_Ar^2i_{\nabla_r} rd_A)$ is positive definite.   
\end{proof}
 We endow $\mathcal{M}(\vec{\lambda},\vec{m})$ with the $L^2$ metric on $ \Ker(d_A^+)\cap \Ker(\delta_A)$.

\section{Isometry}
\label{sec:Isometry}
In this section, we show that the Down transform defines an isometry between two hyperk\"ahler moduli spaces: 
{
the instanton moduli space $\mathcal{M}(\vec{\lambda},\vec{m})$ and the moduli space of the corresponding bow representation defined in \cite[Sec.~3.2]{Second}. 
We recall that the moduli space of the bow representation is the space of gauge equivalence classes of bow solutions \((T,B,Q)\). Its metric is induced by the natural norm \eqref{AffNorm} on the affine space of all bow data (see \cite[Sec.~3]{Second} for more details):
\begin{align}\label{BowNorm}
	\|(\dot{T},\dot{B},\dot{Q} )\|^2 = 
		\sum_{a=0}^3\| \dot T^a\|^2_{L^2([0,l])}+\sum_{m,\sigma}\|\dot B_\sigma^m\|^2+\sum_{\lambda\in \Lambda^0}\|\dot Q_\lambda \|^2
.\end{align}
}

  Since $A_0+a\in \mathcal{A}$,  $a\in L^4$. Consequently Stokes' theorem and the cyclic property of the trace imply that 
$$\int_{\TN_k}\mathrm{tr}F_{A_0+a}\wedge F_{A_0+a}=\int_{\TN_k}\mathrm{tr}F_{A_0}\wedge F_{A_0}.$$ 

 Since   $\vec{\lambda}$ and $\vec{m}$   are constant on $ \mathcal{A}$, $R_0:=\mathrm{ind}(D_{A_0+a})$ is constant for   instantons in $ \mathcal{A}$. 

Let $A(u)$ be a smooth curve  in $\mathcal{M}(\vec{\lambda},\vec{m})$, which we identify with a curve of instanton connections on $\E$ by Theorem \ref{23andme}. Let $A(u,s)$ denote the corresponding family of connections on $\E\otimes \Fe_s.$  Set $\dot A = \frac{d}{du} A(u,s)$. Let $\dot D$  ($=Cl(\dot A))$ denote the corresponding $u-$derivative of the associated family of Dirac operators, $D_{A(u,s)}$, and set $D_s = D_{A(u,s)}$.  
\subsection{Down}
We first compute the variation in the Bow data induced, via the Down transformation, by the variation in $A$. Of course, varying $A$ also varies the image of the Down transform and therefore the derivative of the bow data will include components which are not morphisms between the expected spaces. In computing the variation of the bow data, we will project out these latter summands, effectively replacing derivatives $\frac{d}{du}$ by covariant derivatives $\nabla_{\frac{d}{du}}:= \Pi\frac{d}{du}\Pi$, where $\Pi$ is $\Pi_{D_s}$ or $\Pi_{\D(t,b)}$, depending on context. 
We have 
\begin{align} \frac{d}{du} \Pi_{D_{A(u,s)}} = -D_{s}G_{D_s}\dot D\Pi_{D_s}- \Pi_{D_s}\dot D G_{D_s}D_{s}.
\end{align}

We record relevant derivatives  
\begin{align}\label{Tprime0}
 \frac{d}{du}T^j   =&  \frac{d}{du} \Pi_{D_s}t^j\Pi_{D_s}  \nonumber\\
=&   - \Pi_{D_s}\dot D G_{D_s}\frac{c^j}{\sqrt{V}}\Pi_{D_s} +\Pi_{D_s}\frac{c^j}{\sqrt{V}}G_{D_s}\dot D\Pi_{D_s}\nonumber\\
&- T^j\dot D G_{D_s}D_{s}-D_{s}G_{D_s}\dot D T^j.
\end{align}
\begin{align}\label{Tprime}
\nabla_{\frac{d}{du}}T^j  =& \Pi_{D_s}(\frac{d}{du}T^j)\Pi_{D_s} \nonumber\\
=&       \Pi_{D_s}\dot D I^jG\frac{c^4}{\sqrt{V}}\Pi_{D_s} +\Pi_{D_s} \frac{c^4}{\sqrt{V}}GI^j\dot D\Pi_{D_s}  .
\end{align}
Writing $\nabla_{\frac{d}{ds}} = \frac{d}{ds}+\ii T^0,$ we have 
\begin{align}\label{T0prime}
\nabla_{\frac{d}{du}}T^0  =&  \Pi_{D_s}\dot D G_{D_s} \frac{c^4}{\sqrt{V}} \Pi_{D_s} -\Pi_{D_s} \frac{c^4}{\sqrt{V}}   G_{D_s}\dot D\Pi_{D_s}.
\end{align}
and 
\begin{align}\label{Bprime0}& \frac{d}{du}B_\sigma^j   =  \frac{d}{du} \Pi_{D_{p_\sigma-}}b_\sigma^j\Pi_{D_{p_\sigma+}}  \nonumber\\
&=     \Pi_{D_{p_\sigma-}}\dot D G_{D_{p_\sigma-}}I_a^\dagger\frac{c^4}{\sqrt{V}}\frac{ \ii(I_a^\dagger b_\sigma)^j}{2t_\sigma}\Pi_{D_{p_\sigma+}}  -\Pi_{D_{p_\sigma-}}\frac{ \ii(I_a^\dagger b_\sigma)^j}{2t_\sigma}\frac{c^4}{\sqrt{V}}I_a G_{D_{p_\sigma+}}\dot D\Pi_{D_{p_\sigma+}}\nonumber\\
&\quad - B_\sigma^j \dot D G_{D_{p_\sigma+}}D_{{p_\sigma+}}-D_{p_\sigma-}G_{D_{p_\sigma-}}\dot DB_\sigma^j.
\end{align}
and
\begin{align}\label{Bprime}&\nabla_{\frac{d}{du}}B^j_\sigma   \nonumber\\
&=     \Pi_{D_{p_\sigma-}}\dot D G_{D_{p_\sigma-}}I_a^\dagger\frac{c^4}{\sqrt{V}}\frac{ i(I_a^\dagger b_\sigma)^j}{2t_\sigma}\Pi_{D_{p_\sigma+}}  -\Pi_{D_{p_\sigma-}}\frac{ \ii(I_a^\dagger b_\sigma)^j}{2t_\sigma}\frac{c^4}{\sqrt{V}}I_a G_{D_{p_\sigma+}}\dot D\Pi_{D_{p_\sigma+}}\nonumber\\
&=    \psi_\beta(p_\sigma-)\left\langle  f_j^\dagger\otimes G_{D_{p_\sigma-}}\frac{\ii c^4}{\sqrt{V}}\psi_\alpha(p_\sigma+)\left(\frac{  b_\sigma }{ t_\sigma}\right), \dot D\psi_\beta(p_\sigma-)\right\rangle\langle\cdot, \psi_\alpha(p_\sigma+)\rangle \nonumber\\
&\quad + \psi_\beta(p_\sigma-)\left\langle ((\dot D\psi_\alpha(p_\sigma+))^{c\dagger})^j, G_{D_{p_\sigma+}} \frac{ic^4}{\sqrt{V}}\psi_\beta(p_\sigma-)\left(\frac{    b_\sigma^c }{ t_\sigma}\right)\right\rangle\langle\cdot,\psi_\alpha(p_\sigma+)\rangle.
\end{align}
and
\begin{align}\label{Qprime0}
\frac{d}{du}Q_\lambda   &=  - \ii\frac{d}{du} f_m\otimes D_\lambda f_m^\dagger\otimes    =        - i  f_m\otimes  \dot D  f_m^\dagger\otimes.\\
\label{Qprime}
\nabla_{\frac{d}{du}}Q_\lambda  &=- \ii  f_m\otimes \Pi_{D_{\lambda}}\dot D  f_m^\dagger\otimes.
\end{align}

Similarly, for a curve of bow data, we have (in an obvious notation)
\begin{align} \frac{d}{du} \Pi_{\D_{u,(t,b)}} = -\D_{u,(t,b)}G_{\D_{u,(t,b)}}\dot \D^\dagger\Pi_{\D_{u,(t,b)}}- \Pi_{\D_{u,(t,b)}}\dot \D G_{\D_{u,(t,b)}}\D_{u,(t,b)}^\dagger.
\end{align}
\subsection{Up}
We next compute the variation of the connections induced by the Up transformation and a curve of bow data: 
$$\hat A= \Pi_{\D_{u,(t,b)}}d^H\Pi_{\D_{u,(t,b)}},$$ where $d^H$ denotes the horizontal lift of the exterior derivative to the level set.  

We have 
\begin{align}
\frac{d}{du}\hat A(u)(\Theta_a)  &= \frac{d}{du}\Pi_{\D_{u,(t,b)}}V^{-\frac{1}{2}}W_a\Pi_{\D_{u,(t,b)}}\nonumber\\
&=  \Pi_{\D_{u,(t,b)}}\dot \D G_{\D_{u,(t,b)}}I_a^\dagger[\nabla_{\Theta_4},\D_{u,(t,b)}^\dagger ]\Pi_{\D_{u,(t,b)}}\nonumber\\
&\quad -\Pi_{\D_{u,(t,b)}}[\nabla_{\Theta_4},\D_{u,(t,b)}]I_aG_{\D_{u,(t,b)}}\dot \D^\dagger\Pi_{\D_{u,(t,b)}}\nonumber\\
&\quad -\D_{u,(t,b)}G_{\D_{u,(t,b)}}\dot \D^\dagger A(u)(\Theta_a)
-A(u)(\Theta_a)\dot \D G_{\D_{u,(t,b)}}\D_{u,(t,b)}^\dagger,
\end{align}
and
\begin{align}\nabla_{\frac{d}{du}}\hat A(u)(\Theta_a)  &=  \Pi_{\D_{u,(t,b)}}\dot \D G_{\D_{u,(t,b)}}I_a^\dagger[\nabla_{\Theta_4},\D_{u,(t,b)}^\dagger ]\Pi_{\D_{u,(t,b)}}\nonumber\\
&-\Pi_{\D_{u,(t,b)}}[\nabla_{\Theta_4},\D_{u,(t,b)}]I_aG_{\D_{u,(t,b)}}\dot \D^\dagger\Pi_{\D_{u,(t,b)}}.
\end{align}

\subsection{The Norm Comparison}
In order to show that the Up and Down transformations are isometries, we need to understand better the properties of the two operators 
$G\frac{c^4}{\sqrt{V}}\Pi_{D_s}$ and $G_{\D_{u,(t,b)}}[\nabla_{\Theta_4},\D_{u,(t,b)}^\dagger ]\Pi_{\D_{u,(t,b)}}.$

Let $y\in \text{TN}_k$. Let $\{v_\mu(y)\}_\mu$ be a unitary frame of $\mathcal{E}_y$. Then by Proposition \ref{template}, there exists a   constant $C_D\not =0$, independent of $y$, an equivariant unitary frame $\{v_\mu'(y,b)\}_\mu$ of $\Ker (\D_{(y,b)}^\dagger)$ and therefore a unitary frame 
\begin{align*}
	\left\{v_\mu'(y)=\left(\begin{array}{c} z_\mu(y) 
		\\ w_{\lambda,\mu}(y) 
		\\ n_\mu^-(y)
		\\ n_\mu^+(y)
		\end{array}\right)\right\}_\mu
\end{align*}
of $Down\circ Up(\E)_y$ 
such that 
\begin{align}
\left(\begin{array}{c} \Pi_{D_s}\frac{\ii c^4}{\sqrt{V}}G_{D_s}  v_\mu(y)\delta_y
\\
\frac{1}{2\ii}\langle v_\mu(y),\phi(y)\rangle\phi
\\\Pi_{D_{p_\sigma +}}\frac{c^4}{\sqrt{V}}\frac{b_\sigma^{\dagger}}{t_\sigma}G_{\nabla^*\nabla}  v_\mu(y)\delta_y
\\\Pi_{D_{p_\sigma -}}\frac{c^4}{\sqrt{V}}\frac{b_\sigma^{c\dagger}}{t_\sigma}G_{\nabla^*\nabla}  v_\mu(y)\delta_y
\end{array}\right) = C_Dv_\mu'.
\end{align}

Rewrite this relation as 
\begin{align}\left(\begin{array}{c}  G_{D_s}\frac{\ii c^4}{\sqrt{V}}\psi_\alpha(s,)(y) 
\\\langle v_\mu(y),\phi(y)\rangle\phi
\\   (G_{\nabla^*\nabla}(\frac{c^4}{\sqrt{V}}\psi_\alpha(p_\sigma+,))(\frac{b_\sigma }{t_\sigma}))(y)  
\\  (G_{\nabla^*\nabla}(\frac{c^4}{\sqrt{V}}\psi_\alpha(p_\sigma-,))(\frac{b_\sigma^{c }}{t_\sigma}))(y)       
\end{array}\right) 
= C_D\left(\begin{array}{c} \langle f_j\otimes \psi_\alpha(s,),z_\mu\rangle f_j^\dagger \otimes v_\mu(y)
\\ 2 \ii w_{\lambda,\mu} 
\\ -2\langle \psi_\alpha(p_\sigma+,),n_\mu^-(y)\rangle v_\mu(y)
\\ -2\langle \psi_\alpha(p_\sigma-,),n_\mu^+(y)\rangle v_\mu(y)
\end{array}\right).
\end{align}

Similarly, from Corollary \ref{coro25}, we have  $C_U \not = 0$, independent of $s$,  an equivariant unitary basis $\{u_\alpha(s)\}_\alpha$ of $\mathpzc{E}_s$ and a unitary basis   $\{\psi_\alpha'\}_\alpha$ of $\Ker_{L^2}(D_s)$ such that 
\begin{align}  (G_{\D_{(y,b)}}[\nabla_{\Theta_4},\D_{(y,b)}^\dagger ]v_\mu')(s)
  = C_U\langle c^4f_j^\dagger\otimes v_\mu'(y,\cdot),\psi_\alpha'(s,y)\rangle f_j\otimes  u_\alpha(s),
\end{align}
for any unitary basis $\{v_\mu'\}_\mu$ of $\Ker(\D_{(y,b)})$.
With these preliminaries, we now show 
\begin{theo}\label{IsomThm}
The Up and Down transforms are isometries. 
\end{theo}
\begin{proof}
First consider on the Upside:
\begin{align}\label{upside}&\| \nabla_{\frac{d}{du}}\hat  A\|_{L^2}^2  = \int_{\TN_k} \nabla_{\frac{d}{du}}\hat  A(\Theta_a)_{\mu\nu}(y)\overline{\nabla_{\frac{d}{du}}\hat  A}(\Theta_a)_{\mu\nu}(y)dv\nonumber\\
&= 2\int_{\TN_k} \langle  \nabla_{\frac{d}{du}}\hat  A(\Theta_a)v_\mu(y), \dot \D G_{\D_{u,(t,b)}}\e_a^\dagger[\nabla_{\Theta_4},\D_{u,(t,b)}^\dagger ] v_\mu(y)\rangle dv\nonumber\\
&= 2\int_{\TN_k} \langle  \nabla_{\frac{d}{du}}\hat  A(\Theta_a)v_\mu(y), \dot \D   C_U\langle c^4f_j^\dagger\otimes v_\mu ,\psi_\alpha'\rangle \e_a^\dagger f_j\otimes  u_\alpha\rangle dv\nonumber\\
&= -2C_U\int_{\TN_k} \langle  v_\mu(y), \dot \D      f_j\otimes  u_\alpha\rangle \langle \psi_\alpha', \dot D f_j^\dagger\otimes  v_\mu  \rangle dv.
\end{align}
On the Down side, we have {the norm \eqref{BowNorm} on the bow moduli space
}
\begin{align}\label{downside}
&\|\dot\D\|^2  = \sum_{a=0}^3\| \dot T^a\|^2_{L^2([0,l])}+\sum_{m,\sigma}\|\dot B_\sigma^m\|^2+\sum_{\lambda\in \Lambda^0}\|\dot Q_\lambda \|^2\nonumber\\
&= \sum_{a=0}^32\langle   \Pi_{D_s}\dot D I_a  G_{D_s}\frac{c^4}{\sqrt{V}}\Pi_{D_s},\dot T^a \rangle_{L^2([0,l])}
 \nonumber\\
&+\sum_{m,\sigma} \langle  f_m^\dagger\otimes G_{D_{p_\sigma-}}
\left(\frac{\ii c^4}{\sqrt{V}}\psi_\alpha (p_\sigma+)\right)
\frac{  b_\sigma }{ t_\sigma}, \dot D\psi_\beta(p_\sigma-)\rangle\langle  \psi_\beta(p_\sigma-)   ,\dot B_\sigma^m \psi_\alpha(p_\sigma+)\rangle\nonumber\\
&+  \langle (\dot D\psi_\alpha(p_\sigma+))^{c\dagger} , f_m^\dagger\otimes G_{D_{p_\sigma+}} \left(\frac{\ii c^4}{\sqrt{V}}\psi_\beta (p_\sigma-)\right)\frac{    b_\sigma^c }{ t_\sigma}\rangle\langle\psi_\beta(p_\sigma-) ,\dot B_\sigma^m \psi_\alpha(p_\sigma+)\rangle\nonumber\\
&+\sum_{\lambda\in \Lambda^0}\langle - \ii  f_m\otimes \Pi_{D_{\lambda}}\dot D  f_m^\dagger\otimes,\dot Q_\lambda \rangle\nonumber\\
&= -2C_D\int\langle   \dot D  f_k^\dagger \otimes v_\mu(y),\psi_\beta \rangle \langle  \ii \e_a^\dagger\dot T^a f_k\otimes  \psi_\beta(s,),z_\mu\rangle ds dy
 \nonumber\\
&\quad+\sum_{m,\sigma}2C_D\langle \dot D f_m^\dagger\otimes  v_\mu(y) , \psi_\beta(p_\sigma-)\rangle  \langle -\dot B_\sigma^\dagger f_m\otimes \psi_\beta(p_\sigma-),n_\mu^-(y)\rangle\nonumber\\
&\quad + 2C_D\langle  \psi_\alpha(p_\sigma+)  , \dot D f_m^\dagger\otimes v_\mu(y)\rangle \langle  n_\mu^+(y),\dot B^{c \dagger} f_m\otimes \psi_\alpha(p_\sigma+)\rangle\nonumber\\
&\quad-\sum_{\lambda\in \Lambda^0}\langle \dot D  f_m^\dagger\otimes \phi,\psi_\alpha(\lambda,)\rangle\langle  \dot Q_\lambda^\dagger f_m\otimes \psi_\alpha(\lambda,), \phi \rangle\nonumber\\ 
&= 2C_D\mathfrak{Re}\,\biggl[\int\langle   \dot D  f_k^\dagger \otimes v_\mu(y),\psi_\beta(s,\cdot) \rangle \langle (-\ii\dot T^0+\ii \e_j\dot T^j) f_k\otimes  \psi_\beta(s,),z_\mu\rangle ds dy
 \nonumber\\
&\quad+\sum_{m,\sigma} \langle \dot D f_m^\dagger\otimes  v_\mu(y) , \psi_\beta(p_\sigma-)s\rangle  \langle -\dot B_\sigma^\dagger f_m\otimes \psi_\beta(p_\sigma-),n_\mu^-(y)\rangle\nonumber\\
&\quad +  \langle   \dot D f_m^\dagger\otimes v_\mu(y),\psi_\alpha(p_\sigma+)\rangle \langle \dot B^{c \dagger} f_m\otimes \psi_\alpha(p_\sigma+),n_\mu^+(y)\rangle\nonumber\\
&\quad-  \sum_{\lambda\in \Lambda^0}\langle \dot D  f_m^\dagger\otimes v_\mu ,\psi_\alpha(\lambda,)\rangle\langle   \dot Q_\lambda^\dagger f_m\otimes \psi_\alpha(\lambda,), w_{\lambda,\mu}\rangle\biggr].
\end{align}

Comparing \eqref{downside} to \eqref{upside}, we see that  the moduli spaces are isometric up to scale. 
\end{proof}
  
\appendix
\section*{Appendices}
\addcontentsline{toc}{section}{Appendices}
\renewcommand{\thesubsection}{\Alph{subsection}}
\subsection{Elementary quaternion conventions and identities} \label{ApA} 
We list standard quaternion identities here for easy reference.  Let $I_1,I_2,I_3$ denote a standard basis of imaginary quaternions, with $I_j^2 = -1$, $I_1I_2=I_3$, and cyclic permutations. Set $I_4 = 1$. 
We will use $a,b,c$ to denote basis indices running from $1$ to $4$, and $i,j,k,m$ to denote basis indices running from 1 to 3. Let $\ii=\sqrt{-1}.$  We consider the standard action of the quaternions on $S=\IC^2$  and the contragredient action on $S^*$. Hence for $\phi\in S^*$, we have  $I_a\phi = \phi\circ I_a^\dagger$. We will frequently identify $S$ with $\Hom(\IC,S)$; thus for $z\in S$, $z^\dagger$ denotes the metrically defined dual element in $S^*$.   Fixing a unitary basis of $\IC^2$, we will also use {\em charge conjugation} \cite[Sec.3.1]{First}: 
$$\left(\begin{array}{c}z_1\\z_2\end{array}\right)^c : = \left(\begin{array}{c}-\bar z_2\\\bar z_1\end{array}\right).$$ 
Charge conjugation is basis dependent, however it only enters naturally into computations (in basis-independent pairs). It provides a simple way to encode certain quaternion identities.  We extend these notations in an obvious way to $S\otimes \End(Z)$, for any Hermitian vector space $Z$. With these conventions, we have the following useful elementary identities. 

For $z,w\in \IC^2$, 
\begin{align}\label{chargeconj0}
I_az\otimes w^{\dagger}I_a^\dagger &= 2w^\dagger(z)1_{\IC^2},
\end{align}
which is equivalent to 
\begin{align}\label{chargeconj2}
I_a\left(\begin{array}{cc}a&b\\c&d\end{array}\right)I_a^\dagger &= 2(a+d)1_{\IC^2},
\end{align}
and to 
\begin{align}\label{fierz}
I_aI_bI_a^\dagger &= 4\delta_{b4}I_4.
\end{align}
We also have 
\begin{align}\label{chargeconj1}
I_az\otimes w^{\dagger}I_a &= -2w^{c}\otimes z^{c\dagger}.
\end{align}
For any $X=I_a X^a$ we have $I_a I_b^\dagger\otimes (I_a^\dagger X I_b) = 4(1\otimes X_4 + I_k \otimes X_k)$, since
\begin{align}
I_a I_b^\dagger\otimes (I_a^\dagger X I_b) 
&= \overset{a=b}{1\otimes(I_a^\dagger X I_a)} 
+ \overset{a=0}{I_k^\dagger\otimes(X I_k)} 
+\overset{b=0}{I_k\otimes(I_k^\dagger X)}
 + \overset{a\neq 0\neq b}{\sum_{i\neq j} I_i I_j^\dagger\otimes(I_i^\dagger X I_j)}
 \nonumber\\
 &=I_4\otimes (4 X_4 1) - \sum_{\mathclap{(i,j,k)\circlearrowleft(1,2,3)}} I_k\otimes (X I_k + I_k X - I_i X I_j + I_j X I_i)
  \nonumber\\
 &= 4(I_4\otimes X_4 + I_k \otimes X_k),
 \label{eq:tens}
\end{align}
where we have used \eqref{fierz} and $X I_k + I_k X - I_i X I_j + I_j X I_i = -4 X_k 1_{\IC^2}$ for any $(i,j,k)$ cyclic permutation of $(1,2,3).$
\subsection{Covariant Derivatives of Bifundamentals} 
\label{sec:CovDer}
The bifundamental datum $b_\sigma\in \Hom(\Fe_{p_\sigma+},\calS\otimes \Fe_{p_\sigma-})$ defines a canonical equivariant section of $\mathcal{P}^*\Hom(\Fe_{p_\sigma+},\calS\otimes \Fe_{p_\sigma-})$ on the level set $\bm{\mu}^{-1}(\ii\nu)$ of the moment map and therefore descend to canonical section of $(S^+)^\ast\otimes \Hom(\Fe_{p_\sigma+}, \Fe_{p_\sigma-})$ on $\TN_k$. In this subsection, we use the moment map equations to compute covariant derivatives of these sections and related identities. We recall that the covariant derivative in the direction $X$ of these equivariant sections is computed by taking the derivative of the corresponding equivariant section in the direction of the horizontal lift $X^H$ of $X$. (See \cite[Intro. to Sec.~5]{Second}.) 
 
We view each $\nu_\sigma$ of the moment map level \eqref{nanosec}  
and the coordinate {$t=t_1\e_1+t_2\e_2+t_3\e_3$}  on the base of $\TN_k^\nu$ as  imaginary quaternions. For convenience, we set $t_\sigma:=|t-\nu_\sigma|$ and $\tslash_\sigma:=\ii(t-\nu_\sigma).$ Then, the small bow representation moment map conditions are 
\begin{align}\label{ofsmallmoment}b_\sigma b_\sigma^\dagger=t_\sigma+\tslash_\sigma, \text{ and } b_\sigma^c (b_\sigma^c)^\dagger=t_\sigma-\tslash_\sigma.
\end{align}
 
Using \eqref{ofsmallmoment} we express $\pi_k^*(dt^2)$ as 
\begin{align*}
2d\vec{t}_\sigma^{\, 2}&=\mathrm{tr}\, d \tslash_\sigma^2=\mathrm{tr}(d(b_\sigma b_\sigma^\dagger-t_\sigma))^2
=\mathrm{tr}(db_\sigma b_\sigma^\dagger+b_\sigma db_\sigma^\dagger)^2-2 dt_\sigma d\mathrm{tr}b_\sigma b_\sigma^\dagger+2d t_\sigma^2\\
&=(b_\sigma^\dagger db_\sigma)^2+(d b_\sigma^\dagger b_\sigma)^2
+2 db_\sigma^\dagger db_\sigma (b_\sigma^\dagger b_\sigma)
-4 dt_\sigma^2+2dt_\sigma^2\\
&=4 t_\sigma db_\sigma^\dagger db_\sigma
+(b_\sigma^\dagger d b_\sigma)^2+(d b_\sigma^\dagger b_\sigma)^2-2 \left(d (b_\sigma^\dagger b_\sigma/2)\right)^2\\
&=4t_\sigma d b_\sigma^\dagger d b_\sigma
+\frac{1}{2}(b_\sigma^\dagger d b_\sigma-d b_\sigma^\dagger b_\sigma)^2.
\end{align*}
In particular, we have 
\begin{align*}
db_\sigma^\dagger d b_\sigma=\frac{d\vec{t}{}^{\,2}}{2t_\sigma}
-\frac{(b_\sigma^\dagger d b_\sigma-d b_\sigma^\dagger b_\sigma)^2}{8t_\sigma}=\frac{d\vec{t}{}^{\,2}}{2t_\sigma}+\frac{\hat{\eta}_\sigma^2}{\frac{1}{2t_\sigma}},
\end{align*}
with 
{
$\hat{\eta}_\sigma:=\ii \frac{b_\sigma^\dagger d b_\sigma- d b_\sigma^\dagger b_\sigma}{4 t_\sigma}.$
}

Rewriting \eqref{ofsmallmoment} as 
$$b_\sigma b^\dagger_\sigma -\frac{1}{2}b^\dagger_\sigma b_\sigma 1_{2\times 2} = \ii\tslash_\sigma.$$
we see 
	\begin{align}&{(db_\sigma)} b^\dagger_\sigma +b_\sigma db^\dagger_\sigma -\frac{1}{2}d(b^\dagger_\sigma b_\sigma) 1_{2\times 2} = \ii d\tslash_\sigma\nonumber\\
&\Rightarrow 
{(d b_\sigma)} |b_\sigma|^2 +b_\sigma d(b^\dagger_\sigma)b_\sigma -\frac{1}{2}b_\sigma d(b^\dagger_\sigma b_\sigma)   = d\tslash_\sigma b_\sigma\nonumber\\
&\Rightarrow 
 db_\sigma + b_\sigma\frac{1}{4|t_\sigma|}(d(b^\dagger_\sigma)b_\sigma - b^\dagger_\sigma db_\sigma)    = \frac{\ii d\tslash_\sigma}{2|t_\sigma|} b_\sigma.
\end{align}

Hence 
$0=db_\sigma{-}b_\sigma\ii\hat{\eta}_\sigma  -\frac12\frac{d\tslash}{t_\sigma}b_\sigma,$ 
{so $(d-\ii\hat{\eta}_\sigma)b_\sigma = \frac{1}{2} \frac{d \tslash}{t_\sigma}b_\sigma$.} 
Similarly, for the charge conjugate $(d+\ii\hat{\eta}_\sigma)b_\sigma^c=-\frac12\frac{d\tslash}{t_\sigma}b_\sigma^c.$
Recall from \cite[Eq.~(54)]{Second} that the horizontal lifts of our orthonormal frame $\{\Theta_a\}_a $ are given by 
\begin{align*}
\Theta_j^H=V^{-\frac{1}{2}}W_j&= V^{-\frac{1}{2}}\left(\frac{\p}{\p t^j}-\sum_{\sigma=1}^k\eta_\sigma^j\frac{\p}{\p \phi_\sigma}\right),& &\text{for } j=1, 2, 3
,\end{align*} 
and 
$$\Theta_4^H=V^{-\frac{1}{2}}W_0= V^{-\frac{1}{2}}\left(\sum_{\sigma=1}^k\frac{1}{2r_\sigma}\frac{\p}{\p \phi^\sigma}-\sum_{\sigma=1}^{k+1} \frac{\p}{\p t_0^\sigma}\right).$$  
Consequently,  
{
	\begin{align}\label{ssa}\Theta_a^Hb_\sigma    &= -\frac{\ii \e_a^\dagger}{2\sqrt{V}t_\sigma}b_\sigma
		= - \frac{\ii}{2\sqrt{V} }\frac{b_\sigma}{t_\sigma} I_a,\\
	\Theta_a^Hb_\sigma^c   & =  \frac{\ii \e_a^\dagger}{2\sqrt{V}t_\sigma}b_\sigma^c
	=\frac{\ii}{2\sqrt{V}}	\frac{b_\sigma^c}{t_\sigma} I_a
.\end{align}
}
 From this, we see that $\frac{b_\sigma}{t_\sigma}$ satisfies a Dirac type equation.
\begin{align}
\e_a^\dagger \nabla_{\Theta_a} \frac{b_\sigma}{t_\sigma}  &=
\e_a^\dagger\Theta_a^H\frac{b_\sigma}{t_\sigma}  
  = \frac{\ii \e_j^\dagger \e_j}{2\sqrt{V}t_\sigma^2}b_\sigma-\frac{ \e_j^\dagger(t^j-\nu^j_\sigma)b_\sigma}{\sqrt{V}t_\sigma^3}
-\frac{\ii}{2\sqrt{V}t_\sigma^2}b_\sigma\nonumber\\
&= \frac{\ii }{\sqrt{V}t_\sigma^2}b_\sigma-\frac{ \ii  (b_\sigma b_\sigma^\dagger-t_\sigma)b_\sigma}{\sqrt{V}t_\sigma^3}
= 0.\end{align}
Taking adjoints, we also have 
{
	\begin{align}\label{diracadj}
 I_a^\dagger \nabla_{\Theta_a}\frac{b_\sigma^\dagger}{t_\sigma} 	=\nabla_{\Theta_a}\frac{b_\sigma^\dagger}{t_\sigma} \e_a  =0.
\end{align} 
}

Equation \eqref{ssa} also implies that 
  $b$ is harmonic since 
\begin{align}\nabla^*\nabla b&=(- \nabla_{\Theta_4}^2-\frac{1}{\sqrt{V}}\nabla_{\Theta_j}\sqrt{V}\nabla_{\Theta_j})b\nonumber\\
&=-\frac{-\ii}{\sqrt{V}}\frac{1}{2t}\frac{-\ii}{\sqrt{V}}\frac{b}{2t}-\frac{1}{\sqrt{V}}\nabla_{\Theta_j}\ii \e_j\frac{b}{2t}=0.
\end{align}  
{Using this fact and \eqref{diracadj} we have
\begin{align}\label{Eq:commbG}
	&[G,b] = G\left[\nabla_4^2 + \frac{1}{\sqrt{V } }\nabla_j\sqrt{V} \nabla_j, b\right] G \nonumber
      \\&= 2 G \left((\nabla_4b)\nabla_4 + (\sqrt{V} \nabla_j b) \frac{1}{\sqrt{V} }\nabla_j\right) G 
	= - G \frac{i}{\sqrt{V} } \frac{b}{t} I_a\nabla_a G
.\end{align}}

\printbibliography

\end{document}